\documentclass[10pt,twoside]{article}
\usepackage{lmodern}
\usepackage{cite}
\usepackage[T1]{fontenc}
\usepackage{amsthm}
\usepackage{amsfonts}
\usepackage{amssymb}
\usepackage{amsmath}
\usepackage{graphicx}
\usepackage{mathrsfs} % for \mathscr{F}
\usepackage{mathtools}
\usepackage{hyperref}
\usepackage{array,multirow}
\usepackage{caption}   % for tables
\usepackage{float}   % for tables
\usepackage{microtype}
\usepackage{subfiles}    %for splitting input file

\usepackage{tikz}
\usetikzlibrary{matrix}
\usetikzlibrary{calc} 

   \author{by  Martin~Seysen\\
     }
   \title{A computer-friendly construction of the monster}
   
   \newcommand{\Field}{\mathbb{F}}
   \newcommand{\setZ}{\mathbb{Z}}
   \newcommand{\setR}{\mathbb{R}}
   \newcommand{\setN}{\mathbb{N}}
   \newcommand{\setC}{\mathbb{C}}
   
   \newcommand{\setF}{\mathbb{F}}

   \newcommand{\setM}{\mathbb{M}}

   \newcommand{\Cfrac}[1]{{/ \kern-0.30em / #1 / \kern-0.30em /}}

   \newcommand{\Skip}[1]{}

   \newcommand {\End}{\,{\mathop{\rm E\kern-0.05em{n}\kern-0.05em{d}}}}
   \newcommand {\SL}{\,{\mathop{\rm SL}}}

  \newcommand {\im}{\,{\mathop{\rm im}}\,}

  \newcommand {\Aut}{\,{\mathop{\rm Aut}}\,}
  \newcommand {\Out}{\,{\mathop{\rm Out}}\,}

   \newcommand {\tr}{\,{\mbox{t{r}}}}

   \newcommand{\conj}[1]{{\mkern 1.5mu\overline{\mkern-1.5mu#1\mkern-1.5mu}\mkern 1.5mu}}

   \newcommand{\AutStP}{\Aut_{\!\mbox{\scriptsize St}} \mathcal{P}}
   \newcommand{\Cyclicxyz}{(\& \! \stackrel{\curvearrowright}{xyz})}

   \newcommand {\biscalar}[1]{\left<\! \left< #1 \right> \! \right>}
   \newcommand {\smallbiscalar}[1]{
       \langle\! \langle #1 \rangle \! \rangle}

\DeclareMathSymbol{\mlq}{\mathord}{operators}{``}
\DeclareMathSymbol{\mrq}{\mathord}{operators}{`'}

   \newcommand {\fatline}[1]{\noindent {\bf #1}}

   \newcommand {\proofend}{\noindent $\Box$ }

   \newtheorem {Theorem}{Theorem}[section]
   \newtheorem {Definition}[Theorem]{Definition}
   \newtheorem {Corollary}[Theorem]{Corollary}
   
   \newtheorem {Lemma}[Theorem]{Lemma}

   \numberwithin{equation}{Theorem}
 
   \newcounter{StatementListCounter}
   
    \newenvironment{StatementList}[1]{\begin{list} %
     {(#1.\arabic{StatementListCounter})}{\usecounter{StatementListCounter}}
    }{\end{list}}

\begin{document} {\large}

  \maketitle

\centerline{\large Preliminary version}

  \begin{abstract}
Let  $\mathbb{M}$ be the  monster group which is the largest 
sporadic finite simple  group, and  has first been constructed in 1982 by
Griess. In 1985, Conway has constructed a  196884-dimensional 
representation $\rho$ of $\mathbb{M}$ with matrix coefficients in 
$\mathbb{Z}[\frac{1}{2}]$.
So these matrices may be reduced modulo any (not necessarily prime) odd
number $p$, leading to representations of $\mathbb{M}$ in odd characteristic.
The representation $\rho$  is based on representations of two maximal 
subgroups $G_{x0}$ and $N_0$ of $\mathbb{M}$. In ATLAS notation, $G_{x0}$ 
has structure $2_+^{1+24}.\mbox{Co}_1$ and $N_0$ has structure 
$2^{2+11+22}.( M_{24} \times S_3)$. Conway has 
constructed an explicit set of generators of $N_0$, but not of $G_{x0}$.

This paper is essentially a rewrite of Conway's construction augmented
by an explicit construction of an element of $G_{x0} \setminus N_0$.
This gives us a complete set of generators of $\mathbb{M}$.
It turns out that the matrices of all generators of $\mathbb{M}$ 
consist of monomial blocks, and of blocks which are essentially
Hadamard matrices scaled by a negative power of two. Multiplication
with such a generator can be programmed very efficiently if the
modulus $p$ is of shape $2^k-1$. 

So this paper may be considered a as programmer's reference for 
Conway's construction of the monster group $\mathbb{M}$.
We have implemented representations
of $\mathbb{M}$ modulo 3, 7, 15, 31, 127, and 255.
  \end{abstract}

\section*{}

{\noindent \bf Key Words:}

Monster group, finite simple groups, group representation

\vspace{1ex}

\noindent MSC2020-Mathematics Subject Classification (2020):
20C34, 20D08, 20C11

%%%%%%%%%%%%%%%%%%%%%%%%%%%%%%%%%%%%%%%%%%%%%%%%%%%%%%%%%%%%%%%%%%%%%%%%%%%%%%

\section{Introduction}
\label{section:Introduction}

%%%%%%%%%%%%%%%%%%%%%%%%%%%%%%%%%%%%%%%%%%%%%%%%%%%%%%%%%%%%%%%%%%%%%%%%%%%%%%

Let $\setM$ be the  monster group, which is the largest 
sporadic finite simple  group, and  has first been constructed
by Griess \cite{Griess:Friendly:Giant}. That construction has
been greatly simplified by Conway \cite{Conway:Construct:Monster},
leading to a rational representation $\rho$ of $\setM$ of dimension
196883+1.
This paper is essentially a rewrite of Conway's construction of
$\setM$, augmented by an explicit construction of the representations 
of a complete set of generators of $\setM$. Here  all generators of
$\setM$ have coefficients in $\setZ[\frac{1}{2}]$ as 
in~\cite{Conway:Construct:Monster}, and they can be computed
very efficiently. So a programmer may use this paper as a reference
for implementing the representation $\rho$ of  $\setM$  modulo a small
odd number $p$.
 
The success of the construction in~\cite{Conway:Construct:Monster} 
is due to the fact that the representation 
$\rho(N_{x0})$ of a large subgroup $N_{x0}$ of $\setM$ of structure
$2^{1+24+11}.M_{24}$ can be made explicit and  monomial.
 We use the notation in the ATLAS
\cite{Atlas} for describing the structure of a group, so $M_{24}$ is 
the Mathieu group acting as a permutation group  on 24 elements.
There is a maximal subgroup  $N_0$ of $\setM$ with $N_0 : N_{x0} = 3$.
 For a so-called
triality element $\tau \in N_0 \setminus N_{x0} $ the 
(non-monomial) representation $\rho(\tau)$ is given explicitly
in~\cite{Conway:Construct:Monster}, so that we can 
effectively compute in $\rho(N_{0}) $.

There is another maximal subgroup $G_{x0}$ of $\setM$ of structure
$2^{1+24}.\mbox{Co}_1$. Here $2^{1+24}$ is an extraspecial
2-group, and the simple group $\mbox{Co}_1$ is the automorphism
group of the 24-dimensional Leech lattice $\Lambda$ modulo~2,
see e.g. \cite{Conway-SPLG}.

The group $N_{x0}$ is a maximal subgroup of  $G_{x0}$. 
So an explicit description of 
$\rho(\xi)$ for any $\xi \in G_{x0} \setminus N_{x0}$ 
gives us the capability to compute in 
$\rho(\setM)$.  Holmes and Wilson  \cite{holmes_wilson_2003} have
used similar ideas to construct the representation of the Monster 
in characteristic~3. Their construction uses non-monomial matrices
with blocks of size up to $276 \times 276$.
In~\cite{Conway:Construct:Monster} no explicit
construction of $\rho(\xi)$ is given for any 
$\xi \in \setM \setminus N_0$. 

The essential new result in this paper is a computer-free explicit 
construction of $\rho(\xi)$ for such a generator $\xi$ in
section~\ref{Sect:construct:xi} ff.

The construction \cite{Conway:Construct:Monster} of $\setM$ uses
the Parker loop, which is a non-associative loop of order
$2^{1+12}$.
For computations in the Parker loop a cocycle can be used, see
e.g. Aschbacher \cite{Aschbacher-Sporadic}, Chapter~4. Such a
cocycle is not unique. 
While the choice of that cocycle is not too  important for 
computing in the subgroup $N_0$ of $\setM$,
we have to select a specific  cocycle for the Parker loop in
such a way that the (non-monomial) representation of a certain 
element $\xi$  of $G_{x0} \setminus N_{x0}$ obtains a reasonably 
simple structure. 

In our construction we take the notation from
\cite{Conway:Construct:Monster}, except for two explicitly
stated sign changes required for the construction of a
generator $\xi \in G_{x0} \setminus N_{x0}$.
One of these changes is explained in section~\ref{section:N0}
and motivated by Ivanov's construction of  $\setM$ in
\cite{citeulike:Monster:Majorana}. This leads to simpler
relations in the group $N_0$.

The other sign change is explained in section~\ref{subsection:rep:4096}).
Therefore we remark that
another group $N(4096_x)$ of structure $2^{1+24}.\mbox{Co}_1$, but not
isomorphic to $G_{x0}$, plays an important role in Conway's construction.
$N(4096_x)$ is a subgroup of the Clifford group $\mathcal{C}_{12}$ of 
structure  $2_+^{1+24}.\mbox{O}^+_{24}(2)$. 
Nebe, Rains, and Sloane~\cite{NebRaiSlo2001} have constructed the group  
$\mathcal{C}_{12}$  as the automorphism group of $M_1^{\otimes 12}$,  
where $M_1$ is a certain 2-dimensional $\setZ[\sqrt{2}]$-lattice. 
This construction leads to a 4096-dimensional real representation of 
$\mathcal{C}_{12}$ and also of its subgroup $N(4096_x)$, which is 
called  $4096_x$ in \cite{Conway:Construct:Monster}.
Selecting two orthogonal short rational  vectors in $M_1$ as a basis 
of the vector space   $M_1 \otimes \setR$, we also obtain a basis of the 
real vector space $4096_x$ = $(M_1 \otimes \setR)^{\otimes 12}$
in the usual way.
We change some signs of the basis vectors of $4096_x$ used in
\cite{Conway:Construct:Monster}, so that they are compatible with the
signs of the basis vectors of $4096_x$, when constructed as a tensor 
product.
We will not show this compatibility, since we do not need it
explicitly in this paper. We claim that
this compatibility greatly simplifies the construction of a
$\xi \in G_{x0} \setminus N_{x0}$.

It is worth noting that the complex analogue $\mathcal{X}_{n}$ of the
real Clifford group $\mathcal{C}_{n}$, also defined 
in~\cite{NebRaiSlo2001}, plays an important role in the theory of 
quantum computing. A general quantum  circuit with $n$ qubits is modelled 
using a dense subgroup of the complex unitary group $U(2^n,\setC)$ and 
cannot be simulated in polynomial time on a classical computer. 
There is an important  class of quantum  circuits with $n$ 
qubits called {\em stabilizer circuits},
which can be modelled using the discrete subgroup $\mathcal{X}_{n}$
of  $U(2^n,\setC)$. Here all quantum gates can be modelled 
by tensor products of certain well-behaved monomial and of Hadamard 
matrices. The Gottesman-Knill theorem states that 
stabilizer circuits can be simulated in  polynomial time, see e.g.
 \cite{AaronsonGottesman2004} for background.  
  This `explains' in a way,  why 
choosing the signs of the basis vectors of  $4096_x$ to be compatible  
with~\cite{NebRaiSlo2001} and with  the simulation in 
~\cite{AaronsonGottesman2004} may simplify computations in the subgroup
$2^{1+24}.\mbox{Co}_1$ of~$\mathcal{X}_{12}$.

The hard part of Conway's construction \cite{Conway:Construct:Monster}
was to find an algebra invariant under $\rho(\setM)$, similar to 
the Griess algebra  defined in \cite{Griess:Friendly:Giant}. Therefore
Conway defined an algebra visibly invariant under $\rho(G_0)$ and  
showed that this algebra is also invariant under $\rho(N_{0})$, using
a basis of $\rho$ where $N_{0}$ has a simple representation. In this
paper we need an explicit representation $\rho(\xi)$ of a
$\xi \in G_0 \setminus N_{x0}$. So we adjust the signs of that basis 
in order to simplify $\rho(G_0)$. 

Define a Hadamard step on $\rho(\setM)$ to be
a multiplication with a  matrix that contains just
monomial blocks and blocks of $2 \times 2$-Hadamard 
matrices, i.e. matrices of shape 
$c \begin{psmallmatrix} 1 & 1\\ 1 & -1\end{psmallmatrix}$
for $c \in \{\frac{1}{2}, 1\}$.
We will see that the representations of all our generators of
$\setM$ can be decomposed into at most 6 Hadamard steps plus some
monomial operations.

For any odd natural number $p$ let $\rho_p$ be the representation
$\rho$ of the monster, where all coefficients are taken modulo $p$. 
Computations in $\rho_p$ are very efficient if $p+1$ is a reasonably
small power $2^k$  of two. In this case we may represent an integer
modulo $p$ with $k$ bits, putting $(1,...,1)_2$ = $(0,...,0)_2$ = $0$.
Then negation can be done by complementing
all bits, halving can be done by right rotation, and the carry bit 
of an addition has the same valence as the least significant bit.
Using integer additions, and bitwise and shift operations on a 32-bit 
or 64-bit computer, several components of a vector in $\rho_p$ can be
negated, halved, or added  modulo a small number $p = 2^k-1$ 
simultaneously in a single register.

We have implemented the representations $\rho_p$ of the monster for
$p = 3, 7, 15, 31, 127$, and $255$, see \cite{mmgroup2020}.  
Calculating in the monster modulo different numbers is useful e.g. for
distinguishing between classes in $\setM$ of the same order, see
\cite{Barraclough:2005:CCR}.

Let $g \in \setM$ be the product of an arbitrary element of $G_{x0}$
with a power of $\xi$. On the author's 64-bit Windows computer, the 
operation of  $\rho_p(g)$ on a single vector costs 0.73 ms for $p=3$
and 1.35 ms for $p=255$. That computer has an Intel Core i7-8750H CPU
running at up to 4.0 GHz. These benchmarks are single-threaded.

%%%%%%%%%%%%%%%%%%%%%%%%%%%%%%%%%%%%%%%%%%%%%%%%%%%%%%%%%%%%%%%%%%%%%%%%%%%%%%

%%%%\subfile{chapters/sample_chapter}

% !TeX spellcheck = en_GB

%%%%%%%%%%%%%%%%%%%%%%%%%%%%%%%%%%%%%%%%%%%%%%%%%%%%%%%%%%%%%%%%%%%%%%%
\section{The Golay code $\mathcal{C}$ and its cocode  $\mathcal{C}^*$}
\label{section:Golay}
%%%%%%%%%%%%%%%%%%%%%%%%%%%%%%%%%%%%%%%%%%%%%%%%%%%%%%%%%%%%%%%%%%%%%%%

\subsection{Description of the Golay code $\mathcal{C}$
    and its cocode $\mathcal{C}^*$}
\label{subsection:Golay}

Let $\tilde{\Omega}$ be a set of size~24 and construct the vector
space $\setF_2^{24}$ as $\prod_{i\in \tilde{\Omega}} \setF_2$.
A {\em Golay code} $\mathcal{C}$   is a 12-dimensional linear 
subspace of  $\setF_2^{24}$ whose smallest weight is $8$.
This characterizes the Golay code up to permutation. 
A Golay code has weight distribution
$0^{1}8^{759}12^{2576}16^{759}24^{1}$.  
We identify  the power set of $\tilde{\Omega}$ with
$\setF_2^{24}$ by mapping each subset of $\tilde{\Omega}$ to its
characteristic function, which is a vector in $\setF_2^{24}$.
So we may write $\tilde{\Omega}$ for the Golay code word 
containing 24 ones, and for  elements $d, e$ of $\setF_2^{24}$ 
we write $d \cup e$,  
$d \cap e$ for  their union and intersection
and $d+e$ for their symmetric difference.
Golay code words of length~8 and~12 are called 
{\em octads} and {\em dodecads}, respectively.

The Golay cocode $\mathcal{C}^*$ corresponding to a Golay code
$\mathcal{C}$ is the 12-dimensional quotient space  
$\setF_2^{24} / \mathcal{C}$. For each 
$\delta \in \mathcal{C}^*$ define the weight of $\delta$ to be
the weight of its lightest representative in $\setF_2^{24}$.
Then  $\mathcal{C}^*$   has weight distribution
$0^{1}1^{24}2^{276}3^{2024}4^{1771}$. Here the  lightest
representative is unique if its weight is less than 4.
If $\delta$ has weight 4, there is a set of six mutually disjoint
lightest representatives of $\delta$, and any such subset of
$\setF_2^{24}$  is called a {\em tetrad}.
The parity of a cocode element $\delta$ is defined as the parity 
of its weight $|\delta|$.

There is a natural scalar product $\left<,\right>$ on 
$\mathcal{C} \times \mathcal{C}^*$ given by
$\left<d,\delta \mathcal{C}  \right> = |d \cap \delta| \bmod 2$
for any $d \in \mathcal{C}$, $\delta \in \setF_2^{24}$. 
As usual, a subspace $X^*$ of $\mathcal{C}^*$ is {\em orthogonal}
to a subspace  $X$ of $\mathcal{C}$ if 
$\left<d, \delta\right>=0$
for all $d \in X$, $\delta \in X^*$.

The automorphism group of a Golay code
 $\mathcal{C}$ is the Mathieu group $M_{24}$. 
$M_{24}$ also preserves  $\mathcal{C}^*$, see \cite{Conway-SPLG}. 
For our construction of the Monster we need a
specific instance $\mathcal{C}$ of a Golay code. 
Therefore we assume that the reader is familiar with  
\cite{Conway-SPLG}, Chapter~11, section 1--7 and 9.

There is a 3-dimensional linear code  of length 6 over $\setF_4$
with weight distribution $0^1 4^{45} 6^{18}$ called the
{\em hexacode}. A detailed description of one instance 
$\mathcal{H}_6$  of
the hexacode is given in  \cite{Conway-SPLG}, Chapter 11.
%For our purposes, any instance of the hexacode
%with weight distribution as given above will do.
To be concrete, let $\setF_4 = \{0,1,\alpha, \bar{\alpha}\}$
with $\alpha^2 = \bar{\alpha} = 1 + \alpha$, and let 
$\mathcal{H}_6$ be the subspace of $\setF_4^6$ spanned by 
$({{1001\bar{\alpha}\alpha, 0101\alpha\bar{\alpha}, 
		001111 }})$  as in  \cite{Conway-SPLG}.

We number the elements of $\tilde{\Omega}$ from 0 to 23. 
and arrange them in a table with 4 rows and 6 columns. 
We also assign an element of $\setF_4$ and  a colour to
each row  as follows:
\stepcounter{Theorem}
\begin{equation}
\label{elements:MOG}
\small
\begin{array}{l  c |r|r|r|r|r|r|}
 \cline{3 - 8}
     \mbox{\small white} \mbox{\;} &  0
       & \hphantom{1}0  & \hphantom{1}4 & \hphantom{1}8 
          & 12 & 16 & 20  \\
 \cline{3 - 8}
    \mbox{\small red} &  1
         & 1 & 5 & 9 & 13 & 17 & 21  \\
 \cline{3 - 8}
     \mbox{\small green} &  \alpha 
         & 2 & 6 & 10 & 14 & 16 & 22  \\
 \cline{3 - 8}
    \mbox{\small blue} &  \bar{\alpha}
         & 3 & 7 & 11 & 15 & 19 & 23  \\    
 \cline{3 - 8}
\end{array}  
\end{equation}
This table is called MOG (Miracle Octad Generator) in
\cite{Conway-SPLG}. We start row and column numbers with 0, so element 
$m + 4\cdot n$ of $\tilde{\Omega}$ corresponds to row $m$, 
column $n$. 
This is bad for Fortran and good for C programmers.

Each element $x$ of $\setF_2^{24}$  has
a {\em hexacode value} $\mathfrak{h}(x) \in \setF_4^6$ which is 
calculated  as follows.
We write the entries of $x$ into the MOG. Then for each column
of the MOG we compute a weighted sum of its  nonzero
entries,  where each entry of the MOG  has weight corresponding to the 
element of $\setF_4$ associated with its row. The result
 $\mathfrak{h}(x)$ is a vector in $\setF_4^6$.

Also, for each row or column in the MOG the parity of $x$ in 
that row or column is the number of the nonzero entries in that 
row or column taken modulo~2. 

\begin{Definition}
Let $\mathcal{C}$ be the linear subspace of $\setF_2^{24}$ 
characterized by the following properties:

$x$ is in $\mathcal{C}$ if and only if
$\mathfrak{h}(x) \in  \mathcal{H}_6$ and for all columns of the MOG
the parity of  $x$ in that column is equal to  
 the parity of $x$ in  row 0.
\end{Definition}

In \cite{Conway-SPLG}, Chapter~11 it is shown that 
$\mathcal{C}$ is indeed a Golay code. We call an element of 
$\mathcal{C}$ {\em odd} or  {\em even}, depending
on its parity in  row 0.

\subsection{The 'grey' and the 'coloured' subspaces of
 $\mathcal{C}$ and of  $\mathcal{C}^*$}
\label{subsection:gray:col}

This subsection contains material which is not covered by
\cite{Aschbacher-Sporadic},\cite{Conway:Construct:Monster},
\cite{Conway-SPLG} or
\cite{citeulike:Monster:Majorana}, and which we need
in section~\ref{section_xi} for the first time.
So it may be skipped at first reading.

The construction
of the Golay code in \cite{Conway-SPLG}, Chapter~11.5 
%and also the description of the {\em Three bases subgroup} of
%the automorphism group $\mbox{Co}_1$ of the Leech lattice
%$\Lambda \bmod 2$ in \cite{citeulike:Monster:Majorana},
%Chapter 1.11 
motivates the construction of the
Golay code as a direct sum  $\mathcal{C}$ =  
$\mathcal{G} \oplus  \mathcal{H}$. 
For reasons to be 
explained below, the elements of $\mathcal{G}$ and
$\mathcal{H}$ will be called {\em grey} and {\em coloured},
respectively. We also give a similar decomposition
$\mathcal{C^*}$ =  $\mathcal{G^*} \oplus  \mathcal{H^*}$ of
the Golay cocode $\mathcal{C^*}$ into a direct sum of a grey
and a coloured subspace. 
\Skip{
These decompositions greatly simplify
the explicit construction of an element $\xi$ of the
Monster group $\setM$, which is not contained in the 
maximal subgroup $N_{x0}$ of  $\setM$, in
section~\ref{section_xi}.
}

In section \ref{section:Leech:link}  we will embed both,
the Golay code $\mathcal{C}$ and its cocode  $\mathcal{C}^*$,
into the 24--dimensional Leech lattice $\Lambda$ modulo~2.
Loosely speaking, we will construct an automorphism  $\xi$ of 
$\Lambda/ 2\Lambda$ that fixes $\mathcal{H}$ and $\mathcal{H}^*$, 
and exchanges  $\mathcal{G}$ with $\mathcal{G}^*$. By
construction, $\xi$ is in the automorphism group $\mbox{Co}_1$ 
of $\Lambda/2 \Lambda$, but not in the automorphism group
$M_{24}$ of $\mathcal{C}$. Thus $\xi$ can be lifted to an
element of the monster group in 
$2^{1+24}.\mbox{Co}_1 \setminus 2^{1+24+11}.M_{24}$ 
as required. We remark that our construction of $\xi$ can also be 
achieved by embedding two orthogonal copies of the 12-dimensional 
Coxeter-Todd lattice $K_{12}$ into the Leech lattice, as outlined
in \cite{Conway-SPLG}, Ch.~4.9. Here one copy of $K_{12}$ 
corresponds to $\mathcal{G} \oplus\mathcal{G}^*$ , and the other 
copy to $\mathcal{H} \oplus\mathcal{H}^*$.

We call an element $d$ of $\mathcal{C}$ 
{\em grey} if   in each column of the MOG
all entries of $d$ in that column in rows  1--3 are equal. 
Different columns may have different entries in rows  1--3.
Imagine that each nonzero entry in the MOG switches on
a light with a colour as given for its row in the MOG,
and that the lights in each column are mixed to a common
single colour. Then $d$ is grey if none of these 6 mixed lights 
switched on by $d$ shows any colour apart from white or black.

The even grey elements of $\mathcal{C}$ are precisely those
with equal entries in each column, so that the number of
nonzero columns is even. So there are 32 of them and they
form a 5-dimensional subspace $\mathcal{G}^0$ of $\mathcal{C}$. 
There is also an odd grey element of $\mathcal{C}$, e.g.:
\begin{eqnarray}
\label{odd:gray:Golay:elem}
   g_{0} =  \small
\begin{array}{|r|r|r|r|r|r|}
\hline   0 & 1 & 1 & 1 & 1 & 1  \\
\hline   1 & 0 & 0 & 0 & 0 & 0  \\
\hline   1 & 0 & 0 & 0 & 0 & 0  \\
\hline   1 & 0 & 0 & 0 & 0 & 0  \\
\hline
\end{array}     \; .
\end{eqnarray}

Thus the grey elements form a 6-dimensional subspace $\mathcal{G}$ 
of $\mathcal{C}$.
 
We  define a monomorphism 
$\mathfrak{h}^*: \setF_4^6 \rightarrow \setF_2^{24}$ 
with $\mathfrak{h}(\mathfrak{h}^*(h)) = h$ for $h \in \setF_4^6$
($\setF_4^6$ being considered as a vector space over $\setF_2$)
as follows.
For $h \in  \setF_4^6$ let $\mathfrak{h}^*(h)$ be the unique 
element $x$ of $\setF_2^{24}$ that has zeros in  row 0 of the MOG
and exactly $0$ or $2$ nonzero elements in each column of the MOG
such that the hexacode  value $\mathfrak{h}(x)$ is equal to $h$. 
We put 
$\mathcal{H} = \{\mathfrak{h}^*(x) \mid x \in \mathcal{H}_6  \}$.
Then we have $\mathcal{C} = \mathcal{G} \oplus \mathcal{H}$.
The elements of $\mathcal{H}$ are called {\em coloured}.
So the coloured Golay code words are those without any white light
and with an even number of coloured lights in each column of the MOG.

Let  $\mathcal{G}^*$ be the 
orthogonal complement of $\mathcal{H}$ and let
$\mathcal{H}^*$ be the orthogonal complement of $\mathcal{G}$ 
in  $\mathcal{C}^*$. Elements of $\mathcal{G}^*$ and $\mathcal{H}^*$
are also called grey and coloured, respectively.
The monomorphism $\mathfrak{h}^*$ yields a natural isomorphism 
$\setF_4^6 /  \mathcal{H}_6 \rightarrow \mathcal{H}^*$.
Since more than half of the codewords in  $\mathcal{H}_6$
have weight~4,  the following lemma is obvious:
\begin{Lemma}
\label{lemma:gen:coloured}
$\mathcal{H}$ is generated by  $\mathfrak{h}^*(h)$, where
$h$ runs over the code words in $\mathcal{H}_6$
of weight 4. 
$\mathcal{H}^*$ is generated by   $\mathfrak{h}^*(h)$,
where $h$ runs over the basis vectors of $\setF_4^6$
and their scalar multiples. 
\end{Lemma}

Let $\omega_\infty$ be the element of $\Field_2^{24}$
with entries $1$ in MOG row 0 and zeros in the other
rows. Let $\omega_i$, $i = 0,\ldots,5$ be the element 
of  $\Field_2^{24}$ with entries $1$ in MOG column $i$
and zeros in the other columns. Put 
$g_i = \omega_i + \omega_\infty$. Then $g_0$ is as 
in  (\ref{odd:gray:Golay:elem}) and $g_i$ is obtained 
form $g_0$ by exchanging column~0 with column~$i$ in the
MOG. $g_0,\ldots,g_5$ is a basis of $ \mathcal{G}$.
So $\omega_0,\ldots,\omega_5$ and $\omega_\infty$ represent 
the same element of the cocode $C^*$, which we will demote
by $\omega$. $\omega$ has minimum weight 4 in $C^*$
and the  tetrad corresponding to $\omega$ is
$\{ \omega_0,\ldots,\omega_5  \}$ .
Then  $d \in \mathcal{C}$ is
even if and only if   $\left<d, \omega\right>=0$.

Let $ \gamma_{n} \in \mathcal{C}^*$ 
correspond the vector with an entry 1 in MOG  row 0, column $n$, and
entries 0 elsewhere. Then $( \gamma_0, \dots, \gamma_5)$ is a basis of $ \mathcal{G}^*$.

\begin{Definition}
\label{Def:w}    
For $d \in \mathcal{G}$ and 
$\delta \in \mathcal{G}^*$ let $w(d)$ and $w(\delta)$ be the
weight of $d$ and $\delta$ with respect to the basis 
$( g_0,\ldots, g_5)$ and $( \gamma_0, \dots, \gamma_5)$, respectively.
\Skip{
We extend $w$ from
$\mathcal{G}$ to $\mathcal{C} = \mathcal{G} \oplus \mathcal{H}$  by
decreeing $w(g+h) = w(g)$ for $g \in \mathcal{G}$,
$h \in \mathcal{H}$. 
}
\end{Definition}

\stepcounter{Theorem}

Then  
$\tilde{\Omega} = \sum_{n=0}^5 g_n$, $\omega = \sum_{n=0}^5 \gamma_n$
and we have $w(\tilde{\Omega}) = w(\omega) = 6$. For
$d \in \mathcal{G}$, $\delta \in \mathcal{G}^*$, the weights
$w(d)$ and $w(\delta)$  determine the weight
$|d|$  and the minimum weight $|\delta|$ as follows:
\begin{align}
\label{table_w_d}
\begin{tabular}{|c|c|c|c|c|c|c|c|}
 \hline
   $w(d), w(\delta) $ & 0 & 1 & 2 &  3 &  4 &  5 &  6 \\
 \hline
    $|d|$             & 0 & 8 & 8 & 12 & 16 & 16 & 24 \\
 \hline
    $\min |\delta|$   & 0 & 1 & 2 &  3 &  4 &  3 &  4 \\
 \hline    
\end{tabular}
\end{align}

By definition of $g_0,\ldots,g_5$ and 
$\gamma_0 , \ldots, \gamma_5$  we have: 
\begin{equation}
\label{eqn:g:gamma}
\left<g_m, \gamma_n\right> = 0  \quad 
\mbox{if} \quad m=n \quad \mbox{and} \quad  
\left<g_m, \gamma_n\right> = 1
\quad \mbox{otherwise .} 
\end{equation}
Thus the reciprocal basis of
($g_0,\ldots,g_5$) is ($\gamma_0 + \omega, \ldots \gamma_5 + \omega$) 
and the reciprocal basis of 
($g_0+ \tilde{\Omega},\ldots,g_5+ \tilde{\Omega}$) is
($\gamma_0, \ldots \gamma_5$).

The following fact is rather obvious:
\begin{align}
\label{fact:w:Golay}
\forall \; d \in \mathcal{G}: \quad
w(d + \tilde{\Omega}) = 6 - w(d) \, , \quad
w(d) =  \mbox{parity}(d) \pmod{2} \; .
\end{align}

%%%%%%%%%%%%%%%%%%%%%%%%%%%%%%%%%%%%%%%%%%%%%%%%%%%%%%%%%%%%%%%%%%%%%%%
\section{The Parker loop $\mathcal{P}$}
\label{section:Parker:loop}
%%%%%%%%%%%%%%%%%%%%%%%%%%%%%%%%%%%%%%%%%%%%%%%%%%%%%%%%%%%%%%%%%%%%%%%

\subsection{The definition of the Parker loop}
\label{subsection:Parker:loop}

The Parker loop $\mathcal{P}$ is a non-associative Moufang loop
written multiplicatively and operating on
the set $\mathcal{P} = \mathcal{C} \times \setF_2$, see
\cite{Aschbacher-Sporadic,Conway:Construct:Monster,
	citeulike:Monster:Majorana}.
For any element $d$ of $\mathcal{P}$ we write
$\bar{d}$ for the loop inverse of $d$ in $\mathcal{P}$ and we
write $\tilde{d}$ for the projection of $d$ into $\mathcal{C}$
obtained by dropping the second component of $d$.
This projection is a homomorphism from 
$\mathcal{P}$ to the additive group $(\mathcal{C} , +)$.
We also write $1$,  $-1$, $\Omega$ and $-\Omega$ for the elements
$(0,0)$, $(0,1$), $(\tilde{\Omega},0)$ and $(\tilde{\Omega},1)$ of  
$\mathcal{P}$. Let $d, e, f \in \mathcal{P}$.  

following property characterizes the Parker loop up to isomorphism:
\stepcounter{Theorem}
\begin{align}
\label{Parker:Loop:Comm:Assoc}
      d^2  & =  (-1)^{P(\tilde{d})}
         \qquad \;  \mbox{with} \quad
      P(\tilde{d}) =  {\textstyle\frac{1}{4}}|\tilde{d}| 
                                              \; ,  \nonumber \\
     (de)(ed)^{-1}  & =
         (-1)^{C(\tilde{d} ,\tilde{e})}
         \quad \; \;  \mbox{with} \quad
           C(\tilde{d},\tilde{e}) = 
         {\textstyle\frac{1}{2}} |\tilde{d} \cap \tilde{e}|
                                                   \; , \\
     (d(ef)) ((de)f)^{-1} & = 
       (-1)^{A(\tilde{d},\tilde{e}, \tilde{f})}
            \quad  \mbox{with} \quad
         A(\tilde{d},\tilde{e},\tilde{f}) = 
     {|\tilde{d} \cap \tilde{e} \cap \tilde{f}|}
                                                 \; .  \nonumber
\end{align}
In \cite{Aschbacher-Sporadic} the mappings 
$P: \mathcal{C} \rightarrow \setF_2$, 
$C: \mathcal{C}^2 \rightarrow \setF_2$ and 
$A: \mathcal{C}^3 \rightarrow \setF_2$ given by
(\ref{Parker:Loop:Comm:Assoc})  are called 
{\em power map}, {\em commutator} and  {\em associator},
respectively. Recall that $'+'$ denotes the symmetric difference of 
two sets.
By \cite{Aschbacher-Sporadic},  Lemma~11.1 and Lemma~11.8 or by
direct calculation using the identity
\[
  |\tilde{d} + \tilde{e}| = |\tilde{d}| + |\tilde{e}|
    - 2  |\tilde{d} \cap \tilde{e}| \; ,
\]
for finite sets $\tilde{d}, \tilde{e}$  we obtain:
\begin{Lemma}
\label{eqn:Delta:PCA}
The associator $A$ is a symmetric trilinear form on 
$\mathcal{C}$, and we have:    
\begin{align}
 P(\tilde{d}+\tilde{e}) &=  P(\tilde{d}) + P(\tilde{e}) +
    C(\tilde{d},\tilde{e}) \, ,\\
 C(\tilde{d}+\tilde{e},\tilde{f}) &=  C(\tilde{d},\tilde{f}) +  C(\tilde{e},\tilde{f}) +
    A(\tilde{d},\tilde{e},\tilde{f})  \;      .
\end{align}    
\end{Lemma}

We obviously have
$C(\tilde{d},\tilde{e}) = C(\tilde{e},\tilde{d})$ and by
Lemma~\ref{eqn:Delta:PCA} we have
$A(\tilde{d},\tilde{d},\tilde{e})=0$.
This implies that $\mathcal{P}$ is
{\em diassociative}, i.e. any subloop of  $\mathcal{P}$
generated by two elements is a group, see 
\cite{Aschbacher-Sporadic, citeulike:Monster:Majorana}.
This saves a few brackets in some cases.

From now on we follow the convention in 
\cite{Conway:Construct:Monster}, using the same notation for 
elements of $\mathcal{P}$ and of $\mathcal{C}$. If a function
$F$ has domain $\mathcal{C}^n$ then   $F(d,e,f, \ldots)$
will mean $F(\tilde{d}, \tilde{e}, \tilde{f}, \ldots)$
for $d,e,f,\ldots \in \mathcal{P}$. E.g. $A(d,e,f)$ means
$A(\tilde{d}, \tilde{e}, \tilde{f})$,
$\left<d,\delta\right>$ means $\big<\tilde{d},\delta\big>$
for $\delta \in \mathcal{C}^*$.
But we still  distinguish between $d \in \mathcal{P}$ and
$\tilde{d} \in \mathcal{C}$ whenever we consider $d$ as a
pair $d = (\tilde{d}, \lambda) \in \mathcal{C} \times \setF_2$. 
We use this convention also for functions defined on subsets of
$\tilde{\Omega}$, with the embedding 
$\mathcal{C} \rightarrow  \tilde{\Omega}$. So
$i \in d$ means $i \in  \tilde{d}$,
$\frac{1}{2}|d \cap e|$ means
 $\frac{1}{2}|\tilde{d} \cap \tilde{e}|$,  etc.
%Usually $d \cap e$ is considered as an element of $\mathcal{C}^*$.
We also  use the conventions in \cite{Conway:Construct:Monster} 
for denoting elements of $\mathcal{P}$, $\mathcal{C}$
and  $\mathcal{C}^*$:
\[
\begin{array}{ll}
a,b,c,d,e,f,h  & \mbox{denote elements of $\mathcal{P}$
	    or, loosely, of  $\mathcal{C}$}, \\
\delta, \epsilon, \varphi, \eta &
       \mbox{denote elements of $\mathcal{C}^*$},\\  
i,j,k  & \mbox{denote elements of $\tilde{\Omega}$,
	        also considered  as elements of  $\mathcal{C}^*$ 
             of weight 1}, \\
ij \quad \mbox{and} \quad d \cap e & 
    \mbox{denote the sets $\; \{i,j\}  \;$ and
        $ \; \tilde{d} \cap \tilde{e} \,,\;$ 
          considered  as elements of  $\mathcal{C}^*$ }.          
\end{array}
\]
 
We will also write $\bar{d}$ for the inverse $d^{-1}$ of
$d$ in $\mathcal{P}$. We have $\bar{d} = (-1)^{|d|/4} d$. 
For $d = (\tilde{d}, \mu) \in \mathcal{P} =\mathcal{C} \times \setF_2 $
we put $\mbox{sign}(d) = (-1)^\mu$. Note that this sign
mapping is {\em not} a homomorphism from  $\mathcal{P}$ to
$\{\pm 1\} $.

\subsection{Cocycles for the Parker loop} 
\label{subsection:any:cocycle}
 
Let $V$ be an $n$-dimensional vector space over $\Field_2$ 
and $q: V \rightarrow \Field_2$ be a function with 
$q(0) = 0$. Let $\beta_q: V \times V \rightarrow  \Field_2$
be given by  
$\beta_q (x,y) = q(x+y) + q(x) + q(y)$. If 
$\beta_q$ is bilinear then $q$ is called a 
{\em quadratic form} on $V$ and $\beta_q$ is the
bilinear form {\em associated} with $q$.
Two quadratic forms on $V$ have the same associated
bilinear form if they differ by a linear function on $V$.
The bilinear form $\beta_q$ associated with $q$
is {\em alternating}, i.e.  $\beta_q(x,x) = 0$ for all
$x \in V$. Any alternating bilinear form on $V$ is also
symmetric.

For computations in $\mathcal{P}$ it is  convenient to use a cocycle
$\theta: \mathcal{C} \times \mathcal{C} \rightarrow \setF_2$
such that for  $\mathcal{P} = \mathcal{C} \times \setF_2$ we have:
\begin{equation}
\label{eqn:cocycle:Parker:Loop}
(\tilde{d_1}, \lambda_1 ) \cdot  (\tilde{d_2}, \lambda_2) =
(\tilde{d_2} + \tilde{d_2} , \,
\lambda_1 + \lambda_2 + \theta(d_1,d_2)) \; , \qquad
\tilde{d_1}, \tilde{d_2}
\in \mathcal{C}, \, \lambda_1,\lambda_2 \in \setF_2 \; .
\end{equation}
Cocycles for certain types of loops
have been studied in
\cite{Aschbacher-Sporadic,  DraVoj:CodeLoopBothPar,
    Griess:CodeLoops}.
A cocycle for $\mathcal{P}$ satisfying~\ref{Parker:Loop:Comm:Assoc},
which is quadratic in the first and linear in the second argument, 
has been constructed in \cite{Conway-SPLG}, Ch.~29, Appendix~2.
This can be summarized as follows:

\begin{Lemma}
    \label{Lemma:valid:cocycles}	
    There is a cocycle 
    $\theta:\mathcal{C} \times \mathcal{C} \rightarrow \setF_2$
    for  $\mathcal{P}$ with the following properties:
    \begin{align}
    P(d) &= \theta(d, d)  \; ,  \\
    C(d, e) &=  \theta(d, e) + \theta(e, d) \; , \\
    \theta(d + e, f)  &=  \theta(d,f) +  \theta(e, f)
    + A(d, e, f)  \; , \\
    \theta(d, e + f)  &= 
    \theta(d, e) +  \theta(d, f)  \; .   
    \end{align}
\end{Lemma}

\fatline{Sketch Proof}

Let $b_0,\ldots,b_{11}$ be a basis of $\mathcal{C}$.
Define
\[
    \theta(b_i,b_j) = \left\{  \begin{array}{ll}
       0 & \mbox{if }  i < j \\
       P(b_i)  & \mbox{if }  i = j \\
       C(b_i, b_j)  & \mbox{if }  i > j
    \end{array}    
    \right. \; \,  \quad
\]
Then $\theta(b_i,f)$, $f \in \mathcal{C}$ is  uniquely determined 
by $ \theta(b_i, e + f)  = \theta(b_i, e) +  \theta(b_i, f)$. 
We put
\[
   \theta(\textstyle \sum_i \mu_i  b_i, f)  = 
    \big( \,    \sum_i \mu_i \theta(b_i, f) \, \big) \; + \;
          \sum_{i < j}  \mu_i \mu_j A(b_i, b_j, f) \; . 
\]
By construction of $\theta$ and trilinearity of
$A$, cocycle $\theta$ satisfies (\ref{Lemma:valid:cocycles}.4).
Put $d = \sum_i \mu_i b_i$, $e = \sum_i \nu_i b_i$. By
construction of $\theta$  and trilinearity of $A$ we have
\[
\theta(d + e, f)  +  \theta(d,f) +  \theta(e, f) =
  \textstyle  \sum_{i \neq j}  \mu_i \nu_j A(b_i, b_j, f) \; , 
\]
so together with $A(b_i,b_i,f) = 0$, we obtain
(\ref{Lemma:valid:cocycles}.3).

Put $C'(d,e) = \theta(d,e) + \theta(e,d)$.
(\ref{Lemma:valid:cocycles}.3) and (\ref{Lemma:valid:cocycles}.4) 
imply $C'(d+e,f) = C'(d,f) + C'(e,f) + A(d,e,f)$. We have
$C'(b_i,b_j) = C(b_i,b_j)$. So by
(\ref{eqn:Delta:PCA}.2) and induction over the basis vectors
we obtain $C'=C$, i.e.  (\ref{Lemma:valid:cocycles}.2).  
A similar argument using $\theta(b_i,b_i) = P(b_i)$,
(\ref{eqn:Delta:PCA}.1) and induction over the basis vectors
establishes (\ref{Lemma:valid:cocycles}.1).  

\proofend

An immediate consequence of Lemma~\ref{Lemma:valid:cocycles} is:
	
\begin{Lemma}
\label{Lemma:delta:cocycle:quadratic}   
Let $\theta_1$ be a fixed cocycle satisfying  
Lemma~\ref{Lemma:valid:cocycles}. Then the  cocycles satisfying  
Lemma~\ref{Lemma:valid:cocycles} are precisely the functions
$\theta_1 +\beta$, with $\beta$ an alternating bilinear form on
$\mathcal{C}$.   
\end{Lemma}

Since a cocycle given by Lemma~\ref{Lemma:valid:cocycles} is
linear in its second argument, it may also be interpreted as a function $\mathcal{C} \rightarrow  \mathcal{C}^*$, 
and we let $\theta(d_1)$ be the element of
$ \mathcal{C}^*$ such that  $\left< e, \theta(d) \right>$ 
= $\theta(d,e)$ holds for all $e \in  \mathcal{C}$. 
Similarly, we let $A(d, e)$ be the element of
$\mathcal{C}^*$ such that   $\left< f, A(d, e) \right>$   
= $A(d, e, f)$ holds for all   $f \in  \mathcal{C}$.  
Then for $ d, e \in \mathcal{C}$ 
(or  $ d, e \in \mathcal{P}$) we have:
\stepcounter{Theorem}
\begin{align}
\label{eqn:cocyle:add}
\theta(d + e) =  \theta(d) +  \theta(d) + A(d, e)
\, , \quad \mbox{with} \; \;
A(d, e) =  d  \cap e \; \in \; \mathcal{C}^* \; .
\end{align}

There are cocycles satisfying (\ref{eqn:cocycle:Parker:Loop})
but not Lemma~\ref{Lemma:valid:cocycles}. But we only consider
coycles  which satisfy Lemma~\ref{Lemma:valid:cocycles}.

%%%%%%%%%%%%%%%%%%%%%%%%%%%%%%%% 
\subsection{Selecting a suitable cocycle for the Parker loop}
\label{subsection:cocycle}
%%%%%%%%%%%%%%%%%%%%%%%%%%%%%%%%%%%% 
 
In this subsection we extend the decomposition of $\mathcal{C}$ into a 
grey and a coloured subspace to the Parker loop $\mathcal{P}$. Then we 
construct a cocycle for $\mathcal{P}$ that has some specific properties 
regarding that decomposition as stated in 
Lemma~\ref{Lemma:cocycle:property}.	
This is needed in  section~\ref{section_xi}
for the first time and  may be skipped at first reading.

We also talk about {\em grey} and {\em coloured} elements
of  $\mathcal{P}$. Let  $\mathcal{P}_\mathcal{G}$ and
$\mathcal{P}_\mathcal{H}$ be the set of elements of
$\mathcal{P}$ that are mapped to $\mathcal{G}$ and
$\mathcal{H}$, respectively, by the homomorphism 
$\tilde{\hphantom{.}} : \mathcal{P}  
 \rightarrow  \mathcal{C}$.
\Skip{  
 We put
\[
    \mathcal{P}_\mathcal{G}^+  =
    \left\{ (d, 0) \in \mathcal{P} \mid d \in \mathcal{G}
    \right\} \; .
\]
We also write $g_i$ for  the unique
element $(g_i, 0)$ of  $\mathcal{P}_\mathcal{G}^+$
which maps to $g_i$ in $\mathcal{C}$. 
}

We need a cocycle such that in our representation  of $\setM$
the non-monomial generator $\xi$ of the monster $\setM$
to be constructed in section~\ref{section_xi} becomes as simple as 
possible. More specifically, we will select a cocycle that
is related to the decomposition of $\mathcal{C}$  and
$\mathcal{C^*}$  into a grey and a coloured subspace as
indicated below.

Let $\gamma_{(m,i)}$ be the element of $\setF_2^{24}$ with 
an entry 1 in MOG row $m$, column $i$ and zero elsewhere. 
Thus
$\gamma_{(m,0)} = \gamma_m$.
We define a function
 $\gamma: \setF_2^{24} \rightarrow \setF_2^{24}$ by:
\begin{align}
\label{Def:gamma}
\gamma\left(\sum_{i=0}^5 \sum_{m=0}^3 \mu_{m,i} \gamma_{(m,i)} \right)
 = \sum_{i=0}^5 
    \binom{ \mu_{1,i} + \mu_{2,i} + \mu_{3,i}}{2} \gamma_i
      \; , \qquad 
       \mbox{for} \quad \mu_{m,i} \in \{0,1\} \; .
\end{align}
So $\gamma(x)$ has entry one 
in row~0, column~$n$ of the MOG, if $x$ has at least two nonzero 
entries in column~$n$, ignoring the entry in row 0.
We usually consider $\gamma$ as a function 
$\mathcal{C} \rightarrow  \mathcal{G}^*$ with 
$\gamma(d)$ equal to the coset  $\mathcal{C} \gamma(d)$.
But we occasionally write
$\gamma(d) \cap \gamma(e)$, which is meaningful only for
$\gamma(d)$ and $\gamma(e)$ considered as elements of
$\Field_2^{24}$.
 The restriction of
$\gamma$ to $\mathcal{G}$ is an linear bijection
$\mathcal{G} \rightarrow \mathcal{G}^*$ with 
$\gamma(g_m) = \gamma_m$. We have
\begin{equation}
\label{eqn_gamma_i}
\gamma(g_i) = \gamma(\omega_i) = \gamma_i \; .
\end{equation}

\begin{Lemma}
\label{lemma:sum:gammma:gray:col}   
$\gamma(d+e) = \gamma(d) + \gamma(e)$ for $d \in \mathcal{C}$,
$e \in \mathcal{G}$.
\end{Lemma}

\fatline{Proof}

Let $\phi_j: \Field_2^{24} \rightarrow \Field_2$ be the
function that maps $x \in \Field_2^{24}$ to the coefficient
of the function  $\gamma(x)$ corresponding to $\gamma_j$. 
Then $\gamma = \sum_{j=0}^5 \gamma_j \phi_j$. 
Since $(g_0, \ldots g_5)$ is a basis of $\mathcal{G}$, 
it suffices to show 
\begin{align}
\label{Proof:Lemma:sum:gammma:gray:col}  
\phi_j(d + g_i) =  \phi_j(d) + \phi_j(g_i)  \; , \qquad
\mbox{for} \quad  i,j = 0,\ldots,5 \; .
\end{align}
$\phi_j$ depends only on the co-ordinates of $\Field_2^{24}$
corresponding to  $\gamma_{(m,j)}$ for $m = 1,2,3$,
Since all these coordinates of $g_i$ are zero in case
$i \neq j$, (\ref{Proof:Lemma:sum:gammma:gray:col}) is
established for $i \neq j$.

Assume $i=j$. Then by definition of $\phi_j$ we have:
\[
 \phi_j \left(
    \mu_{1} \gamma_{(1,j)} +  \mu_{2} \gamma_{(2,j)}
      +  \mu_{3} \gamma_{(3,j)}
 \right) =
       \binom{\mu_1 + \mu_2 + \mu_3}{2} \; , \qquad \mbox{for}
   \quad  \mu_1, \mu_2, \mu_3  \in \{0, 1\} \, . 
\]
These three co-ordinates of $g_j$ corresponding to
$\gamma_{(1,j)}, \gamma_{(2,j)}, \gamma_{(3,j)}$ are all 
equal to one.
Thus  (\ref{Proof:Lemma:sum:gammma:gray:col})
follows from $\binom{3-k}{2}$ = 
$\binom{k}{2} + \binom{3}{2} $ $\pmod{2}$, for
$k \in \{0,1,2,3\}$.  
  
\proofend

\stepcounter{Theorem}

Define $w_2: \mathcal{G} \cup \mathcal{G^*}
     \rightarrow  \Field_2$ by
\begin{align}
\label{Def:w_2}
     w_2(d) = \binom{w(d)}{2}  \pmod{2} \; , \qquad 
 \mbox{with $w$ as in Definition~\ref{Def:w}. }
\end{align}
Since $w(\Omega + d) = 6 - w_2(d)$ we have
\begin{align}
\label{w_2_plus_Omega}
w_2(\Omega + d) = w_2(d) +w(d) + 1 \; , \quad 
\mbox{for} \quad d \in \mathcal{G} \; .
\end{align}

A quadratic form $q: \Field_2^n \rightarrow \Field_2$  is called
non-singular if its associated form $\beta_{q}$ is non-singular,  i.e.  $\det(\beta(v_i,v_j)) = 1$ for a basis $(v_1,\ldots,v_n)$
of $\Field_2^n$. A bilinear form which is non-singular and 
alternating is called {\em symplectic}.

\begin{Lemma}
\label{Lemma:scalprod:gamma}
Define $\biscalar{.,.}:
   \mathcal{G} \times  \mathcal{G}\rightarrow \Field_2$  
by  $\biscalar{d,e} =  \left<e, \gamma(d) \right>$.
Then $\biscalar{.,.}$ is symplectic and it is the
bilinear form associated with the quadratic form $w_2$
on $\mathcal{G}$. We also have:
\begin{align}    
\label{Lemma:scalprod:gamma_1}
\left<e, \gamma(d) \right> 
  & =  w(\gamma(d)) \cdot w(\gamma(e)) 
     + w(\gamma(d) \cap \gamma(e))   \qquad
  & \mbox{for} \quad d \in \mathcal{C} , e \in \mathcal{G}    
\end{align}     
\end{Lemma}

\fatline{Proof}

Both sides of (\ref{Lemma:scalprod:gamma_1})  are linear in
$\gamma(d)$ and also in $e$. So it suffices to show
(\ref{Lemma:scalprod:gamma_1}) for $e = g_i$ and  $\gamma(d)$
 being substituted by $\gamma_i$. Thus we have to show
$\left<g_i, \gamma_j \right> = 1 - \delta_{i,j}$, with
$\delta_{i,j}$ the Kronecker delta. But this is an
immediate consequence of (\ref{eqn:g:gamma}). 

Hence $\biscalar{g_i,g_j}=1 - \delta_{i,j}$.
For  $i=1,2$ assume 
$d_i = \sum_{n=0}^5 \mu_{i,n} g_n$, $\mu_{i,n} \in \{0,1\}$. Since $\binom{n}{2} = \sum_{m=0}^{n-1} m $, we have
$
w_2(d_i) = \sum_{0 \leq m < n \leq 5}  \mu_{i,m} \mu_{i,n}. 
$
So $w_2$ is a quadratic from on   $\mathcal{G}$.  Let
$\beta': \mathcal{G} \times  \mathcal{G}\rightarrow \Field_2$ 
be given by 
$(d_1, d_2) \mapsto w_2(d_1) + w_2(d_2) + w_2(d_1 + d_2)$.
Then
$\beta'(d_1, d_2) = \sum_{m \neq n}  \mu_{1,m} \mu_{2,n}.$
Thus $\beta'$ is linear in both arguments and 
we have $\beta'(g_m, g_n) = 1 - \delta_{m,n}$.
Hence $\biscalar{.,.} = \beta'$, i.e. $\biscalar{.,.}$ is
associated with $w_2$.

We have just shown that $\det \beta'$ (with respect to
the basis $g_1\ldots,g_5$)  has entry
$1 - \delta_{m,n}$ in row $m$, column $n$, so direct
calculation yields $\det \beta' = 1 \pmod{2}$.

\Skip{ % Alternative proof for (\ref{Lemma:scalprod:gamma_1})
Put $\delta = \gamma(d) \in \Field_2^{24}$, 
$\epsilon = \gamma(e)\in \Field_2^{24}$. Let 
$\epsilon' = e \cap \omega_\infty$
with $\omega_\infty$ as in section~\ref{subsection:gray:col}.
Then 
$\delta \cup \epsilon \cup \epsilon' \subset \omega_\infty$
and we have $\epsilon' = \epsilon$ if $\epsilon$ is even and
$\epsilon' = \epsilon + \omega_\infty$ if $\epsilon$ is odd.
Thus
\[
   \left< \delta, e \right> = |\delta \cap e|
    =   w(\delta \cap \epsilon') = w(\delta \cap \epsilon) 
         + w(\epsilon) \cdot w(\delta \cap \omega_\infty) 
     = w(\delta \cap \epsilon) 
+ w(\delta) \cdot w(\epsilon) \; . 
\] 
}

\proofend

\vspace{1ex}

We claim the existence of a specific cocycle for 
$\mathcal{P}$ as follows:

\begin{Lemma}
\label{Lemma:cocycle:property}	
There is a cocycle $\theta$ for $\mathcal{P}$  satisfying 
Lemma~\ref{Lemma:valid:cocycles}
with the following properties:
\[
\begin{array}{rclcl}
  \theta(d + \tilde{\Omega}) & = & \theta(d) 
       &\qquad \mbox{for} \qquad &
        d \in  \mathcal{C} \; ,
          \\ 
  \theta(e) & =  & (w(e)-1)  \gamma(e)
             + w_2(e) \omega
          &  \qquad \mbox{for} \quad &
       e \in  \mathcal{G}  \, ,  \quad
      \omega \; \mbox{as in section \ref{section:Golay}},
        \\         
  \theta(e,h) & = & 0 & \qquad  \mbox{for} \qquad &
      e \in  \mathcal{G} , \;  h \in   \mathcal{H} \; ,\\
  \theta(h, e) & = & \left< e, \gamma(h) \right> & 
       \qquad  \mbox{for} \qquad &
e \in  \mathcal{G} , \;  h \in   \mathcal{H} \; . 
\end{array}
\]
\end{Lemma}
	
\fatline{Proof}

Choose a basis $(b_0,...,b_{11})$ of
$\mathcal{C} = \mathcal{G} \oplus \mathcal{H}$ with
$b_m = g_m\in \mathcal{G}$ for
$m=0,\ldots,5$, $b_6,\ldots,b_{11} \in \mathcal{H}$, and
$g_m$ as in section~\ref{section:Golay}.
Since two cocycles for $ \mathcal{C}$ differ by an alternating
bilinear form, a cocycle $\theta$ is uniquely determined by
decreeing the values $\theta(b_m,b_n) $ for all
$0 \leq m < n < 12$. So we put $\theta(b_m,b_n)=0 $ for $m < n$.
%In \cite{Conway-SPLG} a coycle $\theta$ is constructed in the 
%same way for an arbitrary basis of $\mathcal{C}$.
We have 
$ |g_m \cap g_n| = 4$ for
$m \neq n$ So for $0 \leq m < n < 6$ we have:
\[   
  \theta(b_m,b_n) +  \theta(b_n,b_m) 
  = C(g_m, g_n) 
 	 =  {\textstyle{\frac{1}{2}}}
 |g_m   \cap g_n|  = 0 \pmod{2}	 \; ,	
\]	 	
by (\ref{Parker:Loop:Comm:Assoc})
  and Lemma~\ref{Lemma:valid:cocycles}. For $m < 6$ we have
$\theta(b_m,b_m) =   P(g_i)
  =   {\textstyle{\frac{1}{4}}} | g_m |
    = 0 \pmod{2}$.
Thus $\theta(b_m,b_n) = 0$ for all $m <6, n < 12$. Hence 
 $\theta(b_m)=0$  for $m < 6$. Thus by
(\ref{eqn:cocyle:add}) the restriction of $\theta$ 
(as a mapping $\mathcal{C} \rightarrow \mathcal{C}^*$)
 to $\mathcal{G}$ is symmetric under permutations
of the basis vectors $b_0,\ldots, b_5$.
$\theta(e)$ can be computed
for all $e \in \mathcal{G}$ by (\ref{eqn:cocyle:add})
Due to the permutation symmetry of $\theta$ it suffices to
verify the formula for $\theta(e)$  for the cases
 $\theta(e_m)$ with $e_m = \sum_{n=0}^m  g_n$, 
$m=1,\ldots 5$. This can easily be done by hand calculation.
(Row~0 of $\theta(e_m)$  in the MOG is
{\small{000000, 001111, 111111, 111100, 000000}} for
$m=1,\ldots, 5$; the other rows are zero.)

That way we obtain $\theta(\Omega)=0$, so we have 
$\theta(d + \Omega)  =  \theta(d)$ by (\ref{eqn:cocyle:add}).
We have $\theta(e,h) =0$ for $e\in \mathcal{G}$,
$h \in \mathcal{H}$, by construction of $\theta$. Thus
$\theta(h,e)$ = $C(e,h)$ = $\frac{1}{2}|e \cap h|$.
Note that $e \cap h$ has 0 or 2 nonzero entries in each
column of the MOG and zeros in  row~0. Hence  
\[
 \theta(h,e) =
\frac{1}{2}|e \cap h| = |\gamma(e) \cap \gamma(h)| 
= w(\gamma(e) \cap \gamma(h)) =  \left< e, \gamma(h) \right>
   + w(\gamma(e)) \cdot  w(\gamma(h)) \; .
\]   
by definition of $\gamma$ and
(\ref{Lemma:scalprod:gamma_1}). Since $\gamma(h)$ is even,
this proves the  formula for  $\theta(h,e)$.

\proofend

\begin{Corollary}
\label{Corol:cocycle:colored}    
For every $h \in \mathcal{H}$ there is a cocycle $\theta$ on
$\mathcal{C}$ satisfying Lemma~\ref{Lemma:cocycle:property}	
with $\theta(h) = \gamma(h) \in \mathcal{G}^*$.
\end{Corollary}    

\fatline{Proof}

In the proof of  Lemma~\ref{Lemma:cocycle:property}	
we simply choose a basis $b_0,\ldots,b_{11}$ of
$\mathcal{C}$  with $b_6 = h$.

\proofend

(\ref{table_w_d}) and
Lemma~\ref{Lemma:cocycle:property}  imply that
for every $e \in \mathcal{G} ,  h \in \mathcal{H} $  we have:
\begin{align}
\label{diff2:P_gray:col}
\begin{array}{c}
P(e) = w(e) \cdot w_2(e) \, , \;
P(h) = w_2(\gamma(h)) \, , \;
C(e,h) = \left<e, \gamma(h)\right> \; .
\end{array}
\end{align}

From now on we assume that the cocycle for $\mathcal{P}$ satisfies
the conditions in Lemma~\ref{Lemma:cocycle:property}.

% !TeX spellcheck = en_GB

%%%%%%%%%%%%%%%%%%%%%%%%%%%%%%%%%%%%%%%%%%%%%%%%%%%%%%%%%%%%%%%%%%%%%%%
\section{Automorphisms of the Parker loop $\mathcal{P}$}
\label{section:auto:Parker}
%%%%%%%%%%%%%%%%%%%%%%%%%%%%%%%%%%%%%%%%%%%%%%%%%%%%%%%%%%%%%%%%%%%%%%%

The {\em center} $Z(G)$ of a group or a loop $G$
is the set of elements $d$  of $G$ such that $d$ commutes
and associates with all elements of  $G$.
Any element of the Parker loop  $\mathcal{P}$ squares to 
$\pm 1$ and we have $Z(\mathcal{P}) = \{\pm1, \pm \Omega\}$. 
Thus any automorphism $\pi$ of $\mathcal{P}$
fixes $\{\pm1\}$  and maps $\Omega$ to $\pm \Omega$. 
$\pi$ is called {\em even} if $\pi(\Omega) = \Omega$ and {\em odd}
if  $\pi(\Omega) = -\Omega$.
%Let $\pi$ be an automorphism of the Parker loop $\mathcal{P}$.
Since $\pi$ fixes $\{\pm 1\}$, it maps to a unique  automorphism
$\tilde{\pi}$ of  $\mathcal{C} = \mathcal{P}/\{\pm1\} $. 
Then $\tilde{\pi}$ preserves  power map $P$, commutator
$C$ and associator $A$, but it need not preserve the Golay
code on the vector space $\mathcal{C}$.
An automorphism $\pi$  of  $\mathcal{P}$ is called
a {\em standard automorphism} in \cite{Conway:Construct:Monster},
if $\tilde{\pi}$ preserves the Golay code.
E.g. for a nonzero even $\delta \in \mathcal{C}^*$ the mapping
$d \mapsto d \cdot \Omega^{\left< d, \delta \right>} $
is  a non-standard automorphism of  $\mathcal{P}$,
see \cite{citeulike:Monster:Majorana}, section~1.6
for background.
In the sequel we only deal with  the group $\AutStP$ 
of standard automorphisms $\pi$ of 
$\mathcal{P}$. For any   $\pi \in \AutStP$ 
we have  $\tilde{\pi} \in M_{24}$. 
For any $\tilde{\pi} \in M_{24}$ there are 
precisely $2^{12}$ standard automorphisms $\pi$ of
$\mathcal{P}$ mapping to  $\tilde{\pi}$,
see \cite{Conway:Construct:Monster,Conway-SPLG}.

Fix a cocycle $\theta$ for $\mathcal{P}$ that satisfies
Lemma~\ref{Lemma:valid:cocycles}. For any  $\pi \in \mathcal{P}$ define 
$\theta_\pi: \mathcal{C} \times \mathcal{C} \rightarrow \setF_2$
by
\begin{align}
\label{eqn:def_theta_pi}
\theta_\pi(d, e)  \, = \, 
\theta(d^\pi,e^\pi)
+ \theta(d, e) \; .
\end{align}	

Then the following lemma is useful for effective computations in the
group $\AutStP$:	
	
\begin{Lemma}
\label{Lemma:standard:automorphism}	
For any standard automorphism $\pi$ of $\mathcal{P}$ the
mapping $\theta_\pi$ is a alternating bilinear  form on 
$ \mathcal{C}$ depending on $\tilde{\pi}$
only. There is a unique quadratic form 
$q_\pi$ on $\mathcal{C}$ with associated
bilinear form $\theta_\pi$, such that for any
$d = (\tilde{d}, \lambda) \in  \mathcal{P}$, 
$\tilde{d} \in \mathcal{C}$, $\lambda \in \setF_2$ we have 
\[
    (\tilde{d}, \lambda)^\pi = 
      (\tilde{d}^{\pi}, \lambda + q_\pi(\tilde{d})) \; .
\]
\end{Lemma}

\fatline{Proof}

Define 
$\theta^{(\pi)}: \mathcal{C} \times \mathcal{C} 
        \rightarrow \setF_2$ by
$\theta^{(\pi)}(d, e) =   
  \theta(d^\pi,e^\pi)$.
Since $M_{24}$ acts as a group of linear transformations on
$\mathcal{C}$ and  preserves $P, C$ and $A$, the
function $\theta^{(\pi)}$ satisfies the conditions for a 
cocycle in Lemma \ref{Lemma:valid:cocycles}. So
$\theta_\pi =  \theta^{(\pi)} + \theta$ is an alternating bilinear
form by Lemma~\ref{Lemma:delta:cocycle:quadratic}.   
By construction, $\theta_\pi$ depends on $\tilde{\pi}$ only.        
Let $q_\pi$ be any quadratic from with associated form $\theta_\pi$.
Write $(\tilde{d}, \lambda)^{\pi,q_\pi}$ for 
$ (\tilde{d}^{\pi}, \lambda + q_\pi(\tilde{d}))  $.
For $\tilde{d}, \tilde{e} \in \mathcal{C}$ and
$\lambda, \mu \in \Field_2$ we have:
\begin{align*}
  (\tilde{d}, \lambda)^{\pi,q_\pi} \cdot
      (\tilde{e}, \mu)^{\pi,q_\pi}  &=
    \left(
      \tilde{d}^{\pi}+ \tilde{e}^{\pi} ,
   \lambda  + q_\pi(\tilde{d}) + \mu
    + q_\pi(\tilde{e}) +\theta(d^\pi, e^\pi)
    \right) \\
   &=   
    \left(
     \tilde{d}^{\pi}+ \tilde{e}^{\pi}  ,
     \lambda +  \mu +  q_\pi(\tilde{d} + \tilde{e})
     + \theta_\pi(d,e)   +\theta(d^\pi,  e^\pi)
     \right) \\
   &=   
    \left(
    \tilde{d}^{\pi}+ \tilde{e}^{\pi} ,
     \lambda +  \mu  + \theta(d ,e)  
       + q_\pi(\tilde{d} + \tilde{e})
     \right) \\
   &=   
    \left(
       \tilde{d} + \tilde{e} ,
        \lambda +  \mu  +  \theta(d, e)
   \right)^{\pi, q_\pi} \\
   &=  
       (\tilde{f}, \nu)^{\pi,q_\pi} \; ,  
        \qquad \qquad \mbox{with} \quad 
       (\tilde{f}, \nu) = 
      (\tilde{d}, \lambda) \cdot (\tilde{e}, \mu)
      \; .
\end{align*}
So the mapping  
$(\tilde{d}, \lambda) \mapsto (\tilde{d}, \lambda)^{\pi,q_\pi}$ 
is a standard automorphism of $\mathcal{P}$ which maps to the
element  $\tilde{\pi}$ of $M_{24}$.
There are $2^{12}$ different standard automorphisms $\pi$ of
$\mathcal{P}$ mapping to the same element  $\tilde{\pi}$ 
of $M_{24}$, 
and there are $2^{12}$ different quadratic forms on $\mathcal{C}$
with the same associated bilinear form $\theta_\pi$. Hence for any
such $\pi$ there is a unique $q_\pi$ satisfying the lemma.

\proofend

\stepcounter{Theorem}

A standard  automorphism $\delta$ of $\mathcal{P}$ that maps
to the neutral element $\tilde{\delta}=1$ of $M_{24}$ is called 
a {\em diagonal automorphism}. In that case the bilinear
form $\theta_\delta$ in Lemma~\ref{Lemma:standard:automorphism}	
is 0 and a quadratic form associated with the bilinear form 0
is a linear form in 
$\hom(\mathcal{C}, \setF_2) \cong \mathcal{C}^*$, which we also
denote by $\delta$. So a diagonal automorphism
$\delta \in  \mathcal{C}^*$ maps $d$ to
$d \cdot (-1)^{\left< d, \delta \right>}$.
The parity of $\delta$ as a diagonal automorphism agrees with 
the parity of $\delta$ as an element of $\mathcal{C}^*$.

$\AutStP$ is a non-split extension with normal
subgroup $\mathcal{C}^*$ and factor group $M_{24}$. Even if we
assume that calculations in  $\mathcal{C}$, $\mathcal{C}^*$,
$\mathcal{P}$ and $M_{24}$ are easy, the calculations in
$\AutStP $ are still quite technical. 
In the remainder of this section we study such  calculations.

Using Lemma~\ref{Lemma:standard:automorphism},	
 calculation in  $\AutStP$ can be done as follows.
Choose a basis $(\tilde{b}_0,\ldots,\tilde{b}_{11})$ of 
$\mathcal{C}$ and for
$\pi \in M_{24}$ let $[\pi] \in  \AutStP$ be defined by
\[
     (\tilde{b}_i,0)  \stackrel{[\pi]}{\longmapsto}
         (\tilde{b}_i^\pi,0) \; , \quad i = 0, \ldots, 11 \; .
\] 
Then any element of $\AutStP$ can uniquely by written in the form
$\delta \cdot [\pi]$, $\delta \in \mathcal{C}^*$, $\pi \in M_{24}$,
and we have $ [\pi] \cdot \delta^\pi = \delta \cdot [\pi]$.

For $\pi, \pi' \in M_{24}$ we have 
$[\pi \pi'] = \vartheta(\pi, \pi') \cdot [\pi \pi']$, where
$ \vartheta(\pi, \pi') \in \mathcal{C}^*$ is given by
\begin{align}
\label{formula:mult:AutStP}
   \big< \tilde{b}_i,  \vartheta(\pi, \pi') \big> =
      q_{[\pi']}(b_i^{\pi}) \; , \quad i = 0, \ldots, 11 \; ,
\end{align}
and $q_{[\pi']}$ is the unique quadratic form on $\mathcal{C}$ with
associated bilinear form  $\theta_{\pi'}$ satisfying 
$q_{[\pi']}(\tilde{b}_i) = 0$ for $i=0,\ldots,11$. Here
$\theta_{\pi'}$ is as in (\ref{eqn:def_theta_pi}).
Note that $\theta_{\pi'}$ and hence also  $q_{[\pi']}$ can easily
be computed from $\theta(b_i)$ and $\theta(b_i^{\pi'})$, 
$i=0,\ldots,11$.    

The proof of (\ref{formula:mult:AutStP}) is a simple calculation based on the equation
\begin{align}
\label{formula:op:mult:AutStP}
       (\tilde{d}, \lambda)^{[\pi] [\pi']}= 
   (\tilde{d}^{\pi \pi'}, \lambda + q_{[\pi]}(\tilde{d})
     +    q_{[\pi']}(\tilde{d}^\pi)   ) \; , \quad
  %   (\tilde{d}, \lambda) \in \mathcal{P}, \;
         \tilde{d} \in \mathcal{C}, \; \lambda \in \Field_2 , \;
       \pi, \pi' \in \AutStP \, ,
\end{align}
with $q_{[\pi]}$ and $q_{[\pi']}$  as in 
Lemma~\ref{Lemma:standard:automorphism}. 
(\ref{formula:op:mult:AutStP}) is a consequence of 
Lemma~\ref{Lemma:standard:automorphism}.

% !TeX spellcheck = en_GB

%%%%%%%%%%%%%%%%%%%%%%%%%%%%%%%%%%%%%%%%%%%%%%%%%%%%%%%%%%%%%%%%%%%%%%%
\section{The code loop group $N_0$ of the Parker loop $\mathcal{P}$}
\label{section:N0}
%%%%%%%%%%%%%%%%%%%%%%%%%%%%%%%%%%%%%%%%%%%%%%%%%%%%%%%%%%%%%%%%%%%%%%%

We define a group $N$ which acts as a permutation group
on the elements of $\mathcal{P}^3$. It will turn out that 
$N$ is a fourfold cover of a maximal subgroup $N_0$ of the
monster $\mathbb{M}$. The elements
of  $\mathcal{P}^3$ are called {\em triples}.
$N_0$ has structure $2^2.2^{11}.2^{22}.(S_3 \times M_{24})$ and 
is the normalizer of a certain four-group $\{1,x,y,z\}$
in $\mathbb{M}$, see \cite{Conway:Construct:Monster}. 
Following Conway's construction of the Monster $\mathbb{M}$  
in  \cite{Conway:Construct:Monster}, we define maximal subgroups
$G_{x0}, G_{y0}$ and $G_{z0}$ of $\mathbb{M}$, each of structure
$2_+^{1+24}.\mbox{Co}_1$, which centralize the elements
$x,y$ and $z$ of the four-group, respectively.

We use the ATLAS conventions \cite{Atlas} for group structures.
$G_1 . G_2$ denotes a group extension with normal subgroup $G_1$
and factor group $G_2$,  $G_1 : G_2$ denotes a split extension
and   $G_1 \times G_2$ a direct product. All these operators 
associate to the left with the same precedence, and they imply
that the given decomposition is invariant. E.g. 
$G = G_1 \times G_2 : G_3 \, .  \, G_4$ means 
$G = ((G_1 \times G_2) : G_3) \, .  \, G_4$ and implies the existence
of normal subgroups $G_1$ and $G_1 \times G_2$ of $G$. 
$2^n$ is an elementary Abelian group of order $2^n$, $S_n$ is the
symmetric permutation group of $n$ elements, $\mbox{Co}_1$ is
the automorphism group of the Leech lattice $\Lambda$ modulo 2, 
and  $\mbox{Co}_0 = 2.\mbox{Co}_1$ is the automorphism group of
$\Lambda$.

The Leech lattice $\Lambda$ is discussed in 
section~\ref{section:Leech:link}. Following
\cite{Conway:Construct:Monster} we construct a 196884-dimensional 
real monomial representation $196884_x$ of the group
$N_{x0} = N_{0} \cap G_{x0}$, with $|N_0| / |N_{x0}|$ = 3,
in section~\ref{section:rep:Nx1}. 
Then in section~\ref{section:rep:Nx0:N0} we extend 
$196884_x$ to a representation of $N_0$ by adding a
{\em triality} element $\tau$ which (by conjugation)
cyclically exchanges the elements $x,y,z$ of the four-group and 
also their centralizers $G_{x0}, G_{y0}, G_{z0}$.
  
The group $2_+^{1+24}$ is an extraspecial 2-group, which we will
discuss in section~\ref{section_xi}. There we explicitly construct
the representation of an element 
$\xi \in G_{x0} \setminus  N_{x0}$ in the $N_0$-module  $196884_x$.
Since $N_{x0}$ is a maximal subgroup of $G_{x0}$, this extends
$196884_x$ to a representation of a group generated
by  $G_{x0}$ and $N_0$, which is visibly equivalent to
the representation of $\mathbb{M}$ in \cite{Conway:Construct:Monster}.
Essentially, the main added value of this paper compared
to~ \cite{Conway:Construct:Monster}  is an explicit
description of the representation of an element 
$\xi$ of $G_{x0} \setminus  N_{x0}$  
of the Monster, so that a programmer may implement this 
construction with not too much effort.

As in \cite{Conway:Construct:Monster}, 
 we define various subgroups of $\mathbb{M}$  shown 
in Figure~\ref{figure:subgroups:monster} .
All these subgroups of $\mathbb{M}$  agree with the corresponding
subgroups of $\mathbb{M}$ in \cite{Conway:Construct:Monster}. 
In our construction we change a few signs compared to  
\cite{Conway:Construct:Monster} in order to simplify the
representation of the element $\xi$. 

\newcommand{\verylonglrarrow}{\xleftrightarrow{\phantom{mm}}}

\begin{figure}[!h]
    \centering
\begin{tikzpicture}
\matrix (m) [matrix of math nodes,row sep=1.0em,column sep=1.0em,minimum width=2em]
{  
    &            & \mathbb{M} &           &           &
    \verylonglrarrow   & \mbox{The Monster}                   \\
    G_{x0}  &            &            &  G_{y0}   &  G_{z0}   &    
    \verylonglrarrow  & 2_+^{1+24}.\mbox{Co}_1                \\   	
    \\
    &            &   N_0     &           &           &
    \verylonglrarrow  &
    2^2 . 2^{11} . 2^{22} . (S_3 \times M_{24})   \\   	
    N_{x0}  &            &            &  N_{y0}   &  N_{z0}   &    
    \verylonglrarrow  &      2^{1+24} . 2^{11} . M_{24}      \\   	
    & N_{xyz0}   &            &           &           &
    \verylonglrarrow  &     2^2 . 2^{11} . 2^{22} . M_{24}   \\
    &            &   Q_0     &           &           &
    \verylonglrarrow  &  2^2 . 2^{11} . 2^{22} . S_3         \\
    Q_{x0}  &            &            &  Q_{y0}   &  Q_{z0}   &    
    \verylonglrarrow  &      2_+^{1+24}                      \\   	
    & Q_{xyz0}   &            &           &           &
    \verylonglrarrow  &     2^2 . 2^{11}                     \\
};         	
\path[] 
(m-2-1) edge (m-1-3)
(m-2-4) edge (m-1-3)
(m-2-5) edge (m-1-3)
(m-4-3) edge (m-1-3)  % from N_0 to \mathbb{M}
(m-5-1) edge (m-2-1)
(m-5-4) edge (m-2-4)
(m-5-5) edge (m-2-5)
(m-5-1) edge (m-4-3)
(m-5-4) edge (m-4-3)
(m-5-5) edge (m-4-3)
(m-6-2) edge (m-5-1)  % lines from N_{xyz0} upward
(m-6-2) edge (m-5-4)
(m-6-2) edge (m-5-5)
(m-7-3) edge (m-4-3)  % from Q_0 to  N_0 
(m-8-1) edge (m-5-1)
(m-8-4) edge (m-5-4)
(m-8-5) edge (m-5-5)
(m-8-1) edge (m-7-3)
(m-8-4) edge (m-7-3)
(m-8-5) edge (m-7-3)
(m-9-2) edge (m-6-2)  % lines from Q_{xyz0} upward
(m-9-2) edge (m-8-1)
(m-9-2) edge (m-8-4)
(m-9-2) edge (m-8-5)
;
) ;
\end{tikzpicture} 
    \caption{Some subgroups of the Monster $\mathbb{M}.$ }  
    \label{figure:subgroups:monster}    
\end{figure}

We start with the definition of a group $\hat{N}$ which has
a simpler structure than $N$ and we will define $N$ as
a subgroup of $\hat{N}$ of index~2.
We define the following generators
$x_d, y_d, z_d, \nu_\pi, \tau, x_\tau, y_\tau, z_\tau$ of $\hat{N}$ 
by their action on a triple  $(a,b,c)$ in $\mathcal{P}^3$:
\begin{align*}
&
x_d: (\bar{d} a d, \bar{d} b, c d)    \qquad
y_d: (a d, \bar{d} b d, \bar{d} c)    \qquad
z_d: (\bar{d} a,  b d, \bar{d} c d) \, ,  \\
&
\nu_\pi :  (a^\pi, b^\pi, c^\pi)  \qquad  \tau: (c, a, b) \, , \\
&    
x_\tau: (\bar{a}, \bar{c}, \bar{b})   \qquad   
y_\tau: (\bar{c}, \bar{b}, \bar{a})   \qquad  
z_z: (\bar{b}, \bar{a}, \bar{c}) \, ,    
\end{align*}
with
$d$ ranging over $\mathcal{P}$ and $\pi$ ranging over $\AutStP$.
We put
\[
x_\pi = x_\tau^{|\pi|} \nu_\pi, \quad  y_\pi = y_\tau^{|\pi|} \nu_\pi,
\quad  z_\pi = z_\tau^{|\pi|} \nu_\pi, \quad
\mbox{with $|\pi| =1$ for odd and $|\pi| =0$ for even $\pi$} .
\]
So e.g. $x_\pi$ maps $(a,b,c)$ to 
$(\bar{a}^\pi,\bar{c}^\pi,\bar{b}^\pi)$ for odd $\pi$. We define
$N$ to be the subgroup of $\hat{N}$ generated by 
$x_d, y_d, z_d, x_\pi, y_\pi, z_\pi$.

Next we show some relations in  $\hat{N}$. 
For elements $u, v$ of a group we write $[u,v]$ for the
commutator $u^{-1} v^{-1}uv$ and  $u^v$ for $v^{-1}uv$.
We abbreviate $x_{-1} , y_{-1}$
and $z_{-1}$ to $x, y$ and $z$, respectively, as in
\cite{Conway:Construct:Monster}.
The groups $\hat{N}$ and  $N$ have a visible symmetry with 
respect to cyclic permutations of the letters $x, y$ and $z$, 
which we will call {\em triality}. 
Conjugation with $\tau$ just performs such a cyclic permutation. 
We often just write only one of the three possible cyclic 
permutations in a formula or a definition, and we state that
the others are obtained by triality.    We use the symbol
$\Cyclicxyz$ to indicate that a relation remains valid
if the letters $x, y, z$ are cyclically exchanged, e.g.:
\[
x_d y_d= {z}_{d}  z^{P(d)} \; \Cyclicxyz  
\qquad \mbox{means} \quad  
x_d y_d= {z}_{d} z^{P(d)} , \;
y_d z_d= {x}_{d} x^{P(d)} , \;
z_d x_d= {y}_{d} y^{P(d)} \, .
\]

\begin{Theorem}
\label{Thm:relations:N}
The following relations define the group $\hat{N}$:
\[
\left.
\begin{array}{clll}
& x_d x_e = x_{de} x_{A(d,e)}, &
\nu_\delta \nu_\epsilon = \nu_{\delta\epsilon},
 \qquad \hphantom{.} & 
 [x_d, \nu_\delta] = x^{\left<d,\delta\right>}, 
\\
& [x_d, y_e] = \nu_{A(d,e)} z^{C(d,e)}, \; \hphantom{.} &
x_d y_d z_d = 1,
\\
& x_d \nu_\pi = \nu_\pi x_{d^\pi} , &
\nu_{\pi} \nu_{\pi'}  = \nu_{\pi\pi'} , 
\\
& \tau^3 = x_\tau^2 = 1, & \tau = y_\tau x_\tau, & 
   x_\tau \tau = \tau^2 x_\tau, 
\\ 
& [\nu_\pi, \tau] = 1, & x_d \tau = \tau y_d, &
   x_d \tau^2 = \tau^2 z_d,
\\
&
[\nu_\pi,x_\tau] = [x_d, x_\tau] = 1, &
x_\tau y_d = z_d x_\tau, & x_\tau z_d = y_d x_\tau \, .
\end{array}
\right\} \; \Cyclicxyz.
\]
  		
\end{Theorem}

\fatline{Proof}

Most of these relations can be shown by calculations within
associative subloops of $\mathcal{P}$ generated by two 
elements, which we leave to the reader. We first show
$x_d x_e = x_{de} x_{A(d,e)}$. We have
\begin{align*}
   (a,b,c) & \stackrel{x_d x_e}{\longmapsto}
    \left(   \bar{e}(\bar{d} a d) e, \bar{e}(\bar{d} b),
      (c d) e \right) \\
    & \; =  \;
    \left( (-1)^{C(a, d)+C(a,e )} a, \,
     (-1)^{A(b,d,e)}  \conj{(de)}  b, \,
       (-1)^{A(c,d,e)}   c(d e) \right)\\
    & \; =  \;
     \left( (-1)^{C(a, d+e)+A(a,d,e )} a, \,
     (-1)^{A(b,d,e)}  \conj{(de)}   b, \,
      (-1)^{A(c,d,e)} c(d e) \right) \\
    & \; =  \;
     \left( (-1)^{A(a,d,e )}  \conj{(de)} a (de), \,
     (-1)^{A(b,d,e)}  \conj{(de)}   b, \,
      (-1)^{A(c,d,e)} c(d e) \right) \; .     
\end{align*}
For the last step, note that 
 $C(a, d)+C(a,e ) = C(a, d+e)+A(a,d,e )$
 by (\ref{eqn:Delta:PCA}.2)
Clearly,  $x_{de} x_{A(d,e)}$ maps
$(a,b,c)$ to the same triple.

We also show $[x_d, y_e] = z_{A(d,e)} z^{C(d,e)}$. We have
\begin{align*}
(a,b,c) & \stackrel{[x_d, y_e]}{\longmapsto}
 \left( (-1)^{C(a,d)+C(a+e,d)}  a,  (-1)^{C(b+d,e)+C(b,e)}  b,
         \bar{e}\left((e(c\bar{d}))d\right) \,  \right)  \\
     & \; = \;   
 \left(  (-1)^{C(d,e)+A(a,d,e)} a,  (-1)^{C(d,e)+A(b,d,e)} b,
     (-1)^{A(c,d,e)}  c
     \right)   \; .
\end{align*}
Note that $A(c+d,d,e)=A(c,d,e)$. Clearly, 
$z_{A(d,e)}z^{C(d,e)}$ maps $(a,b,c)$ to the same triple. 

Using the relations shown we may convert any word in the
 generators of $\hat{N}$ to the form:
\[
   x_d y_e x_\pi  \tau^m  x_\tau^n, \quad  
   0 \leq m < 3, \, 0 \leq n < 2 \; .
\]
Such a word with $(m,n) \neq (0,0)$ does not fix a triple
$(1,\Omega,d)$ with $d \notin \{\pm 1, \pm \Omega\}$. A word
of shape $ x_d y_e x_\pi$, permutes the elements of 
each component $\mathcal{P}$ in $\mathcal{P}^3$. The 
stabilizers of components 1, 2 and 3 of $\mathcal{P}^3$ are
$\{x_{\pm 1}, x_{\pm\Omega}\}$, $\{y_{\pm 1}, y_{\pm\Omega}\}$
and $\{z_{\pm 1},  z_{\pm\Omega}\}$.  
The intersection of these stabilizers  is $\{1\}$.

\proofend

\stepcounter{Theorem}
An immediate consequence of Theorem~\ref{Thm:relations:N} is:
\begin{equation}
\label{consequence:Thm:relations:N}
   [x_d, x_e] =  x^{C(d,e)},  \quad x_d^2 =x^{P(d)}, \quad
   \nu_\delta^2 = [\nu_\delta, \nu_\epsilon] = 1\, . \quad \Cyclicxyz
\end{equation}
Note that $
[x_d, x_e] =(x_e x_d)^{-1} x_d x_e =
x_{A(d,e)} (x_ {ed})^{-1}x_{de}   x_{A(d,e)} 
%= x_{A(d,e)} x^{C(d,e)}   x_{A(d,e)}
= x^{C(d,e)}$ .

In order to make calculations in $\hat{N}$ and $N$ easy,  
we have stated more relations than necessary. From
Theorem~\ref{Thm:relations:N} and its proof we see that
$\hat{N}$ has structure $2^2.2^{12}.2^{24}.(S_3 \times M_{24})$.

Assigning ''odd'' parity to $x_\tau, y_\tau, z_\tau$ and to 
$\nu_\pi$, $\pi$ odd, and ''even'' parity to the other generators 
of $\hat{N}$, we see that the generators of $N$ have even parity 
and that the relations in  Theorem~\ref{Thm:relations:N} 
preserve that parity. Thus $N$ is the subgroup of the 
''even'' elements of $\hat{N}$ and we have $|\hat{N}/N|=2$.
Note also that $\tau \in N$.

Our generators of $N$ agree with Conway's generators in
\cite{Conway:Construct:Monster} with the following exceptions:
Our generators  $x_d, y_d, z_d$,  correspond to 
$x_d \cdot z_{-1}^{P(d)}$, $y_d \cdot x_{-1}^{P(d)}$
$z_d \cdot y_{-1}^{P(d)}$  in \cite{Conway:Construct:Monster}.
But our group $N$ agrees with that
in \cite{Conway:Construct:Monster}.
Our definition of $x_d, y_d, z_d$  leads to simpler relations 
in $N$, and it agrees with Ivanov's construction of $G_2$ in
\cite{citeulike:Monster:Majorana}, section~2.7,
with $G_2$ in \cite{citeulike:Monster:Majorana} corresponding
to Conway's and our group $N_0$.

As in \cite{Conway:Construct:Monster}, we define 
$\mathcal{K}_1 = y_\Omega z_{-\Omega}$, 
$\mathcal{K}_2 = z_\Omega x_{-\Omega}$, 
$\mathcal{K}_3 = x_\Omega y_{-\Omega}$, and we put:
\[
  K_1 = \{1,\mathcal{K}_1 \}, \quad
  K_2 = \{1,\mathcal{K}_2 \}, \quad
   K_3 = \{1,\mathcal{K}_3 \}, \quad
   K_0 =  \{1,\mathcal{K}_1,\mathcal{K}_2,\mathcal{K}_3 \}, \; .
\]
Then $K_0$ is normal in $N$ and $N_0 = N/K_0$ is the normalizer
of a certain fourgroup in $\mathbb{M}$. For a subgroup $\Gamma$
of $N$ we define $\Gamma_{\!n} = \Gamma /  \Gamma \cap K_n$,
$n= 0,1,2,3$. Note that $K_0$ is not normal in $\hat{N}$. 
The symbol $\Cyclicxyz$ will also indicate that a statement 
remains valid if the indices $1,2,3$ in the definitions above 
are cyclically permuted together with $x, y, z$. E.g.
\[
    \qquad\mathcal{K}_1=x_\Omega z \in K_1 \; \Cyclicxyz \qquad
    \mbox{implies}  \quad
      \mathcal{K}_2=y_\Omega x \in K_2 \quad \mbox{and} \quad
        \mathcal{K}_3=z_\Omega y \in K_3 \; .
\]

\stepcounter{Theorem}

As a corollary of the proof of Theorem~\ref{Thm:relations:N} 
we obtain the following structure of $N$:
\begin{equation}
\label{eqn:structure:N}
\begin{array}{ccccccccccccl}
N= & {2^2} & \! \! \times \!  \! & {2^2} &
     \! \!   .  \! \!  & 2^{11} &  \! \!  .  \! \!  &
    2^{22}&   \! \! . \! \!  & 
   (S_3  &  \! \! \times  \! \! & \,  M_{24})  \; . \\ 
&\scriptsize{\overbrace{\mathcal{K}_1,\mathcal{K}_2,\mathcal{K}_3}}&   
   & \scriptsize{\overbrace{x,y,z}}  &     
 & \scriptsize{\overbrace{x_\delta\!= \!y_\delta\!=\!z_\delta}} &
     
   & \scriptsize{\overbrace{x_d,y_d,z_d}}  &     
   & \scriptsize{\overbrace{x_\delta,y_\delta,z_\delta}}  &     
   & \scriptsize{\overbrace{x_\pi\!=\!y_\pi\!=\!z_\pi}}   \\
   &&&&& {\scriptstyle \delta} \, \mbox{\scriptsize even} &   
   &&&   {\scriptstyle \delta} \, \mbox{\scriptsize odd} &
   & {\scriptstyle \pi} \, \mbox{\scriptsize even}  
\end{array}
\end{equation}
By omitting the $M_{24}$  generators we obtain a normal
subgroup $Q$ and by omitting the $S_3$ generators we obtain a 
normal subgroup $N_{xyz}$ of $N$. 
We put $N_{xyz0} = N_{xyz} / K_0$.

The centralizer of $x$ in $N$ is called $N_x$ and has structure
\begin{equation}
\label{equ:structure:Nx}
\begin{array}{ccccccccccl}
N_x= & {2} & \times  & {2} & : & 2_+^{1+24} &
. & 2^{11} & .   & M_{24} \, & 
\quad   \Cyclicxyz .  \\ 
\scriptsize 
&\scriptsize{\overbrace{\mathcal{K}_1}}&  &
 \scriptsize{\overbrace{\mathcal{K}_2,\mathcal{K}_3}}&   &
 \scriptsize{\overbrace{\hphantom{mi}x_d, x_\delta\hphantom{mi}}}&     
& \scriptsize{\overbrace{y_d, z_d}}  &     
& \scriptsize{\overbrace{\hphantom{x} x_\pi \hphantom{x} }} \\
& & & & &  \mbox{\scriptsize (This is $Q_x$)}     
\end{array}
\end{equation}

Let   $Q_x$ be the subgroup of $N_x$ generated by 
$x_d, x_\delta$ with $d \in \mathcal{P}, \delta \in \mathbb{C}^*$. 
Then 
$Q_x \cap K_0 = 1$, so we have $Q_x \cong Q_{x1} \cong Q_{x0}$.
$Q_x$ is not normal in $N_x$, but $K_1 Q_x$ is. To see this,
it suffices to verify $[x_\delta,\mathcal{K}_1] = 1$ and
$\mathcal{K}_2 x_\delta = x_\delta \mathcal{K}_3$ for
odd $\delta$  (apart from trivial checks).

The defining relations in the group $Q_{x} $ given by
Theorem~\ref{Thm:relations:N} and
(\ref{consequence:Thm:relations:N}) simplify to:
\begin{align}
\label{relations:Qx}
x_{d}^2 = x^{P(d)}, \; x_{\delta}^2 = 1 , \; 
x_d x_e = x_{de} x_{A(d,e)}, \;
x_{\delta} x_{\epsilon} =   x_{\delta+\epsilon} , \;
[ x_d, x_\delta] = x^{\left<d,\delta\right>}  \; ,
\end{align}
so every element of $Q_{x}$ can uniquely be
written in the form $x_d x_\delta$. We often write $x_r,x_s,...$
for  elements of $Q_{x}$. We put
\begin{equation}
\label{eqn:def:tilde:xd}
 \tilde{x}_{d \vphantom{\tilde{d}}}
    = x_{-1}^{\scriptsize \mbox{sign}(d)}  x_{d} x_{\theta(d)}
 \quad \mbox{for} \;  d  \in \mathcal{P}, \;
 \mbox{and cocyle} \;  \theta, \; \mbox{i.e.} \;
    \theta(d) \in \mathcal{C}^* \, .
\end{equation}
$\tilde{x}_{d}$  does not depend on the sign of $d$, so
$\tilde{x}_d$ is also well defined for $d \in \mathcal{C}$.
Then the group $Q_{x}$ generated by $x_d, x_\delta$, 
$d \in  \mathcal{P}, \delta \in \mathcal{C}^*$ is also
generated by $\tilde{x}_d, x_\delta$. It is 
visibly extraspecial of type $2_+^{1+24}$, satisfying
the even simpler relations:
%\stepcounter{Theorem}
\begin{align}
\label{relations:Qx:extraspecial}
\tilde{x}_{d}^2 = x_{\delta}^2 = 1 , \; 
\tilde{x}_{d} \tilde{x}_{e} = \tilde{x}_{de} , \;
x_{\delta} x_{\epsilon}  =   x_{\delta\epsilon} , \;
[ \tilde{x}_{d}, x_{\delta} ] 
= x^{\left<\tilde{d},\delta\right>}  \; ,
\end{align}
which are easy consequences of (\ref{eqn:cocycle:Parker:Loop})
(\ref{Lemma:valid:cocycles}) and (\ref*{relations:Qx}).
We will discuss extraspecial 2-groups in 
section~\ref{section:extra:special:2}.
So the structure description 
(\ref{equ:structure:Nx}) of  $N_x$ should now be clear.

Similarly, we let $N_y$ and $N_z$ be the centralizers of 
$y$ and $z$ in $N$. We let  $Q_y$ and  $Q_z$ be the subgroups 
of  $N_y$ and $N_z$ generated by $y_d, y_\delta$ and
$z_d, z_\delta$, respectively.

Our definition of $Q_x$ differs from that in 
\cite{Conway:Construct:Monster}, since the definition of~$N$ 
in \cite{Conway:Construct:Monster} does not easily lead to a 
split extension corresponding to our extension $K_0:Q_x$.
We put 
\[
Q_{xyz} = K_0 Q_x \cap K_0 Q_y \cap  K_0 Q_z \; .
\]
$Q_{xyz}$ is elementary Abelian of structure $2^2.2^2.2^{11}$
and generated by $K_0, x ,y, z$ and $x_\delta$, $\delta$ even.

%%%%%%%%%%%%%%%%%%%%%%%%%%%%%%%%%%%%%%%%%%%%%%%%%%%%%%%%%%%%%%%%%%%%%%%
\section{The relation between $Q_{x}$ and the Leech lattice
	 $\Lambda$}
\label{section:Leech:link} 
%%%%%%%%%%%%%%%%%%%%%%%%%%%%%%%%%%%%%%%%%%%%%%%%%%%%%%%%%%%%%%%%%%%%%%%

\subsection{The homomorphism from $Q_x$ to the $\Lambda / 2\Lambda$}
\label{subsection:Qx:Leech}

The 24-dimensional Leech lattice $\Lambda$ is defined as follows. 
Consider an Euclidean vector space $\setR^{24}$ with scalar product 
$\left<.,.\right>$ and basis vectors labelled by the elements of
 $\tilde{\Omega}$. Assume that a Golay code $\mathcal{C}$ is given on  
$\Field_2^{24} = \prod_{i \in \tilde{\Omega}} \Field_2$. Then the 
Leech lattice is the set of vectors $u$ in $\setR^{24}$
with co-ordinates $(u_i, i \in \tilde{\Omega})$, such that 
there is an $m \in \{0,1\}$ and a $d \in \mathcal{C}$ with 
\[
\begin{array}{rllcl}
  \forall \, i \in \tilde{\Omega}  : & u_i
          &= &m + 2\cdot \left< d, i \right> &\pmod{4} \; ,\\
  \textstyle \sum_{i \in \tilde{\Omega}} & u
                 &=&  4m &\pmod{8}  \; . 
\end{array}
\]
See e.g. \cite{Conway-SPLG}, Chapter~4, section~11 for background. 
The Leech lattice $\Lambda$ is an even unimodular lattice, i.e. 
$\left<u, u\right>$ is even for all 
 $u \in \Lambda$ and $\det \Lambda = 1$. 
Therefore we have to scale  the basis vectors of the underlying space
$\setR^{24}$ so that they have norm $1/\sqrt{8}$, not $1$.  Then for
vectors $u, v$ with co-ordinates $u_i, v_i$ (where $i$ ranges over 
$\tilde{\Omega}$) we have $\left< u, v\right>$ = 
$\frac{1}{8}\sum u_i v_i$. For $u \in \Lambda$ we define
$\mbox{type}(u)$ = $\frac{1}{2}\left<u, u\right>$ 
=   $\frac{1}{16}\sum u_i^2$. 
The Leech lattice can be characterized as the (unique) 
24-dimensional even unimodular lattice without any vectors
of type~1, see e.g. \cite{Conway-SPLG}.

\newcommand{\InvSqrtEight}{ 
	 { \frac{\scriptstyle 1}{ \sqrt{\scriptstyle 8}}}}
 
\newcommand{\SmallTextMbox}[1]{\mbox{\rm{\scriptsize #1}}}

\begin{Theorem}
\label{Theorem:Q:Leech}	
There is an isomorphism from $Q_{x}/\{1, x\}$ to
$\Lambda/2 \Lambda$ given by:
\begin{align*}
  x_i  & \mapsto  \lambda_i =  ( -3_{\SmallTextMbox{on} \; i},
           1_{\SmallTextMbox{elswehere}}  )  \; ,\\
  x_d &  \mapsto  \lambda_d =  (  2_{\SmallTextMbox{on} \; d},
     0_{\SmallTextMbox{elswehere}}  )  \qquad 
        \mbox{\rm if} \; |d| = 0 \bmod{8} \; , \\
  x_d & \mapsto  \lambda_d =  (  0_{\SmallTextMbox{on} \; d},
    2_{\SmallTextMbox{elswehere}}  )  \qquad 
   \mbox{\rm if} \; |d| = 4 \bmod{8} \; , 
\end{align*} 
\[
  \mbox{which satisfies} \quad 
  [x_r, x_s]  =  x^{\left<\lambda_r, \lambda_s\right>} \, , \quad
  x_r^2 = x^{\SmallTextMbox{type}(\lambda_r)}  \, , \quad  
   \mbox{\rm with} \;\; \mbox{\rm type}(\lambda_r) = 
   {\textstyle   \frac{1}{2}  }
       \left<\lambda_r, \lambda_r\right> \; .
\]       
\end{Theorem}

Our homomorphism $Q_{x} \rightarrow \Lambda/2 \Lambda$
differs slightly from the corresponding homomorphism in
\cite{Conway:Construct:Monster}, Theorem~2, where $x_d$ is mapped
to $(2_{\SmallTextMbox{on} \; d}, 0_{\SmallTextMbox{elswehere}})$
for all $d \in \mathcal{P}$. But it agrees with the isomorphism 
$W_{24} \rightarrow  \Lambda/2 \Lambda$
in \cite{citeulike:Monster:Majorana}, Lemma~ 1.8.3,
which leads to a simpler set of relations in the group~$N$. 
Note that $W_{24}$ in \cite{citeulike:Monster:Majorana}
corresponds to $Q_x/\{1,x\}$ in this paper. 

\vspace{1.0ex}

\fatline{Proof of Theorem~\ref{Theorem:Q:Leech}}
  	
The proof is along the lines of the proof of  Theorem~2 in. 
\cite{Conway:Construct:Monster}.
The definitions of  $\lambda_i, \lambda_d$ imply:
\begin{align*}
x_{\Omega} & \mapsto -2 \lambda_i + \lambda_\Omega =
(8_{\SmallTextMbox{on} \; i}, 0_{\SmallTextMbox{elswehere}})
 \qquad  \mbox{for any} \; i \in \tilde{\Omega}   \; ,\\
  x_d x_\Omega^{|d|/4}  & \mapsto
 (2_{\SmallTextMbox{on} \; d}, 0_{\SmallTextMbox{elswehere}}) \; , \\
x_{\delta} x_{\Omega}^{|\delta|/2}   & \mapsto   (
4_{\SmallTextMbox{on} \; \delta},
0_{\SmallTextMbox{elswehere}}  ) \qquad
\mbox{for any} \; \delta \subset \tilde{\Omega}   \; \mbox{with} \; 
  |\delta| \; \mbox{even} \; .  
\end{align*}
For any $d \in \mathcal{C}$ the value  $|d|/2$ is even, so that
the product  $\Pi_{i\in d}x_i$  maps
to  $2 \cdot \lambda_d$ or to $2 \cdot \lambda_{d+\Omega}$.
Thus for any $\delta \in \mathcal{C}^*$ the image of 
$x_\delta$ is well defined in $\Lambda / 2\Lambda$, even if
$\delta$ is defined modulo $\mathcal{C}$ only. To show that the 
mapping $x_r \mapsto \lambda_r$ is a homomorphism,
it suffices to check that it preserves the relations  
\ref{relations:Qx}. Since $x_{-1}$ is mapped to 0, all these 
checks are easy except for the relation 
$x_d x_e = x_{de} x_{A(d,e)}$.
We have
\[
   x_d x_\Omega^{|d|/4} x_e  x_\Omega^{|e|/4} \mapsto
    (4_{\SmallTextMbox{on} \; d \cap e} , \;
      2_{\SmallTextMbox{on} \; d + e} , \;
         0_{\SmallTextMbox{elswehere}}  ) \; .
\]
Clearly, $x_{d+e} x_\Omega^{|d+e|/4} x_{d \cap e}
   x_\Omega^{|d\cap e|/2}$ maps to the same element of
$\Lambda/2\Lambda$. Thus the relation
$x_d x_e = x_{de} x_{A(d ,e)} x_\Omega^m$ is preserved for
$m = |d+e|/4 + |d|/4 + |e|/4 + |d \cap e|/2 = 0 \pmod{2}$. 
The mapping $x_r \mapsto \lambda_r$  is surjective by 
construction and both,  its origin and its image, have 
size $2^{24}$; hence it is an isomorphism.

In an extraspecial 2 group the commutator is bilinear and we
have $\mbox{type}(u+v) = \mbox{type}(u) + \mbox{type}(v)
   +\left<u,v\right>$ and $(x_rx_s)^2 = x_r^2 x_s^2[x_r,x_s]$.
So it suffices to check
$x_r^2 = x^{\mbox{\scriptsize type}(\lambda_r)}$ and
$ [x_r, x_s]  =  x^{\left<\lambda_r, \lambda_s\right>}$,
with $x_r$, $x_s$ running through all the generators 
 $x_d, x_\delta$ given in the Theorem. 
All these checks can be done by easy calculations 
using~\ref{relations:Qx}.  

\proofend

An immediate consequence of Theorem~\ref{Theorem:Q:Leech} is
$P(d) = \mbox{type}(\lambda_d) \pmod{2}$ for all
$d \in \mathcal{P}$.

We will abbreviate $x_d \cdot x_\delta$ to $x_{d \cdot \delta}$
and  $\lambda_d + \lambda_\delta$ to $\lambda_{d \cdot \delta}$
for $d \in \mathcal{P}, \delta \in \mathcal{C}^*$ as in
\cite{Conway:Construct:Monster}.
For $r = (d,\delta) \in \mathcal{P} \times \mathcal{C}^*$
we will define  $x_r = x_{d \cdot \delta}$, and for 
$f \in \{\pm 1, \pm\Omega\}$ we put 
$x_{f \cdot r} = x_f  \cdot x_r$. (The more relaxed condition 
$f \in \mathcal{P}$ could have the funny effect 
$x_{f \cdot d}  = x_f x_d \neq x_{fd})$.

\subsection{Short vectors in $\Lambda / 2\Lambda$}
\label{subcection:short:vectors}

 A vector in $\Lambda/2\Lambda$
is called {\em short}, if it is congruent modulo $2\Lambda$ to a
vector of type~2 in~$\Lambda$, i.e. to a shortest nonzero vector
in~$\Lambda$. $x_r$ is called short if its image $\lambda_r$
in  $\Lambda/2\Lambda$ is short. We define
$x_r^+$ to be $x_{\Omega \cdot r}$ if this is short and
$x_r$ otherwise.
The short vectors in $\Lambda / 2 \Lambda$ are given as follows, 
see
\cite{Conway:Construct:Monster} or \cite{Conway-SPLG}, Ch. 4.11.    
\begin{align}
 \lambda_{ij} &:
(4_{\SmallTextMbox{on}\,i}, -4_{\SmallTextMbox{on}\,j}, 
    0_{\SmallTextMbox{else}}),
\nonumber \\
\label{eqn:short:vectors}
 \lambda_{ij}^+ = \lambda_{\Omega \cdot ij} &:
   (4_{\SmallTextMbox{on}\,i,j},  0_{\SmallTextMbox{else}}),
   \\
\lambda^+_{d \cdot \delta} = \lambda_{\Omega^n d \cdot \delta} &:
 (-2_{\SmallTextMbox{on}\,\delta},
     2_{\SmallTextMbox{on}\,d \setminus \delta}, 
       0_{\SmallTextMbox{else}} ), \quad
         \; n = |\delta|/2, |d|=8, \mbox{and} \; \delta
         \; \mbox{is an even subset of} \; d,
        %  \subset d,   \delta \, \mbox{even},
 \nonumber   \\
\lambda^+_{d \cdot i}
= \lambda_{\Omega^m d\cdot i}  & :
  (\mp 3_{\SmallTextMbox{on}\,i}, \pm 1_{\SmallTextMbox{else}}),
   \; m={\left< d,i \right>+ P(d)},  \;
    \mbox{and the lower sign is taken on} \; d .
\nonumber      
\end{align}

For $d$ in $\mathcal{C}$ (or in $\mathcal{P})$ with 
$|d| \in \{8,16\}$ we define
$A(d,\mathcal{C})$ = $\{ A(d,c) \mid c \in \mathcal{C}\}$. 
In case  $|d|=8$ it follows from standard facts about the Golay
code that $A(d,\mathcal{C})$ is the set of all even elements of 
$\mathcal{C}^*$
which can be represented as a subset of $d$. An element of
$A(d,\mathcal{C})$  has exactly two representatives $\delta$ and
$d +\delta$ which are subsets of $d$. If
$\delta, \epsilon$ are given as elements of $A(d,\mathcal{C})$, 
then the expressions ${|\delta|/2}$ and 
${|\delta \cap \epsilon|}$ are well defined
(modulo~2)  under the assumption
that subsets of $d$ are chosen as representatives of $\delta$ and
$\epsilon$. Thus $\lambda^+_{d \cdot \delta}$, $|d|=8$, $\delta$
even, is short if and only if $\delta \in A(d,\mathcal{C})$. 
For each octad $d$ (i.e.  $d \in \mathcal{C}$, $|d|=8$) there 
are 64 short vectors of that shape. Note that
$A(d\Omega,\mathcal{C}) = A(d,\mathcal{C})$; so
${|\delta|/2}$ and 
${|\delta \cap \epsilon|}$ are also well defined for
$\delta, \epsilon \in  A(d,\mathcal{C})$, $|d|=16$.

We obviously have $\lambda_{-r}$ = 
$\lambda_{r}$ for all short vectors 
$\lambda_{r}$. We also have 
$\lambda^+_{d\cdot i}$ =  $\lambda^+_{\Omega d\cdot i}$.
$\mathcal{C}$ contains 759 octads. 
So the numbers of the short vectors  as given above are:
\[
  |\lambda_{ij}| =  |\lambda^+_{ij}| =
     {\textstyle \binom{24}{2} } = 276, \;
  |\lambda^+_{d \cdot \delta}|  = 759 \cdot 64 = 48576, \;
       % (\mbox{with} \, |d|=8), \;
     |\lambda^+_{d \cdot i}| = 2048 \cdot 24 = 49152 \, . 
\]
Altogether we have $98280$ short vectors in $\Lambda/2\Lambda$.

%%%%%%%%%%%%%%%%%%%%%%%%%%%%%%%%%%%%%%%%%%%%%%%%%%%%%%%%%%%%%%%%%%%%%%%
\section{Representations of the groups $N_{x0}$, $N_{y0}$, $N_{z0}$}
\label{section:rep:Nx1}
%%%%%%%%%%%%%%%%%%%%%%%%%%%%%%%%%%%%%%%%%%%%%%%%%%%%%%%%%%%%%%%%%%%%%%%

In this section we construct a 196884-dimensional 
real monomial representation $196884_x$ of the group
$N_{x0} = N_{0} \cap G_{x0}$. We give a basis of 
the vector space  $196884_x$  and we state the operations
of the generators of $N_{x0}$.
In section~\ref{section:rep:Nx0:N0} we extend  $196884_x$  to $N_0$
and in section~\ref{section_xi}  we extend $196884_x$  to $G_{x0}$
so that we eventually obtain a representation 
of the monster  $\mathbb{M}$. Similar representations
$196884_y$ and $196884_z$ of $N_{y0}$ and $N_{z0}$
may be obtained in the same way.

\subsection{Representations $98280_x$,  $98280_y$,  $98280_z$ of
          $N_{x0}$,  $N_{y0}$,  $N_{z0}$}

Since $Q_{x0} \cong Q_x$ and
 $Q_{x0}$ is normal in $N_{x0}$, the group $N_{x0}$ acts
on  $Q_{x}$  via conjugation, preserving the paring between 
$x_{r}$ and $x_{-r}$ in $Q_{x}$.
This leads to a representation of $N_{x0}$ on a real vector
space $98280_x$
 of dimension~98280 spanned by short vectors $X_r$,
$r \in \mathcal{P} \times \mathcal{C}^*$ with the relation
$-(X_r) = X_{-r}$, so that $N_{x0}$ acts on $X_r$ in the same
way as on $x_r$ by conjugation.
We will use the same symbol $98280_x$ for that
vector space and the action of the group $N_{x0}$ on that
vector space. Since $N_{x0}$ is a quotient of $N_{x1}$ and $N_x$,
this gives us also a representation on  $N_{x1}$ and $N_x$.

We give  $98280_x$  the structure of an Euclidean space
by assigning  norm $1$ to all vectors $X_r$,
and by declaring all pairs of such vectors perpendicular
unless they are equal or opposite. 

Considering $98280_x$ as a representation of $N_x$, its kernel
is the group generated by $K_0$ and $x$. The elements $y$ and $z$ 
of $N_x$ act on the basis vector $X_{d \cdot \delta}$ of $98280_x$ 
as $(-1)^{|\delta|}$. Using triality, we may define similar
representations $98280_y$ and  $98280_z$ of $N_y$ and $N_z$.

\subsection{Representations $4096_x$, $4096_y$, $4096_z$  of
    $N_{x1}$,  $N_{y2}$,  $N_{z3}$}
\label{subsection:rep:4096}

In this section we construct representations
$4096_x$, $4096_y$ and $4096_z$ of
the groups $N_{x1}$, $N_{y2}$ and $N_{z3}$. On the way we also
obtain a representation  ${(6 \! \cdot \! 2048)}\!_N$ of $N$ which 
will be useful for combining the first three representations to a
representation of $N_0$. 

We construct these representations from the action of $N$ on
$\mathcal{P}^3$ given in section~\ref{section:N0}. We augment
$\mathcal{P}$ by an element $0$ by decreeing 
$0 \cdot d = d \cdot 0 = 0 \cdot 0 = 0$, $d \in \mathcal{P}$ and 
we write $\mathcal{P}_0$ for $\mathcal{P} \cup 0$ with that
operation. Put 
$\mathcal{P}_0^{(3)} 
 = \{(d,0,0), (0,d,0), (0,0,d) \mid d \in \mathcal{P}\}$.
By definition $N$ acts as a permutation group on 
$\mathcal{P}_0^{(3)}$.
 The action of $N$ on $\mathcal{P}_0^{(3)}$
preserves the pairing between $(a,b,c)$ and $(-a,-b,-c)$. So we 
have a monomial representation of $N$ on the real vector space
${(6 \! \cdot \! 2048)}\!_N$
spanned by the basis vectors in  $\mathcal{P}_0^{(3)}$ with 
the identification $-(a,b,c) = (-a, -b, -c)$. Let 
$K_0$, $K_1$, and $K_2$ be as in section~\ref{section:N0}.
We prefer a different basis of 
 ${(6 \! \cdot \! 2048)}\!_N$, so that we can easily identify the
subspaces of  ${(6 \! \cdot \! 2048)}\!_N$ where the
representation of $K_0$, $K_1$, or $K_2$ is trivial.
For all $d \in \mathcal{P}$ we define
\begin{align*}
  d_1^- = (0,d,0) + (0,-\Omega d, 0) \, , \qquad
  d_1^+ = (0,0,\bar{d}) + (0,0,\Omega \bar{d}) \, , \\
  d_2^- = (0,0,d) + (0,0,-\Omega d) \, , \qquad
  d_2^+ = (\bar{d},0,0) + (\Omega \bar{d},0,0) \, , \\
  d_3^- = (d,0,0) + (-\Omega d,0,0) \, , \qquad
  d_3^+ = (0,\bar{d},0) + (0,\Omega \bar{d},0) \, .
\end{align*}

Here we make another deviation from Conway's original definition
in \cite{Conway:Construct:Monster}: $d_m^+$ in 
\cite{Conway:Construct:Monster} corresponds to our
$\bar{d}_m^+$. The motivation for this change is given in
section~\ref{section:Introduction}. This change  greatly 
simplifies the representation of the element $\xi$ of
$G_{x0} \setminus N_{x0}$ in $4096_x$, as stated in
Corollary \ref{eqn:transfrom_dpm_xi}. 

\Skip{
The reason for this modification is that the
representation  $4096_x$ defined below is actually a
representation of the real Clifford group $\mathcal{C}_{12}$,
as defined in \cite{NebRaiSlo2001}. There is a standard
construction of the $2^n$-dimensional matrix representation of
the Clifford group  $\mathcal{C}_{n}$ (and also of its complex
analogue), see e.g. \cite{NebRaiSlo2001}, which is 
also used in theory of quantum computing. Our matrix 
representation $4096_x$ is equivalent to the standard
construction of $\mathcal{C}_{12}$ if we define $d_m^+$ and 
$d_m^-$ as above. Using that standard construction greatly 
simplifies the representation of the element $\xi$ of
$G_{x0} \setminus N_{x0}$ in $4096_x$, as given by
Lemmas~\ref{eqn:transfrom_dpm_xigamma}
and~\ref{eqn:transfrom_dpm_xig}. 
}

The set $(d_m^\pm)$, $d \in \mathcal{P}$, $m=1,2,3$ is a basis of 
$(6 \! \cdot \! 2048)\!_N$ with the obvious identifications 
\stepcounter{Theorem}
\begin{align}
\label{eqn:d:Omega}
 (d)_m^- = -(-d)_m^-= -(\Omega d)_m^- =  (-\Omega d)_m^- \, ,  \;
 (d)_m^+ = -(-d)_m^+= (\Omega d)_m^+ =  -(-\Omega d)_m^+ \, . 
\end{align}
For $m=1,2,3$ we write $\left<d_m^+\right>$ and $\left<d_m^-\right>$ 
for the subspaces of  $(6 \! \cdot \! 2048)\!_N$ generated by the 
corresponding basis vectors. All these 2048-dimensional subspaces 
are invariant under the action of $N_{xyz}$. They are permuted
by the cosets of $N_{xyz}$ as indicated in  
Table~\ref{table:rep:6:2048:N} for the coset 
representatives~$x_i$ and $\tau$.

We give  $(6 \! \cdot \! 2048)\!_N$  the structure of an Euclidean space
by assigning the norm $1/\sqrt{2}$ to all vectors $(d,0,0)$, $(0,d,0)$,
and $(0,0,d)$ and by declaring all pairs of such vectors perpendicular
unless they are equal or opposite. Then all vectors $(d_m^\pm), m=1,2,3$
have norm~1 and pairs such vectors are also perpendicular
unless they are equal or opposite. 
In the sequel  $(6 \! \cdot \! 2048)\!_{N}$ means the representation
of $N$ on this Euclidean space with orthonormal basis vectors
taken from the vectors $(d_m^\pm), m=1,2,3$.
Then  $(6 \! \cdot \! 2048)\!_N$ is an orthogonal and monomial
representation of $N$. 

We define the representations $4096_x$, $4096_y$ and $4096_z$ of
$N_x$, $N_y$ and $N_z$, respectively, by their natural action on 
the following subspaces of $(6 \! \cdot \! 2048)\!_N$: 
\[
4096_x : \left<d_1^-\right>  \oplus \left<d_1^+\right> \, , \quad   
4096_y : \left<d_2^-\right>  \oplus \left<d_2^+\right> \, , \quad   
4096_z  :\left<d_3^-\right>  \oplus \left<d_3^+\right> \, ,  
\]
where  $\oplus$ denotes the direct sum. 
The kernels of $4096_x$, $4096_y$, $4096_z$,  and elements acting as
$-1$ are displayed in~Table~\ref{table:rep:6:2048:N}.
From these kernels we see that $4096_x$, $4096_y$, $4096_z$
are also representations of  $N_{x1}$, $N_{y2}$, $N_{z3}$,
respectively.

%\vspace{1ex}

\begin{table}[H]
\centering	
\begin{tabular}{|l|c|c|c|c|c|c|}
	\hline
      \small{Subspace of $(6 \! \cdot \! 2048)\!_{N} \!\!  $} 
    & \multicolumn{2}{c|}{$4096_x$}
    & \multicolumn{2}{c|}{$4096_y$}
    & \multicolumn{2}{c|}{$4096_z$} \\
	\cline{2-7}
  & $ \vphantom{A^{A^A}_g}
      \left<d_1^- \right> $ & $\left< d_1^+ \right> $
  & $ \left< d_2^- \right> $ & $\left< d_2^+ \right> $
  & $ \left< d_3^- \right> $ & $\left< d_3^+ \right> $\\
    \hline 
 	\small{Kernel of $4096_{x,y,z}$}    
   & \multicolumn{2}{c|}{$\vphantom{A^{A^A}_g}
   	1, \mathcal{K}_1, y_\Omega, z_{-\Omega}$}
   & \multicolumn{2}{c|}{$1, \mathcal{K}_2, z_\Omega, x_{-\Omega}$}
   & \multicolumn{2}{c|}{$1, \mathcal{K}_3, x_\Omega, y_{-\Omega}$}
    \\
 \hline 
 \small{Elements acting as -1}   
   & \multicolumn{2}{c|}{$\vphantom{A^{A^A}_g}
                x_{-1}, \mathcal{K}_2, \mathcal{K}_3 $}
   & \multicolumn{2}{c|}{$ y_{-1}, \mathcal{K}_1, \mathcal{K}_3 $}
   & \multicolumn{2}{c|}{$ z_{-1}, \mathcal{K}_1, \mathcal{K}_2 $}
    \\
  \hline 
  \small{$x_i$ maps subspace to}  
    & $ \vphantom{A^{A^A}_g}
    \left<d_1^+ \right> $ & $\left< d_1^- \right> $
    & $ \left<d_3^+ \right> $ & $\left< d_3^- \right> $
    & $ \left<d_2^+ \right> $ & $\left< d_2^- \right> $  \\
  \hline 
    \small{$\tau$ maps subspace to}  
  & $ \vphantom{A^{A^A}_g}
  \left< d_2^- \right> $ & $\left< d_2^+ \right> $
  & $ \left< d_3^- \right> $ & $\left< d_3^+ \right> $
  & $ \left<d_1^- \right> $ & $\left< d_1^+ \right> $  \\
  \hline 
\end{tabular}
\caption{The kernels of and the action of $S_3$ on the
	subspaces  of  ${(6 \! \cdot \! 2048)}\!_N$ }
\label{table:rep:6:2048:N}	
\end{table}

%\vspace{1ex}

\subsection{Representations $24_x, 24_y, 24_z$  of
    $N_{x1}$,  $N_{y2}$,  $N_{z3}$}

\newcommand{\GroupNAb}{{N_{\!\mbox{\scriptsize Ab}}}}

In this section we construct representations
$24_x$, $24_y$ and $24_z$ of the groups 
$N_{x1}$, $N_{y2}$ and $N_{z3}$.
On the way we also obtain a representation 
${(3 \! \cdot \! 24)}\!_N$ of $N$ which will be useful for 
combining the first three representations to a representation 
of $N_0$. 

\Skip{
The structure of $N$ given by (\ref{eqn:structure:N}) can also
be written as:
\begin{equation}
\label{eqn:structure:N_Ab}
\begin{array}{ccccccccccl}
N=  & {2^2} &
\! \!   .  \! \!  & 2^{11} &  \! \!  .  \! \!  &
2^{24}&   \! \! . \! \!  & 
(S_3  &  \! \! \times  \! \! & \,  M_{24})  \; . \\ 
& \scriptsize{\overbrace{x,y,z}}  &     
& \scriptsize{\overbrace{x_\delta\!= \!y_\delta\!=\!z_\delta}} &
& \scriptsize{\overbrace{x_d,y_d,z_d}}  &     
& \scriptsize{\overbrace{x_\delta,y_\delta,z_\delta}}  &     
& \scriptsize{\overbrace{x_\pi\!=\!y_\pi\!=\!z_\pi}}   \\
&&& {\scriptstyle \delta} \, \mbox{\scriptsize even} &   
&&&   {\scriptstyle \delta} \, \mbox{\scriptsize odd} &
& {\scriptstyle \pi} \, \mbox{\scriptsize even}  \\
\end{array}
\end{equation}
Let  $\GroupNAb$ be the Abelian normal subgroup of $N$ generated 
by  $x,y,z$ and $x_\delta=y_\delta=z_\delta$, $\delta$ even. 
In $N/\GroupNAb$  we have  $x_d = x_{-d} \; 
\Cyclicxyz$, and $N/\GroupNAb$ has structure:
\[
N/\GroupNAb = 2^{24} : (S_3 \times M_{24}), \quad 
\mbox{with $2^{24}$ elmentary Abelian and
      generated by $x_d, y_d, z_d$}.
\]
} % \Skip
Now we construct a representation   ${(3 \! \cdot \! 24)}\!_N$ of $N$.
For each $i \in \tilde{\Omega}$ we define $i_1$ by:
\[
i_1   =    \sum_{d \in \mathcal{P}}  
  (-1)^{ \left<d,i\right> }  (d,0,0) \, , 
\]
without any identification of $(d,0,0)$ and $(-d,0,0)$. 

We define $i_2$, $i_3$ similarly by using the triples
$(0,d,0)$ and $(0,0,d)$. Then it is easy to check that
the subgroup  $Q_{xyz}^{(-)}$ of $Q_{xyz}$ generated by 
$x_\delta, x, y, z$, $\delta$ even, preserves $i_1, i_2, i_3$
and that $N$ operates on $i_1$ as follows:  
\[
i_1  \stackrel{x_d}{\longmapsto}  i_1 \, ,  \quad
i_1 \stackrel{y_d, z_d}{\longmapsto}
(-1)^{\left<d,i\right>} i_1 \, ,  \quad
i_1 \stackrel{x_\pi}{\longmapsto} i_1^\pi \, , \quad 
i_1 \stackrel{\tau}{\longmapsto} i_{2} \, , \quad 
\Cyclicxyz \, .
\]
Here the symbol $\Cyclicxyz$ means that $x, y, z$ and
indices $1, 2, 3$  must be cyclically permuted.
We define $24_x$ to be the representation of $N_x$ with basis 
$i_1, i \in \tilde{\Omega}$.
Representations $24_z$ and $24_y$ are defined similarly, based on 
$i_2$ and $i_3$. Then
\[
{(3 \! \cdot \! 24)}\!_N = 24_x \oplus 24_y \oplus 24_z 
\]
is a representation of $N$ with kernel $Q_{xyz}^{(-)}$.
$x_\delta$ preserves $i_1$ and exchanges $i_2$ with $i_3$ for
odd $\delta$, and  $\tau$ cyclically permutes
$i_1$, $i_2$ and $i_3$.

We declare $\{i_m | i \in \tilde \Omega, m = 1,2,3\}$ to be an 
orthonormal basis of the Euclidean space ${(3 \! \cdot \! 24)}\!_N$.
Then $N$ acts orthogonally and monomially on that space.

Representations $24_x, 24_y, 24_z$ are invariant under the action
of $N_{xyz}$. But they are permuted by the cosets of $N_{xyz}$ 
in $N$, as indicated by the action of the coset representatives 
$x_i$ and $\tau$ in the following Table~\ref{table:rep:3:24:N}. 
This table also displays the
kernels of these representations and the elements acting as $-1$.

%\vspace{1ex}

\begin{table}[H]
\centering	
\begin{tabular}{|l|c|c|c|l}
 \hline
 \small{Subspace of $(3 \! \cdot \! 24)\!_{N} \!\!  $} 
 & $   24_x $ 
 & $  24_y $ 
 & $  24_z $\\
 \hline
 \small{basis vectors} &  
 $ \vphantom{O^{O^O}}   i_1, i \in \tilde{\Omega}$
   & $i_2$ & $i_3$  \\
 \hline 
 \small{Kernel generated by}    
 & $\! \mathcal{K}_1, x_d, x_\delta,y,z$  
 & $\! \mathcal{K}_2, y_d, y_\delta,x,z$  
 & $\! \mathcal{K}_3, z_d, z_\delta,x,y$  \\
 \hline 
 \small{Elements acting as -1}   
 & { $  \mathcal{K}_2, \mathcal{K}_3, y_\Omega, z_\Omega $}
 & { $  \mathcal{K}_1, \mathcal{K}_3, x_\Omega, z_\Omega $}
 & { $  \mathcal{K}_1, \mathcal{K}_2, x_\Omega, y_\Omega $ } \\
 \hline 
 \small{$x_i$ maps subspace to}  
 & $ 24_x $ & $ 24_z $ & $ 24_y $ \\
 \hline 
 \small{$\tau$ maps subspace to}  
 & $ 24_y $ & $ 24_z $ & $ 24_x $ \\
 \hline 
\end{tabular}
\caption{Kernels of and action of $S_3$ on the
	subspaces  of  ${(3 \! \cdot \! 24)}\!_N$ }  
\label{table:rep:3:24:N}	
\end{table}

\Skip{
Occasionally we write $\left<i_1\right>$, $\left<i_2\right>$,
$\left<i_3\right>$ for $24_x, 24_y, 24_z$, respectively. So
$\tau$ maps $\left<i_m\right>$ to $\left<i_{m+1} \right>$, with
$m$ to be taken modulo~3.   
}

%\vspace{1ex}

\subsection{Representation $196884_x$,  $196884_y$,  $196884_z$  of
    $N_{x0}$,  $N_{y0}$,  $N_{z0}$}
\label{subsection:196884}

We will now construct a representation $196884_x$ of subgroup 
$N_{x0} = N_x/K_0$ of the monster  $\mathbb{M}$ from the 
representations $24_x$, $4096_x$ and  $98280_x$. 
We define the representation
\[
    196884_x = 300_x \oplus 98280_x \oplus 98304_x
\] 
of $N_{x0}$ to be the direct sum of representations
$300_x$, $98280_x$ $98304_x$, where
\begin{align*}
    300_x & \mbox{ is the symmetric tensor square 
  $24_x \otimes_{\mbox{\tiny sym}} 24_x$  of } \, 24_x \; , \\
  98280_x & \mbox{ is defined as above} \, , \\
  98304_x & \mbox{ is the tensor product } \; 4096_x \otimes 24_x \, .
\end{align*}  
$98280_x$ is a representation of $N_{x0}$ by construction. For both
representations, $300_x$ and $98304_x$, the group $K_1$ is in their
kernel, and $\mathcal{K}_2, \mathcal{K}_3$ act as $-1$.
Thus $K_0$ is contained in the kernel of the tensor products
$24_x \otimes 24_x$ and $4096_x \otimes 24_x$, 	so that
$300_x$ and $ 98304_x$ are indeed representations of $N_{x0}$. 

We introduce some abbreviations for the basis vectors that we will
use for $196884_x$: 
\begin{align*}
\mbox{for } \hphantom{00}300_x: & \quad
    (ii)_1 = i_1 \otimes i_1  \, , \qquad  
        \mbox{of norm 1}, \\    
     & \quad
    (ij)_1 = i_1 \otimes j_1 + j_1 \otimes i_1 \quad ( i \neq j)
     \, ,   \;   \mbox{of squared norm 2}, \\    
\mbox{for } 98280_x: & \quad X_r \quad (r \;  \mbox{short})
    \, ,   \;   \mbox{of norm 1}, \\ 
\mbox{for } 98304_x: & \quad
     d^\pm \otimes_1 i =   d_1^\pm \otimes i_1 \quad
       (d \in \mathcal{P}, i \in \tilde{\Omega}  ) 
          \, ,   \;   \mbox{of norm 1}. 
\end{align*}
These basis vectors are mutually orthogonal, except when equal or
opposite. For representations  $196884_y$ and  $196884_z$ we use
the corresponding notation $\Cyclicxyz$
obtained by cyclically permuting $(x,y,z)$ and $(1,2,3)$.
The action of the generators of $N$ on the basis vectors
of $196884_x$, $196884_y$ and $196884_z$ can be obtained from
Table~\ref{table:action:196884x}.

\vspace{1ex}

\newcommand{\VPhantomTA}{$\vphantom{d_{g_{g_g}}^{P^{\bar{P}}}}$}

\newcommand{\XYZplus}[2]{#1^{\vphantom{g_g}+}_{\vphantom{1^1}#2}}

\begin{table}[H]
\centering	
\begin{tabular}{|l|c|c|c|c|c|}
	\hline
 $g$ &  	 \multicolumn{5}{l|}{
       \small Action of $g$ on the  basis vectors} \\
  \cline{2-6}
   &  $i_1$ \VPhantomTA &
        $d_1^-$ & $d_1^+  $ &
        $ X_{d \cdot \delta} \, ,\;  \delta $ \small even  & 
        $ X^+_{d \cdot i}\, , \; i \in \tilde{\Omega}$ \\
      \hline
  $x_e$  \VPhantomTA
    &  $i_1$   &   $(\bar{e} d)_1^-$    & $(\bar{e} d)_1^+$ &
        $(-1)^{n} X_{d \cdot \delta}$ 
     &
         $(-1)^{\left<e, i \right>} 
               X^+_{\bar{e} d e \cdot i}$   \\
     \hline 
   $y_e$  \VPhantomTA
    &  $ (-1)^{\left<e,i\right>} i_1$
    &  $ (\bar{e} d e)_1^- $  &  $ (de)_1^+ $  
    &  $(-1)^{\left<e,\delta \right>}
             X_{\Omega^n d \cdot \delta \delta'}$  
     &    $X^+_{de \cdot i} $ \\
     \hline 
   $z_e$  \VPhantomTA
    &  $ (-1)^{\left<e,i\right>} i_1 $
    &  $ ( d e)_1^- $  &  $ (\bar{e} d e)_1^+ $  
    & $(-1)^{C(d,e)}
       X_{\Omega^n d \cdot \delta \delta'}$   
    & $(-1)^{\left<e,i\right>} X^+_{ \bar{e}d \cdot i}$ \\
    \hline 
   $x_\pi , \;
          \pi $  \small even\VPhantomTA
          & $ (i^\pi)_1$  & $(d^\pi)_1^-$ &  $(d^\pi)_1^+$
    &   $ X_{d^\pi \cdot \delta^\pi} $ 
    &   $ X^+_{d^\pi \cdot i^\pi} $ 
    \\
    \hline 
     $x_\pi, \; \pi $ \small odd \VPhantomTA
    &  $ (i^\pi)_1$  & $({d}^\pi)_1^+$ &  $({d}^\pi)_1^-$ 
    &  $  X_{d^\pi \cdot \delta^\pi} $
    &  $ (-1)^{m}
         \XYZplus{X}{d^\pi \cdot i^\pi} $     
    \\   
    \hline
     $y_\pi, \; \pi $ \small odd \VPhantomTA
    &  $ (i^\pi)_3$  & $({d}^\pi)_3^+$ &  $({d}^\pi)_3^-$ 
    &  $  Z_{d^\pi \cdot \delta^\pi} $
    &  $ (-1)^{m}
         \XYZplus{Z}{d^\pi \cdot i^\pi} $     
    \\   
    \hline
     $z_\pi, \; \pi $ \small odd \VPhantomTA
    &  $ (i^\pi)_2$  & $({d}^\pi)_2^+$ &  $({d}^\pi)_2^-$ 
    &  $  Y_{d^\pi \cdot \delta^\pi} $
    &  $(-1)^{m}
           \XYZplus{Y}{d^\pi \cdot i^\pi} $  
    \\
    \hline
     $ \tau $ \VPhantomTA & $ i_2 $ & $ d_2^- $ & $ d_2^+ $ 
     & $ Y_{d \cdot \delta}  $ & $ Y^+_{d \cdot i} $  
    \\   
    \hline
     $ \tau^2 $ \VPhantomTA & $ i_3 $ & $ d_3^- $ & $ d_3^+ $ 
     & $ Z_{d \cdot \delta}  $ & $ Z^+_{d \cdot i} $  
    \\   
\hline
\end{tabular}
\captionsetup{justification=centering}
\caption{Action of the generators of $N$ on 
 $24_x$, $4096_x$ and $ 98280_x $. \newline
  \small Notation: $n = C(d,e)+ \left<e,\delta\right> , \;
\delta' = A(d,e) \, , \;
X^+_{d \cdot i} 
      = X^{}_{\Omega^{m} d \cdot i} \, , \;
  \mbox{\small with} \; m = P(d) +\left<d, i \right> \; .$ 
}  
\label{table:action:196884x}	   
\end{table}

The action of the generating elements $g$ of $N$ on the basis 
vectors $i_2, d_2^\pm, Y_{d \times \delta}$ and
$i_3, d_3^\pm, Z_{d \times \delta}$   can by obtained by
applying the triality operation $\Cyclicxyz$ to the
corresponding entries in Table~\ref{table:action:196884x}.
E.g. from the action of $z_e$ on $d_1^+$ and $i_1$  we may deduce
\[
d_2^+ \otimes i^{\vphantom{+}}_2 \stackrel{x_e}{\longmapsto} 
 (-1)^{\left<e,i\right>} (\bar{e}de)_2^+ \otimes i^{\vphantom{+}}_2 
     \qquad \mbox{and} \qquad
d_3^+ \otimes i^{\vphantom{+}}_3 \stackrel{y_e}{\longmapsto} 
(-1)^{\left<e,i\right>} (\bar{e}de)_3^+ \otimes  i^{\vphantom{+}}_3 
\]
by applying the triality operation  $\Cyclicxyz$. 
In the sequel the phrase ''from Table \ref{table:action:196884x}
we obtain \ldots'' means that the reader has to apply the triality
operation $\Cyclicxyz$ by himself, if necessary.

The action of $g$ on the basis vectors 
$X_{d \cdot \delta}, X_{d \cdot i}$ is obtained by conjugation 
of the corresponding basis vector with $g$ and taking
the result modulo the kernel $K_1$. 

The basis vector
$X^+_{d \cdot i}$  in Table~\ref{table:action:196884x}
is short for all $ d \in  \mathcal{P}, i \in \tilde{\Omega}$
and all short basis vectors $X_{d \cdot \delta}$ for odd
$\delta$ are of that shape.  We have $X^+_{d \cdot i}$ =
$X^+_{\Omega d \cdot i}$. The notation $X^+_{d \cdot i}$
in Table~\ref{table:action:196884x}  is consistent
with the notation $x^+_{ d \cdot i}$ 
in section~\ref{subcection:short:vectors}.

The following table is helpful for computing the last two
columns in Table~\ref{table:action:196884x}:

\begin{table}[H]
    \centering	
    \begin{tabular}{|l|c|c|c|c|}
        \hline
        $g$ &  	 \multicolumn{4}{l|}{
            \small  Value of $\; g^{-1} B g \;$ for } \VPhantomTA  \\
        \cline{2-5}
        &  $ B = x_d $ \VPhantomTA  
        &  $ B = x_\delta, \delta $  \small even
        & $ B = x_i, i \in \tilde{\Omega} $
        & Remarks \\ 
        \hline
        $x_e$ \VPhantomTA &  $x^{C(d,e)} x_d $   
        & $ x^{\left<e,\delta\right>} x_\delta  $
        &  $x^{\left<e,i\right>} x_i$ & \\
        \hline       
        $y_e $ \VPhantomTA  &  $  z^{C(d,e)} x_d x_{A(d,e)} $  & 
        $ y^{\left<e,\delta\right>} x_\delta  $ 
        & $z^{P(e) + \left<e, i \right>} x_e x_i$
        & $ y = x_{-\Omega}, z  = x_\Omega$ \\
        \cline{1-4}     
        $z_e $ \VPhantomTA  &  $  y^{C(d,e)} x_d x_{A(d,e)} $  
        & $ z^{\left<e,\delta\right>} x_\delta  $
        & $y^{P(e) + \left<e, i \right>} x_e x_i$
        & $  \pmod{K_1}$ \\
        \hline
        $x_\pi$ &   $x_{d^\pi}$ &   $x_{\delta^\pi}$
        &   $x_{i^\pi}$ & \\
        \hline     
        $y_\pi $ &  $ z_{d^\pi} $   &$x_{\delta^\pi}$  
        & $ z_{i^\pi} $ &  $\pi$ \small odd\\
        \hline      
        $z_\pi $ &  $ y_{d^\pi} $   &$x_{\delta^\pi}$  
        & $ y_{i^\pi} $ & $\pi$  \small odd \\
        \hline      
    \end{tabular}
    \captionsetup{justification=centering}
    \caption{Action of generators of $N$ by conjugation 
    }  
    \label{table:action:196884x:aux}	   
\end{table}

%%%%%%%%%%%%%%%%%%%%%%%%%%%%%%%%%%%%%%%%%%%%%%%%%%%%%%%%%%%%%%%%%%%%%%%
\section{Extending the representation $196884_x$ from
	$N_{x0}$ to $N_0$}
\label{section:rep:Nx0:N0}
%%%%%%%%%%%%%%%%%%%%%%%%%%%%%%%%%%%%%%%%%%%%%%%%%%%%%%%%%%%%%%%%%%%%%%%

\subsection{The action of the triality element $\tau$ on $196884_x$}
\label{subsection:action:tau}

In order to extend the representation  $196884_x$ from
$N_{x0}$ to $N_0$ it suffices to state the action of the triality
element $\tau$ on the basis vectors of $196884_x$. We first state the
result, which is the most important information for a
programmer. Proofs are deferred to the next two subsections. 

The operation of $\tau$ on the basis vectors
$(ij)_1$, $X_{ij}$ and  $X^+_{ij}$ is given by:
\stepcounter{Theorem}
\begin{align}
\label{triality:ii}
(ii)_1  \, & \stackrel{\tau}{\longmapsto} \, (ii)_1 \, , \\
\label{triality:ij}
(ij)_1  \, & \stackrel{\tau}{\longmapsto} \,
X_{ij} - X_{ij}^+   \, \stackrel{\tau}{\longmapsto} \,
X_{ij} + X_{ij}^+   \, \stackrel{\tau}{\longmapsto} \, (ij)_1  \, ,
\quad i \neq j  \; .
\end{align}
The operation of $\tau$ on the basis vectors
$X^+_{d\cdot i} $ and $ d^\pm \otimes_1 i$ is given by:
\begin{align}
\label{triality:Xdi}
  X^+_{d \cdot i} 
  = X_{\Omega^{\left<d,i \right> +P(d)}d,i}  
  \, \stackrel{\tau}{\longmapsto} \,
      (-1)^{\left<d,i \right>} \cdot d^- \otimes_1 i
  \, \stackrel{\tau}{\longmapsto} \,
     \bar{d}^+ \otimes_1 i
  \, \stackrel{\tau}{\longmapsto} \,
       X^+_{d\cdot i} \, . 
\end{align}

Now we specify the action of  $\tau$ on the remaining basis vectors
$X_{d \cdot \delta}$ of $196884_x$. We put $X^+_{d \cdot \delta}$ =
$X_{\Omega^n d \cdot \delta}$ for $n = |\delta|/2$, 
 $|d|=8$, $\delta \in A(d,\mathcal{C})$,
with $ A(d,\mathcal{C})$  = $\{ A(d,c) \mid c \in \mathcal{C}\}$,
in accordance with the notation in  
section~\ref{subcection:short:vectors}. Clearly all remaining
basis vectors are of that shape. Let $V_T$ be the subspace of
$196884_x$ generated by these basis  vectors. We extend $V_T$ from 
a representation of $N_{x}$, with kernel $K_T$ generated by 
$x_f,y_f,z_f$ for $f \in \{-1, \pm \Omega\}$, to a representation 
of $\hat{N}$. Then $V_T$ is also
a representation of $N$. It is also a representation of $N_0$,
because $K_0$ is in the kernel of  $V_T$. 

The operations of the generators 
$x_e, y_e, z_e, \nu_\pi, x_\tau, y_\tau$
of $\hat{N}$ on the basis vector $X^+_{d \cdot \delta}$ of $V_T$
are given by:
\begin{equation}
\label{eqn:operations:V_T}
\begin{array}{llll}
   x_e: & 
(-1)^{C(d,e) + \left<e,\delta\right>}   X^+_{d \cdot \delta}
    \,,  \qquad \quad &
   y_e:  & 
(-1)^{ \left<e,\delta\right>}    X^+_{d \cdot \delta  \delta'}
    \, ,\qquad \mbox{with} \quad \delta' = A(d,e) \, ,\\
  z_e:  & 
  (-1)^{C(d,e)}   X^+_{d \cdot \delta \delta'}  
       \,,  \qquad &
  \nu_\pi: &   (-1)^{|\delta|\cdot |\pi|/2}
X^+_{d^\pi \cdot \delta^\pi} \, , \\
   x_\tau: &
   (-1)^{|\delta|/2}  X^+_{d \cdot \delta}
       \,,  \qquad &
   y_\tau: &
          {\textstyle \frac{1}{8}} \,
   { \sum_{\epsilon \in A(d,\mathcal{C}) } \;
       (-1)^{|\delta \cap \epsilon|} X^+_{d \cdot \epsilon}}  \, .    
\end{array}
\end{equation}
Here $|\delta|/2$ and $|\delta \cap \epsilon|$ are defined
(modulo 2) for $\delta$ and $\epsilon$ as even subsets of $d$ 
as in section~\ref{subcection:short:vectors}, and
$|\pi|$ is equal to 1 for odd and to 0 for even $\pi$.
The sum in the expression for $y_\tau$ runs over all even subsets 
$\epsilon$ of $d$, identifying $d \setminus \epsilon$ with
$\epsilon$, so that it has 64 terms.     
$x_\tau$ is monomial and the action of $y_\tau$ resembles that of a 
Hadamard matrix, thus allowing a very efficient implementation.

The operations of $x_e, y_e, z_e$ and $x_\pi$
are obviously equal to the corresponding operations in 
Table~\ref{table:action:196884x}. 
In section~\ref{Subsection:proof:V_T} we will show that the
remaining operations in $V_T$ a are consistent with the
relations in $\hat{N}$. By Theorem~\ref{Thm:relations:N}
we have 
$\tau = y_\tau x_\tau$ and $\tau^2 = \tau^{-1} = x_\tau y_\tau$.
We obtain the operation of $x_\pi$ in 
Table~\ref{table:action:196884x} by using 
$x_\pi = x_\tau^{|\pi|} \nu_\pi$.

\vspace{1em}

The operation of the triality element $\tau$ leads to an identification
of the three spaces $196884_x$,  $196884_y$ and  $196884_z$,
which is called the {\em dictionary} in 
\cite{Conway:Construct:Monster},  Table~2. In our notation
we have the following identification of these three spaces:

\begin{table}[H]  
    \centering	
    \begin{tabular}{clccccc}
    	Subspace  & \quad & $198664_x$ &&  $198664_y$ &&  $198664_z$
    	 \vspace{5pt} \\
    	$V_A$ &&
        $(ij)_1$  & = & $Y_{ij} + Y_{ij}^+$ & = & $Z_{ij} - Z_{ij}^+$
        \\  
        $V_B$ && 
        $X_{ij} - X_{ij}^+$ & = & $(ij)_2$  & = & $Z_{ij} + Z_{ij}^+$
        \\
        $V_C$ &&
        $X_{ij} + X_{ij}^+$ & = &   $Y_{ij} - Y_{ij}^+$ & = & $(ij)_3$
        \\
        $V_D$ &&
        $(ii)_1$  & = &   $(ii)_2$    & = &    $(ii)_3$
        \\
        $V_X$ &&
        $X^+_{d \cdot i}$  & = &  $ \bar{d}^+ \otimes_2 i $
        & = & $ (-1)^{\left<d,i\right>} d^- \otimes_3 i $ 
        \\  
        $V_Y$ &&   
        $ (-1)^{\left<d,i\right>} d^- \otimes_1 i $  & = &
        $Y^+_{d \cdot i}$  & = &  $ \bar{d}^+ \otimes_3 i $   
        \\
        $V_Z$ &&
        $ \bar{d}^+ \otimes_1 i $  & = &        
        $ (-1)^{\left<d,i\right>} d^- \otimes_2 i $  & = &   
        $Z^+_{d \cdot i}$
        \\
        $V_T$ &&
        $X^+_{d\cdot\delta}$& = &
        $\textstyle \frac{1}{8}\sum (-)_{\delta,\epsilon}  
            Y^+_{d\cdot\epsilon}$  & = &
        $\textstyle \frac{1}{8}\sum (-)_{\epsilon,\delta}  
               Z^+_{d\cdot\epsilon}$ 
        \\  
        $V_T$ &&   
        $\textstyle \frac{1}{8}\sum (-)_{\epsilon,\delta}  
        X^+_{d\cdot\epsilon}$ & = &
        $Y^+_{d\cdot\delta}$& = &
        $\textstyle \frac{1}{8}\sum (-)_{\delta,\epsilon}  
        Z^+_{d\cdot\epsilon}$  
        \\
        $V_T$ &&
        $\textstyle \frac{1}{8}\sum (-)_{\delta,\epsilon}  
        X^+_{d\cdot\epsilon}$ & = &  
        $\textstyle \frac{1}{8}\sum (-)_{\epsilon,\delta}  
        Y^+_{d\cdot\epsilon}$ & = &
        $Z^+_{d\cdot\delta}$ 
    \end{tabular} 
    \captionsetup{justification=centering}
    \caption{The dictionary \\
    \small{Here    $(-)_{\delta,\epsilon} = 
        (-1)^{|\delta\cap\epsilon|+|\delta|/2 }$ 
         and the sum in the last three lines runs over
          $\epsilon \in A(d,\mathcal{C})$.}}        
    \label{table:dictionary}	   
\end{table}

Following \cite{Conway:Construct:Monster}, we also give names to the
subspaces of $196884_x$ as indicated in table~\ref{table:dictionary}. 

\subsection{Proofs for the monomial and almost-monomial actions of 
  $\tau$}

The purpose of this subsection is to establish
(\ref{triality:ii}), (\ref{triality:ij}) and (\ref{triality:Xdi}).

Using the information in Table~\ref*{table:action:196884x}, it is 
easy to check that the identification $(ii)_1 = (ii)_2 = (ii)_3$ 
in line~4 of Table~\ref{table:dictionary}
is compatible with the action of $N_{x0}$.
Together with the obvious action
$ (ii)_1 \stackrel{\tau}{\mapsto}  (ii)_2 \stackrel{\tau}{\mapsto}
(ii)_3 \stackrel{\tau}{\mapsto}  (ii)_1$ we obtain
(\ref{triality:ii}). 

Next we show (\ref{triality:ij}). 
The group ring $\setR Q_{xyz0}$ is the real algebra with basis
vectors labelled by $Q_{xyz0}$, and multiplication of the
basis vectors given by the group operation in $Q_{xyz0}$. 
Since $Q_{xyz0}$ is normal in $N_0$, the group $N_0$
operates on   $\setR Q_{xyz0}$ by conjugation.

Let $V$ be the subspace of 
$98280_x \oplus 98280_y \oplus 98280_z$ spanned by
$X_{f \cdot \delta}, Y_{f \cdot \delta}, Z_{f \cdot \delta}$, 
$f \in \{\pm1, \pm \Omega\}$, $\delta$ even. We define
a linear mapping
$V \rightarrow \setR Q_{xyz0}$ by:
\[
    X_{f \cdot \delta} \rightarrow
        x_{f \cdot \delta} - x_{-f \cdot \delta} \; ,\quad  
    Y_{f \cdot \delta} \rightarrow
       y_{f \cdot \delta} - y_{-f \cdot \delta} \; ,\quad  
    Z_{f \cdot \delta} \rightarrow 
       z_{f \cdot \delta} - z_{-f \cdot \delta} \; . 
\] 
The restriction of this mapping to $98280_x$, to $98280_y$ and to
$98280_z$ is injective and  the mapping preserves the
operation of $N_0$. So we may identify the basis vectors of
$V$  with their images in $\setR Q_{xyz0}$.
Using the relations in 
$ Q_{xyz0} =  Q_{xyz}/K_0$ given by $K_0$ and 
$x_{(ij)}=y_{(ij)}=z_{(ij)}$  we obtain:
\[
\begin{array}{ccccccc}
   Y_{ij} +  Y^+_{ij} =  Y_{ij} +  Y_{\Omega \cdot ij}  & =
   &  y_{1\cdot (ij)}    + y_{\Omega\cdot (ij)}
  \, - \,  y_{-1\cdot (ij)}   -   y_{-\Omega\cdot (ij)}  & \\
& =         
   &  z_{1\cdot (ij)}  + z_{-\Omega\cdot (ij)}
   \, - \,   z_{\Omega\cdot (ij)}  -   z_{-1 \cdot (ij)}  
 & = &    Z_{ij} -  Z^+_{ij} \;  \\
 & = &   x_{1\cdot (ij)}  +  x_{-1\cdot (ij)}  
    -   x_{-\Omega \cdot (ij)}  -  x_{\Omega \cdot (ij)}   & .      
\end{array}
\]
Hence we may identify $ Y_{ij} + Y^+_{ij}$ with $ Z_{ij} - Z^+_{ij} $.
A similar calculation yields  $Z_{ij} + Z^+_{ij} = X_{ij} - X^+_{ij}$
and $X_{ij} + X^+_{ij} = Y_{ij} - Y^+_{ij}$.
Put $(X_{ij})$ = $  x_{1\cdot (ij)}  + x_{-1\cdot (ij)}  
- x_{\Omega \cdot (ij)}  - x_{-\Omega \cdot (ij)} $, and
define $(Y_{ij})$, $(Z_{ij})$  similarly, using triality.
From Table~\ref{table:action:196884x}
we see that the action of the generators $x_d$, $x_\pi$ and $\tau$ of 
$N_{x0}$ on the vectors $(X_{ij})$,  $(Y_{ij})$ and  $(Z_{ij})$ 
is the same as their action on $(ij)_1$,  $(ij)_2$ and
$(ij)_3$, respectively. So we have just established the identifications
in the first three lines of  Table~\ref{table:dictionary}. Using these
identifications and the obvious action
$ (ij)_1 \stackrel{\tau}{\mapsto}  (ij)_2 \stackrel{\tau}{\mapsto}
(ij)_3 \stackrel{\tau}{\mapsto}  (ij)_1$ we obtain
(\ref{triality:ij}).

Next we show (\ref{triality:Xdi}). This is a consequence of
lines 5--7 of the dictionary in Table~\ref{table:dictionary}
and of the obvious operation of the triality element $\tau$.
Lines  5--7 in Table~\ref{table:dictionary} are direct consequences 
of the following two lemmas.

\begin{Lemma}
    \label{Lemma:rep:3:49152}	
    The group $N_0$ has a  monomial representation
    $(3 \! \cdot \! 49152 )\!_N$ with basis vectors
    $d^\pm \otimes_m i$, $d \in \mathcal{P}$, $i \in \tilde{\Omega}$,
    $m=1,2,3$, and the
    identifications
    \[
    (-d)^\pm =  -(d^\pm), \quad (\pm \Omega d)^\pm = d^\pm, \quad
    \bar{d}^+ \otimes_m i =
    (-1)^{\left<d,i\right> }  d^- \otimes_{m+1} i \; ,
    \]
    where $m$ is to be taken modulo~3, and the action of the generators
    of $N_0$ is given by table~\ref{table:action:196884x} and triality.   
\end{Lemma}

\fatline{Proof}

For $m=1,2,3$ let $\left<d^+ \otimes_m i \right>$ and 
$\left<d^- \otimes_m i \right>$ be the vector spaces spanned
by  the basis vectors  $d^+ \otimes_m i$ and  
$d^- \otimes_m i$, respectively, for 
$d \in \mathcal{P}$, $i \in \tilde{\Omega}$,
with the identifications given by (\ref{eqn:d:Omega}).
Define ${(6 \! \cdot \! 49152 )}\!_N$ by:
\[
{(6 \! \cdot \! 49152 )}\!_N = 
\left<d^+ \otimes_1 i \right> \oplus \left<d^- \otimes_1 i
\right> 
\, \oplus \, 
\left<d^+ \otimes_2 i \right> \oplus \left<d^- \otimes_2 i
\right> 
\, \oplus \, 
\left<d^+ \otimes_3 i \right> \oplus \left<d^- \otimes_3 i
\right> \; 
\]
The six subspaces in the definition of $ {(6 \! \cdot \! 49152 )}\!_N$
are invariant under $N_{xyz}$, and the action of the kernel elements
in tables~\ref{table:rep:6:2048:N} and~\ref{table:rep:3:24:N} show
that $K_0$ acts as identity everywhere, so
$ {(6 \! \cdot \! 49152 )}\!_N$ is a representation of 
$N_{xyz0}$. From the actions of $\tau$ and $x_i$ in  these
tables we see that both of them preserve 
$ {(6 \! \cdot \! 49152 )}\!_N$ so that this is indeed 
a representation of~$N_0$. 

Let  $\imath$ be the involution on  $ {(6 \! \cdot \! 49152 )}\!_N$
given by:
\[
d^\pm \otimes_m i \stackrel{\imath}{\longmapsto}
(-1)^{\left<d,i\right>}
\bar{d}^{\,\mp} \otimes_{m\pm 1} i
\, , \quad \mbox{with $m$ to be taken modulo $3$} \; ,
\]
Involution $\imath$ trivially preserves the identification
between $(-d)^\pm \otimes_m i $ and
$- (d^\pm \otimes_m i)$.
Since $\left<\Omega,i\right>=1$, we have
\begin{align*}
\imath((\Omega d)^+ \otimes_m i)  &
=  (-1)^{\left<\Omega d, i\right> }  
\Omega \bar{d}^- \otimes_{m+1} i
= (-1)^{\left<d, i\right> } 
(-\Omega \bar{d})^- \otimes_{m+1} i \\
& =  (-1)^{\left<d, i\right> }
\bar{d}^- \otimes_{m+1} i  
\quad =  \; \imath(d^+ \otimes_m i) \; .
\end{align*}
A similar calculation shows
$\imath(-\Omega d_m^- \otimes i_m) = \imath(d_m^- \otimes i_m)$. 
Thus the identifications $(\Omega d)_m^+ = d_m^+$ and  
$(-\Omega d)_m^- = d_m^-$ are also preserved, so that $\imath$ is 
well defined on $ {(6 \! \cdot \! 49152 )}\!_N$.

For the proof of the Lemma is suffices to show that $\imath$
commutes with the action of the generators of $N_0$.
Using Table~\ref{table:action:196884x},
it is trivial to check that $\imath$ commutes with $\tau$ and
$x_\pi$ for even $\pi$. So it remains to check that
$\imath$  commutes with $x_e$ and $x_\epsilon$ for 
odd $\epsilon$,  which is a
tedious but straightforward calculation using
Table~\ref{table:action:196884x}. The relevant action of $x_e$
is given in Table~\ref{table:action:196884xyz:xe}.

\proofend

\begin{Lemma}
The subspace of $196884_x \oplus 196884_y \oplus 196884_z$ spanned by
the basis vectors $X_{d\cdot i}$, $Y_{d\cdot i}$ and $Z_{d\cdot i}$,
$d \in \mathcal{P}$, $i \in \tilde{\Omega}$, $i$ odd,
$\left<d,i\right> = P(d)$  is a monomial representation of $N_0$.
There is an isomorphism from that space into
$(3 \! \cdot \! 49152 )\!_N$ given by
\[
   X_{d\cdot i} \mapsto \bar{d}^+ \!\otimes_2 i \, , \quad 
   Y_{d\cdot i} \mapsto \bar{d}^+ \!\otimes_3 i \, , \quad 
   Z_{d\cdot i} \mapsto \bar{d}^+ \!\otimes_1 i \, . 
\]
\end{Lemma}

\fatline{Proof}

Let $V$ be the space spanned by   $X_{d\cdot i}$, $Y_{d\cdot i}$ and 
$Z_{d\cdot i}$.The last column in Table~\ref{table:action:196884x}
shows that $V$  is  a representation of $N_0$. So we have
to show that the mapping 
$V \rightarrow (3 \! \cdot \! 49152 )\!_N$ given in the
Lemma preserves the action of the generators, say, 
$x_e, x_\pi$ and $\tau$ of $N_0$. This is trivial for $\tau$ and $x_\pi$. 
The operation of $x_e$ on $V$ and on $(3 \! \cdot \! 49152 )\!_N$
is stated in Table~\ref{table:action:196884xyz:xe}. 

\proofend

\begin{table}[H]
    \centering	
    \begin{tabular}{|l|l|l|}
        \hline
         \multicolumn{3}{|l|}{ 
     \small         Operation    $ B \mapsto  B {x_e}$ of the element
           $x_e$ of  $N_0$ on basis vectors $B$  } \\
    \hline        
    \small $ \quad  B \mapsto  B {x_e}$ 
    & \small  $ \qquad B \mapsto  B {x_e}$ 
    & \small  $ \qquad B \mapsto  B {x_e}$     
    \\
    \hline    
    $X^+_{d \cdot i} \! \mapsto  
    (-1)^{\left<e,i\right>} X^+_{\bar{e}de \cdot i}    $     
    &
    $\vphantom{d^{d^d}} \bar{d}^+ \!  \otimes_2 i \mapsto
     (-1)^{\left<e,i\right>} (\bar{e} \bar{d} e)^+ \!  \otimes_2 i $
    &  
    $ d^- \!  \otimes_3 i  \mapsto 
    (-1)^{\left<e,i\right>} (\bar{e} d e)^- \! \otimes_3 i $
    \\
    \hline
    $Y^+_{d \cdot i}  \mapsto  
    (-1)^{\left<e,i\right>} Y^+_{\bar{e}d \cdot i}    $  
    &
    $\vphantom{d^{d^d}} \bar{d}^+ \! \otimes_3 i   \mapsto  
    (-1)^{\left<e,i\right>} (\bar{d} e)^- \!  \otimes_3 i $
    &  
    $d^- \!  \otimes_1 i\mapsto  (\bar{e} d)^- \! \otimes_1 i $
    \\ 
    \hline
    $Z^+_{d \cdot i} \mapsto   Z^+_{de \cdot i}    $    
    &
    $ \vphantom{d^{d^d}}
   \bar{d}^+ \! \otimes_1 i \mapsto
    (\bar{e} \bar{d})^+ \!  \otimes_1 i $ 
    &
    $d^- \! \otimes_2 i  \mapsto 
    (-1)^{\left<e,i\right>}   (de)^- \!  \otimes_2 i $
    \\
    \hline
    \end{tabular}
    \captionsetup{justification=centering}
    \caption{Operation of  $x_e$ on some basis vectors
      of   $196884_x \oplus 196884_y \oplus 196884_z $.\\
     \small{Entries in the same row are equivalent (up to sign)
         with respect to the dictionary in 
             Table~\ref{table:dictionary}.}    
    }  
    \label{table:action:196884xyz:xe}	   
\end{table}

\subsection{The proof for the non-monomial action of $\tau$}
\label{Subsection:proof:V_T}

In this subsection we will show that $V_T$ is a representation
of $\hat{N}$ with the operation of the generators of $\hat{N}$
given by  (\ref{eqn:operations:V_T}). We already know from
section~\ref{subsection:action:tau} that $V_T$ represents
$N_x$ with kernel $K_T$. In the sequel the generators
$x_d, y_d, z_d, x_\delta, x_\tau, y_\tau$ denote the linear 
operation of the corresponding generator on $V_T$, and we put
$\nu_\pi = x_\pi x_\tau^{|\pi|}$. Our goal is to establish
the  defining relations for $\hat{N} / K_T$ in $V_T$, with
the relations between the generators of $N_x$ already
 established in $V_T$.

By construction $x_\tau$ and $y_\tau$ act as real symmetric orthogonal
matrices on $V_T$, establishing $x_\tau^2 = y_\tau^2 = 1$. Note that
$y_\tau$ consists of $64 \times 64$ blocks, one for each $d$, and
each of these blocks can easily be checked to by symmetric and
orthogonal. 

From $\nu_\pi = x_\pi x_\tau^{|\pi|}$ and the operation of 
$x_\pi$ and $x_\tau$ in (\ref{eqn:operations:V_T}) we conclude that
$\nu_\pi$ maps $X_{d \cdot \delta}^+$ to 
$X_{d^\pi \cdot \delta^\pi}^+$ also for odd $\pi \in \AutStP$,
Thus for any odd diagonal automorphism 
$\varphi \in \mathcal{C^*}$ the element $\nu_\varphi$ acts by
multiplication with $\pm1$ on each 64-dimensional block of
$V_T$ corresponding to a fixed $d$. Thus  $\nu_\varphi$ commutes
with  $x_d, y_d, z_d, x_\delta$ for all $d \in \mathcal{C}$,
$ \delta \in \mathcal{C^*}$ and also with
$x_\tau$ and $y_\tau$. 
From  the operation of $\nu_\pi$ on $X_{d \cdot \delta}^+$ in $V_T$
we conlude $\nu_\pi \nu_{\pi'}$ = $\nu_{\pi \pi'}$ for all 
$\pi, \pi' \in \AutStP$. 
With this information we easily check that
the involution  $\nu_\varphi$ acts on the representation
$V_T$ of $N_x$ by conjugation in the same way the element
$\nu_\varphi$ of $\hat{N}$ acts on the group $N_x$. Thus $V_T$
represents the subgroup $\hat{N}_x$ of $\hat{N}$ generated by
$N_x$ and  $\nu_\varphi$, $\varphi$ odd.

So $V_T$ also represents the normal subgroup $\hat{N}_{xyz}$ of
$\hat{N}$ generated by $N_{xyz}$ and  $\nu_\varphi$, $\varphi$ odd.
We claim that the involution $y_\tau$ acts on the representation $V_T$ 
of $\hat{N}_{xyz}$ by conjugation in the same way as the element 
$y_\tau$ of $\hat{N}$ acts on the normal subgroup $\hat{N}_{xyz}$ of 
$\hat{N}$. By (\ref{eqn:operations:V_T}), $y_\tau$ commutes with all  
$x_\pi$, $\pi$ even. We have already shown that $y_\tau$ commutes 
with an odd $x_\varphi$, hence it commutes with all $\nu_{\pi}$,
establishing the claim for all generators $\nu_\pi$.

To prove our claim, we will show that $y_\tau x_e = z_e y_\tau$ holds 
in $V_T$. To see this, note that
\[
X^+_{d \cdot \delta} \stackrel{x_e y_\tau z_e}{\longmapsto}
{\textstyle\frac{1}{8}} \sum_{\epsilon \in A(d,\mathcal{C})}
(-1)^{\left<e,\delta\right> + |\delta \cap \epsilon|}
X^+_{d \cdot \epsilon A(d,e)}
= {\textstyle\frac{1}{8}} \sum_{\epsilon \in A(d,\mathcal{C})}    
(-1)^{m(\epsilon)}   X^+_{d \cdot \epsilon}  \, ,
\]
with $m(\epsilon) = 
|\delta \cap (\epsilon + A(d,e))| + \left<e,\delta\right>
=  |\delta \cap \epsilon| + |\delta \cap d \cap e| 
+ |\delta \cap e|
=  |\delta \cap \epsilon| \pmod{2}$.
Hence $x_e y_\tau z_e = y_\tau$, implying
 $y_\tau x_e = z_e y_\tau$. Together with the relations already
established in $V_T$ we obtain $y_\tau x_e y_\tau = z_e$,  
$y_\tau z_e y_\tau = x_e$ and $y_\tau y_e y_\tau = y_e$, implying 
our claim.  So $V_T$ also represents the group $\hat{N}_y$ 
generated by  $\hat{N}_{xyz}$ and $y_\tau$.

Thus $V_T$ represents an extension with normal subgroup
$\hat{N}_{xyz}$ and factor group generated by the involutions
$x_\tau$ and $y_\tau$, where $x_\tau$ and $y_\tau$ operate on 
$\hat{N}_{xyz}$ by conjugation in $V_T$ in the same way as in 
$\hat{N}$. $\hat{N}$ has structure $\hat{N}_{xyz} : S_3$ with $S_3$ 
the symmetric permutation group of three elements generated by 
$x_\tau$ and $y_\tau$. So in order to show that $V_T$ represents 
$\hat{N}$  it suffices to show  $\tau^3 = 1$ 
(with $\tau = y_\tau x_\tau$) in $V_T$. 

Since $\tau^3=1$ in  $\hat{N}$ and conjugation with $\tau$ in $V_T$ 
is the same as in $\hat{N}$, the matrix $\tau^3$ in $V_T$ 
centralizes all generators in $V_T$.
By construction the generators $x_e, y_e$ and $\tau$ consist of
blocks of $64 \times 64$ matrices, with one block for each octad $d$. 
Thus for each $d$ the block of $\tau^3$ corresponding to $d$ 
centralizes the blocks of all matrices $x_e$ and $y_e$ corresponding 
to the same $d$. It is easy to see that for any such $d$ the 
corresponding blocks of the matrices $x_e$ and $y_e$ generate the 
complete ring of all $64 \times 64$ matrices. So the
block of $\tau^3$ corresponding to $d$ is a multiple of the identity
matrix, and since $\tau$ is orthogonal, that block of matrix $\tau^3$  
is equal to ${\pm1}$. The trace of that block of matrix $\tau$ is
equal to $\frac{1}{8} \sum_{\epsilon \in A(d,\mathcal{C})} 
(-1)^{|\epsilon|/2}$ = $\frac{1}{8}(36-28)$ =  $1$.
It is easy to check that a real $64 \times 64$-matrix~$\tau$  with 
$\tau^3=-1$ has eigenvalues $-1, (1\pm \sqrt{-3})/2$ and hence trace 
equal to $-1 \pmod{3}$. This  establishes $\tau^3=1$ in $V_T$.

\proofend

%%%%%%%%%%%%%%%%%%%%%%%%%%%%%%%%%%%%%%%%%%%%%%%%%%%%%%%%%%%%%%%%%
%%%%%%%%%%%%%%%%%%%%%%%%%%%%%%%%%%%%%%%%%%%%%%%%%%%%%%%%%%%%%%%%%
\section{Extending representation $196884_x$ from
	$N_{x0}$ to  $2_+^{1+24}.\mbox{Co}_1$  }
\label{section_xi}
%%%%%%%%%%%%%%%%%%%%%%%%%%%%%%%%%%%%%%%%%%%%%%%%%%%%%%%%%%%%%%%%%
%%%%%%%%%%%%%%%%%%%%%%%%%%%%%%%%%%%%%%%%%%%%%%%%%%%%%%%%%%%%%%%%%

In this section we extend the representation of  $196884_x$ from
$N_{x0}$ to a maximal subgroup $G_{x0}$ (with structure 
 $2_+^{1+24}.\mbox{Co}_1$) of the monster $\mathbb{M}$. 
It turns out that such an extension is possible for each of the
building blocks $24_x$,  $4096_x$ and $98280_x$ of $196884_x$. 
For an optimized computer construction of $\mathbb{M}$ we need an
explicit representation of a $\xi \in G_{x0} \setminus N_{x0} $ in
all these building blocks such that $N_{x0}$  and $\xi$ generate
$G_{x0}$.  Note that \cite{Conway:Construct:Monster} contains
an explicit construction of $N_0$, but no  construction any 
specific $\xi \in \mathbb{M} \setminus N_0$. 
Here the hardest part is the extension of $N_0$ to $G_{x0}$
and the construction of a suitable 
element $\xi$ in the representation $4096_x$.

We first describe the representation theory of the extraspecial
2~group $2_+^{1+2n}$ and of certain extensions of such groups
called {\em holomorphs} in \cite{Griess:12:sporadic} and 
\cite{citeulike:Monster:Majorana}. In
section~\ref{subsection:Nx0_fiber:product} we present $N_{x0}$ as
a central quotient of a fiber product corresponding to
the faithful representation  $4096_x \otimes 24_x$ of $N_{x0}$.
In  section \ref{subsection:subgroup:Gx0}
we extend that fiber product to
the group $G_{x0}$. In  sections~\ref{Sect:construct:xi}~ff.
we give an explicit construction of a suitable 
$\xi \in G_{x0} \setminus N_{x0}$  of order~3 and of its 
representation. The action of $\xi$ and $\xi^2$ on the components
of $196884_x$ is stated in Lemma~\ref{Lemma:x_dh:xi}   for
$98304_x$, in Corollary~\ref{Corr:op:xi:24x} for $24_x$, and in
in Corollary~\ref{eqn:transfrom_dpm_xi} for $4096_x$.

%%%%%%%%%%%%%%%%%%%%%%%%%%%%%%%%%%%%%%%%%%%%%%%%%%%%%%%%%%%%% 
\subsection{Extraspecial 2-groups and their representations} 
\label{section:extra:special:2}
%%%%%%%%%%%%%%%%%%%%%%%%%%%%%%%%%%%%%%%%%%%%%%%%%%%%%%%%%%%%%

Recall the definition of quadratic, bilinear any symplectic
forms on vector spaces over $\Field_2$ from 
section~\ref{subsection:cocycle}.
Symplectic bilinear or non-singular
quadratic forms exist only on spaces $\Field_2^{2n}$ of even
dimension. A quadratic form $q$ on $\Field_2^{2n}$  is 
of  {\em plus type} if it is non-singular and there is an 
$n$-dimensional subspace $V$ of $\Field_2^{2n}$ 
with $q(x) = 0$ for all $x \in V$. It is known  that all 
quadratic forms of plus type on  $\Field_2^{2n}$ are 
equivalent under the linear group $\SL_{2n}(2)$. 

A finite 2-group $E$ is said to be {\em extraspecial} if its
center $Z(E)$ has order~2, the factor group $E/Z(E)$ is elementary
abelian, and the commutator group  $[E,E]$ of $E$ is contained 
in $Z(E)$.
We usually write $-1$ for the non-identity element in $Z(E)$, 
so  $Z(E) = \{\pm1 \}$.
The elementary Abelian 2-group $E/ \{\pm1 \}$  may be regarded as
vector space over $\setF_2$. 

Define $P: E/ \{\pm1 \} \rightarrow \setF_2$ by  
$x^2 = (-1_E)^{P(x)}$, $x \in E$.
Then $P$ is well-defined on  $E/ \{\pm1 \}$, and by  (23.10) 
in \cite{Aschbacher-FGT}, $P$ is a non-singular quadratic form on 
the vector space  $E/ \{\pm1 \}$, and we have 
$[x,y] = (-1)^{\beta_P(x,y)}$ for all $x,y \in E$,
with  $\beta_P$  symplectic.
Thus $E$ has order $2^{2n+1}$.

An extraspecial 2-group $E$ is of plus type if the corresponding  
quadratic form $P$ on $E/ \{\pm1 \}$ is of plus type.
For any $n \in \setN$ there is a unique extraspecial 2-group
of order $2^{1+2n}$ of plus type denoted by~$2_+^{1+2n}$.

The remainder of this subsection deals with the representation
theory of $2_+^{1+2n}$ and of certain extensions of that group. 

Let $E$ be an extraspecial 2-group of type $2_+^{1+2n}$.
We write  $\setR E$ for the group ring of $E$ over the real
field $\setR$. So  $\setR E$ is a real algebra with the elements
of $E$ as basis vectors, and the basis vectors are multiplied as
in $E$.
Let $\setR E^\pm$ be the quotient algebra of  $\setR E$, where we
identify the real number $-1$ with element $-1$ of $Z(E)$. 
It is well known that the irreducible complex representations of 
the group $2_+^{1+2n}$ consist of a unique
faithful real $2^n$-dimensional  representation and of
$2^{2n}$ real one-dimensional representations of the
elementary Abelian 2 group $G/\{\pm1\}$, see 
e.g.~\cite{Griess:Monster:Algebra}. 
So it is obvious that
$\setR E^\pm$ is isomorphic to the unique faithful irreducible
real representation of $E$.

A {\em holomorph} of a group $E$ is an extension with
normal subgroup $E$ and factor group $\Out E$, where $\Out E$ is
the group of outer automorphisms of $E$, see
\cite{Griess:Monster:Algebra,Griess:12:sporadic, citeulike:Monster:Majorana}.
We remark that the meaning of the term 'holomorph' in
\cite{Griess:12:sporadic} and  \cite{citeulike:Monster:Majorana} differs 
from the meaning of that term in the older group-theoretic literature.
If $E$ is of type  $2_+^{1+2n}$,
the  group $\Out E$ is the
orthogonal group  $\mbox{O}_{2n}^+(2)$. Here the orthogonal
group  $\mbox{O}_{2n}^+(2)$ preserves the quadratic form $P$ of
plus type on the vector space $E/\{\pm 1\}$  defined above.
There is a unique holomorph $\mathfrak{H}(E)$ of the group
$E = 2_+^{1+2n}$ 
which has a faithful $2^n$-dimensional real representation, 
called the {\em standard holomorph} of  $E$,
see \cite{Griess:Monster:Algebra}, Appendix~1
or \cite{citeulike:Monster:Majorana}, Lemma~1.4.2.
Since  $\setR E^\pm$ is isomorphic to the unique faithful
irreducible representation of $E$ of dimension~$2^n$, it is
the restriction of the $2^n$--dimensional real representation of 
$\mathfrak{H}(E)$  to $E$. Schur's lemma implies that an element 
of $\setR E^\pm$ which acts as an automorphism on~$E$ is determined
by that action up to a scalar multiple, which must be $\pm1$.

The following lemma helps us to construct elements of the
standard holomorph $\mathfrak{H}(E)$:

\begin{Lemma}
\label{Lemma:commutator:holomorph}    
 \label{lemma:construct:holomorph}	
Let $E = 2_+^{1+2n}$  and $H \subset E$ an elementary
Abelian 2-group not containing $-1$. We consider $H$ and
$ E/\{\pm1\}$ as vector spaces over $\Field_2$.
Let 
\[ 
  \xi = |H|^{-1/2} \sum_{z \in H} (-1)^{q(z)} z
\]
be an element of $\setR E^{\pm}$, with $q$ a
non-singular quadratic form on $H$ with associated bilinear form
$\beta_q$.
	
Then $\xi^2=1$ and there is a unique
linear mapping $\phi: E/\{\pm1\}  \rightarrow H$ with
$[x,y] = (-1)^{\beta_q(\phi(x),y)}$ for all $x \in E$,
$y \in H$. For all $x \in E$ we have 
\[
  [x,\xi] = (-1)^{q(\phi(x))} \phi(x) \in Z(E) \cdot H  \; .
\]
\end{Lemma}

%Note that Lemma~\ref{Lemma:commutator:holomorph} implies 
%$\xi^x \in E$ for all $x \in E$ and hence $\xi \in  \mathfrak{H}(E)$.

\fatline{Proof}
\begin{align*}
|H| x^{-1}\xi x \xi
& =   \sum_{y,z \in H}
x^{-1}  (-1)^{q(y)} y x  (-1)^{q(z)} z
=    \sum_{y,z \in H} 
(-1)^{q(y)+q(z)} [x,y] yz \\
&  =    \sum_{y,z \in H} 
(-1)^{q(y)+q(yz)} [x,y] z
= \sum_{z \in H} (-1)^{q(z)} z 
\sum_{y \in H} (-1)^{\beta_q(z,y) + C(x,y)}   \; ,     
\end{align*}
where $C$ is the bilinear form on $E/\{\pm1\}$ with
$[x.y] = (-1)^{C(x,y)}$.
Since both, ${\beta_q(z,y)}$ and $C(x,y)$ are linear in
$y$, the last sum over $y$ is equal to $|H|$ if
$\beta_q(z,y)= C(x,y)$ for all $y \in H$ and zero otherwise. 
Since $\beta_q$ is non-singular, it follows from linear algebra 
that for each $x \in E/\{\pm1\}$
there is a unique $z \in H$ with $\beta_q(z,y)= C(x,y)$ for 
all $y \in H$,  
and  that the mapping $\phi$ which maps
each $x \in E/\{\pm1\}$ to that value $z$ is linear.
This proves $x^{-1}\xi x \xi = (-1)^{q(z)} z$, with
$z = \phi(x)$. 
In case $x=1$ we have $\phi(x)=0$, obtaining $\xi^2=1$,
and hence  $[x, \xi] = (-1)^{q(z)} z$.

\proofend
\stepcounter{Theorem}

Conjugation with an element $\xi$ of $\setR E^{\pm}$ 
constructed in Lemma~\ref{lemma:construct:holomorph} is an
automorphism of~$E$, since $[x,\xi] \in E$ for all $x \in E$.
Hence $\xi$ represents an element of $\mathfrak{H}(E)$.
$\setR E^{\pm}$ is a faithful irreducible representation of $E$ and
also of $\mathfrak{H}(E)$. So
we may identify the element $\xi$ of $\setR E^{\pm}$
with the element of $\mathfrak{H}(E)$ represented by $\xi$.

\fatline{Remark}

In Lemma~\ref{Lemma:commutator:holomorph},
$\xi$ operates on the vector space $V = E/\{\pm1\}$ by conjugation
as an orthogonal transformation  with $\im (\xi -1) = H$.
The coset $\xi E$ of $E$ can be considered as an element
of  $\mathfrak{H}(E) / E \cong O_{2n}^+$, and the bilinear
form $\beta_q$ is uniquely determined by the  coset $\xi E$.
For each element $X$ of an orthogonal group operating on a vector
space $V$, Wall \cite{Wall:1963} 
has defined a nondegenerate bilinear form $F_X$ on the
image of $X - 1$ in $V$ called the {\em parametrization} of $X$.
It can be shown that $\beta_q$ is just Wall's parametrization 
of $\xi E$.

%%%%%%%%%%%%%%%%%%%%%%%%%%%%%%%%%%%%%%%%%%%%%%%%%%%%%%%%%%%%%
\subsection{$N_{x0}$ is a central quotient of a fiber product}
\label{subsection:Nx0_fiber:product}
%%%%%%%%%%%%%%%%%%%%%%%%%%%%%%%%%%%%%%%%%%%%%%%%%%%%%%%%%%%%%

In section~\ref{section:rep:Nx1} we have constructed
representations $4096_x$ and $24_x$ of the group $N_x$
which both are also representations of $N_{x1} = N_x / K_1$.
By Tables~\ref{table:rep:6:2048:N} and \ref{table:rep:3:24:N}
the kernels $K(4096_x)$ and $K(24_x)$ of these two 
representations of $N_x$ intersect in $K_1$.
Let $N(4096_x) = N_x /K(4096_x)$ and $N(24_x) = N_x /K(24_x)$
be the quotients of $N_x$ for which the representations 
$4096_x$ and $24_x$ are faithful. 
Put $N_x^*$ =  $N_x//(K(4096_x) K(24_x))$. Then
\[
  N_x^*  \; \cong \;  N(4096_x) / (K(24_x)/K_1)
        \;  \cong  \;  N(24_x) / (K(4096_x)/K_1) \, ,
\]
where e.g. the natural injection from $K(24_x)/K_1$ into
$N(4096_x)$ is obtained by extending the coset $x K_1$ to
$x K_{4096}$ for $x \in K(24_x)$.

Then $N_{x1}$ is isomorphic to the fiber product
$N(4096_x) \bigtriangleup_{N^*_x} N(24_x)$.
If $G_1, G_2$ are groups with a common factor group $H$
and homomorphisms $\phi_i: G_i \rightarrow H$, $i=1,2$, 
then the fiber product $G_1 \bigtriangleup_{H} G_2$ is
the subgroup of the direct product $G_1 \times G_2$
defined by:
\[
 G_1  \bigtriangleup_{H} G_2 \; = \;  
    \left\{(x,y)   \in G_1 \times G_2  
     \mid  \phi_1(x) = \phi_2(y)  \right\} \; .
\]  
If $G_i$ has a center $\{\pm1_{i}\}$ of order 2 for $i=1,2$,
then we write $\frac{1}{2}  (G_1  \bigtriangleup_{H} G_2)$ for 
the quotient  of $G_1  \bigtriangleup_{H} G_2$ by
$\{(1_1, 1_2) , (-1_1, -1_2)\}$ as in \cite{Atlas}.
Since $\mathcal{K}_2$ and $\mathcal{K}_3$ act as $-1$ in both,
$4096_x$ and $24_x$, the group $N_{x0}$ is isomorphic to
$\frac{1}{2}(N(4096_x) \bigtriangleup_{N^*_x} N(24_x))$.

The following diagram, which is essentially a copy of
Fig. 2 in \cite{Conway:Construct:Monster}, depicts the 
homomorphisms from $N_x$ to the various factor groups 
of~$N_x$ defined above.
For each homomorphism we show a generating system 
and the structure of the corresponding factor group. 

\vspace{1.0em}

\begin{tikzpicture}
\matrix (m) [matrix of math nodes,row sep=2.5em,column sep=3em,minimum width=2em]
{
     &         &[+1em]  N(4096_x)  &[+1em]        &      \\
     N_x & N_{x1}  &[+1em]             &[+1em]  N_x^*   & 
         N(24_x) / E  & 
         1_{\vphantom{x}}     \\
     \hphantom{x}  &[+1em]             & 
     \hphantom{0}  N(24_x)  \hphantom{0}   &[+1em]  &\\	
};         	
\path[-stealth] 
   (m-2-1) edge node [above]
        {$\scriptstyle \mathcal{K}_1$}
               node [below] {$\scriptstyle 2$}  
   (m-2-2)
   (m-2-2) edge node [above]
   {$\scriptstyle y_{\Omega^{\vphantom{2}}},
      z_{-\Omega^{\vphantom{2}}} \hphantom{mi}$} 
      node [below] {$\scriptstyle 2^{\vphantom{2}}$}   
   (m-1-3)
   (m-2-2) edge node [above]
      {$\scriptstyle \hphantom{m} 
        x_{d^{\vphantom{2}}},
         x_{\delta^{\vphantom{2}}}$} 
        node [below] {$\scriptstyle 2_+^{1+24^{\vphantom{2}}}$} 
   (m-3-3)       
   (m-1-3) edge node [above]
       {$\scriptstyle \hphantom{m} 
       x_{d^{\vphantom{2}}},
       x_{\delta^{\vphantom{2}}}$} 
       node [below] {$\scriptstyle 2_+^{1+24^{\vphantom{2}}}$}
   (m-2-4)       
   (m-3-3) edge node [above]
      {$\scriptstyle y_{\Omega^{\vphantom{2}}}, 
      z_{\Omega^{\vphantom{2}}} \hphantom{mi}$} 
      node [below] {$\scriptstyle 2^{\vphantom{2}}$}   
   (m-2-4)
   (m-2-4) edge node [above]
      {$\scriptstyle y_{d^{\vphantom{2}}}, z_{d^{\vphantom{2}}} $} 
      node [below] {$\scriptstyle 2^{11}$}   
   (m-2-5)
   (m-2-5) edge node [above]
     {$\scriptstyle x_{\pi^{\vphantom{2}}} $} 
      node [below] {$\scriptstyle M_{24}$}   
   (m-2-6)
      ;
        ) ;
\end{tikzpicture} 
 
From (\ref{equ:structure:Nx}) we see that   
$N(4096_x)$ has structure $ 2_+^{1+24}.2^{11}.\mbox{M}_{24}$,
$N(24_x)$ has structure $2^{12}:\mbox{M}_{24}$
and $N_x^*$ has structure $2^{11}:\mbox{M}_{24}$.
The relations (\ref{relations:Qx:extraspecial}) show that the 
subgroup $Q_x$ of $N_x$ generated by $x_d, x_\delta$ has
structure  $ 2_+^{1+24}$. Note that $|Q_x \cap K_1| = 1$, 
so $Q_x$ is also isomorphic to a subgroup of $N_{x1}$.

By construction, $N(24_x)$ acts as a matrix group on
$\setR^{24}$ 
with coordinates labelled by $\tilde{\Omega}$.  $N(24_x)$  has
an elementary Abelian normal subgroup $E$ of order  $2^{12}$.
An element of $E$ corresponds to the negation of the coordinates
given by a codeword of $\mathcal{C}$, hence $E \cong \mathcal{C}$. 
The permutation representation of  $\mbox{M}_{24}$ acts as
a complement  of $E$ in $N(24_x)$.

%%%%%%%%%%%%%%%%%%%%%%%%%%%%%%%%%%%%%%%%%%%%%%%%%%%%%%%%%%%%% 
\subsection{The maximal subgroup
	$ G_{x0} = 2_+^{1+24}.\mbox{Co}_1$ of $\mathbb{M}$}
\label{subsection:subgroup:Gx0}
%%%%%%%%%%%%%%%%%%%%%%%%%%%%%%%%%%%%%%%%%%%%%%%%%%%%%%%%%%%%%

As in \cite{Conway:Construct:Monster}, we will enlarge 
$N(4096_x)$ and $N(24_x)$ to larger groups of 
$G(4096_x)$ and $G(24_x)$ of structure
$2^{1+24}.\mbox{Co}_1$ and $2.\mbox{Co}_1$,
respectively. From this we obtain the group
$G_{x0} = \frac{1}{2}
     ( G(4096_x) \bigtriangleup_{\small \mbox{Co}_1} G(24_x) )$,
which is a maximal subgroup of the monster.

The group $N(4096_x) = N_{x1}/K(4096_x)$ defined in 
section~\ref{subsection:Nx0_fiber:product} has the faithful 
irreducible real representation $4096_x$.  $4096_x$ is also
faithful and irreducible for the extraspecial subgroup $Q_{x}$ of
type $2^{1+24}_+$ of $N(4096_x)$  and hence also
for the standard holomorph $\mathfrak{H}(Q_x)$ of $Q_{x}$ . 

By Theorem~\ref{Theorem:Q:Leech} the quotient of $Q_{x}$ by its 
center is isomorphic to $\Lambda / 2 \Lambda$, and the squaring
map in $Q_x$ is equal to the mapping
$\Lambda / 2 \Lambda \rightarrow \Field_2$ given by 
$\lambda \mapsto \mbox{type}(\lambda)
   = \left< \lambda,\lambda \right> /2$. 
This 'type'  mapping is a quadratic form on $\Lambda / 2 \Lambda$
by construction, which is invariant under the automorphism group
$\mbox{Co}_1$ of $\Lambda / 2 \Lambda$. It is of plus type by
(\ref{relations:Qx:extraspecial}), so  $\mbox{Co}_1$ 
is a subgroup of $O_2^+(24)$. The factor group $N_x^*$ of
$ N(4096_x)$ with structure $2^{11}.M_{24}$ is isomorphic to
the monomial subgroup of the automorphism group
$\mbox{Co}_1$ of $\Lambda / 2 \Lambda$.
This leads to a chain of inclusions
\[
\begin{array}{lccccccc}
     &  Q_x
   & \subset &   N(4096_x)
   & \subset &   G(4096_x)
   & \subset & \mathfrak{H}(Q_x) 
  \\
\mbox{with structures} \hphantom{m} &  2_+^{1+24}
   & \subset &  2_+^{1+24}.2^{11}.M_{24}
   & \subset &   2_+^{1+24}.\mbox{Co}_1 
   & \subset &   2_+^{1+24}.\mbox{O}^+_{24}(2) \, . 
\end{array}     
\] 

There is a similar chain of inclusions of groups
\[  
\begin{array}{lccccccc} 
   &  E
   & \subset &   N(24_x)
   & \subset &   G(24_x)
   & \subset &  \mbox{S0}_{24} (\setR) \\
\mbox{with structures $\qquad$} & \!\!\! 2^{12}
   & \subset &  2^{12}.M_{24}
   & \subset &   2.\mbox{Co}_1 \quad . 
   % & \subset &   ?? .
\end{array}
\]
Here the group $G(24_x)$ is the automorphism group of the 
Leech lattice $\Lambda$, and $ N(24_x)$ is the monomial
subgroup of that group. So the groups $G(4096_x)$
and $G(24_x)$ possess natural homomorphisms onto
$\mbox{Co}_1$ which are extensions of  the homomorphisms
from  $N(4096_x)$ and $N(24_x)$ onto
$N^*_x = 2^{12}:M_{24}$. We can therefore extend the
fiber product $N_{x1}$ of   $N(4096_x)$ and $N(24_x)$
to  the fiber product
\[
    G_{x1} = G(4096_x) 
     \, \bigtriangleup_{\mbox{\small Co}_1}  \,  G(24_x) \; .
\]  
This group has representations of degrees 4096 and 24 extending
the representations $4096_x$ and $24_x$, which we will also
call  $4096_x$ and $24_x$. The tensor product
$4096_x \otimes 24_x$ identifies the centers $\{\pm 1\}$ of both 
of its factors and is hence is a representation of
\[
    G_{x0} =  {\textstyle \frac{1}{2} } \left( 
  G(4096_x)  \, \bigtriangleup_{\mbox{\small Co}_1} \, G(24_x)    
    \right) \; .
\] 
The group $G_{x0}$ is the maximal subgroup of the monster 
$\mathbb{M}$ constructed in \cite{Conway:Construct:Monster}.
The groups $G_{x0}$ and $ G(4096_x)$ are both of structure
$2^{1+24}.\mbox{Co}_1$, but  not isomorphic, see
e.g. \cite{Griess:Monster:Algebra, citeulike:Monster:Majorana}.

We have found representations $4096_x$ and $24_x$ of $G_{x1}$.  
This leads to representations $98304_x$ = $4096_x \otimes 24_x$ 
and $300_x$ = $24_x \otimes_{\mbox{\tiny sym}} 24_x$ of
$G_{x0}$. $G_{x0}$ (as a extension  of its normal
subgroup $Q_x$) permutes the short elements $x_r$ of
$Q_{x}$ via conjugation. Thus the representation $98280_x$ of
$N_{x0}$ can also be extended to a monomial representation of
$G_{x0}$, which we will also call $98280_x$. So we may build a 
representation $196884_x$ of $G_{x0}$ from its components in
the same way as in the case of $N_{x}$.
 
%%%%%%%%%%%%%%%%%%%%%%%%%%%%%%%%%%%%%%%%%%%%%%%%%%%%%%%%%%%%% 
\subsection{Construction of a 
	$\xi \in G_{x0} \setminus N_{x0}$ and  operation
 of $\xi$  on $98280_x$} 
\label{Sect:construct:xi}
%%%%%%%%%%%%%%%%%%%%%%%%%%%%%%%%%%%%%%%%%%%%%%%%%%%%%%%%%%%%%

It remains to construct a specific element 
$\xi \in G_{x0} \setminus N_{x0}$. By \cite{Atlas}, 
$2^{11}:M_{24}$ is maximal in $\mbox{Co}_1$, so $ N_{x0}$ is
maximal in $G_{x0}$, and hence any 
$\xi \in G_{x0} \setminus N_{x0}$ together with $N_{x0}$ 
generates  $G_{x0}$.

$4096_x$ is also a faithful irreducible representation of
$Q_x = 2^{1+24}$, so  $\hom(4096_x, 4096_x)$ is isomorphic to 
$\setR Q_x ^\pm$ as an algebra. Thus we may define 
$\xi$ as an element $\xi = \xi_{4096} \otimes \xi_{24}$ of
 $\setR Q_x^\pm \otimes \hom(24_x, 24_x)$. Conjugation
of $Q_x$ with $\xi_{4096}$ may be computed by 
Lemma~\ref{Lemma:commutator:holomorph}. The action of
$\xi_{4096}$ determines the action of $\xi$ in $98280_x$
uniquely and  in $24_x$ up to sign. We will construct a
specific element $\xi_{4096}$ or order~3. Then we let
$\xi_{24}$ be the unique element of order~3 corresponding
to  $\xi_{4096}$, i.e. $-\xi_{24}$ has order~6.   
In the sequel we will abbreviate $\xi_{4096}$ to $\xi$. 

The decompositions $\mathcal{C} = \mathcal{G} \oplus \mathcal{H}$
and $\mathcal{C}^* = \mathcal{G}^* \oplus \mathcal{H}^*$
of $\mathcal{C}$ and $\mathcal{C}^*$ into
grey and coloured subspaces discussed in 
section~\ref{subsection:gray:col} is useful for describing
the action of $\xi$.

Recall from 
section~\ref{subsection:gray:col}
that $\mathcal{G}$ has a natural
basis $g_0,\ldots,g_5$ and that $\mathcal{G}^*$ has a natural
basis $\gamma_0,\ldots,\gamma_5$. 
Let $w: \mathcal{G} \cup \mathcal{G^*} \rightarrow \setZ$ be the
weight of a vector in the corresponding natural basis, as
in Definition~\ref{Def:w}. Let
$w_2: \mathcal{G} \cup \mathcal{G^*} \rightarrow \Field_2$
as in (\ref{Def:w_2}). So 
$w_2(d) = \binom{w(d)}{2} \pmod{2}$.
Let  $\gamma: \mathcal{C} \rightarrow \mathcal{G}^*  $ 
be as in (\ref{Def:gamma}).
By Lemma~\ref{Lemma:scalprod:gamma}, $w_2$ is a non-singular
quadratic form on $ \mathcal{G}$ 
with associated bilinear form
\begin{align}
\label{lemma:w_2}  
   \beta_{w_2} = \biscalar{.,.} \,, \qquad
\mbox{where}  \quad   
 \biscalar{d,e} = \left<d, \gamma(e)\right> =
  \left<e, \gamma(d)\right>  \; .
\end{align}
By (\ref{eqn_gamma_i}) the restriction of  $\gamma$ to $\mathcal{G}$
is an isomorphism  $\mathcal{G} \rightarrow \mathcal{G}^*$ with
$\gamma(g_i)=\gamma_i$. 
Thus the  mapping $w_2:  \mathcal{G}^* \rightarrow \Field_2$ 
is also a non-singurlar quadratic form on $  \mathcal{G}^*$.
We also define the bilinear form  $\biscalar{.,.}$ on
$\mathcal{G^*}$ by decreeing
$\biscalar{\delta,\epsilon}$  = $\beta_{w_2}(\delta,\epsilon)$.

For $d \in \mathcal{P}$ (or $d \in \mathcal{C}$) let
$\tilde{x}_d   = x_{-1}^{\scriptsize \mbox{sign}(d)} x_d  x_{\theta(d)}$
 be as in (\ref{eqn:def:tilde:xd}). Then
$\{x_\epsilon \mid \epsilon \in \mathcal{G}^* \}$ and 
$\{\tilde{x}_e \mid e \in \mathcal{G} \}$ are
elementary Abelian subgroups of $Q_x$ not containing
the central element $x_{-1}$ of $Q_x$.
In order to construct the {\em three-bases} element $\xi$  
we  put $\xi = \xi_\gamma \xi_g  \in \setR G^\pm$ with
\stepcounter{Theorem}
\begin{align}
\label{eqn:def:xi}
  \xi_\gamma = \frac{1}{8} \sum_{\epsilon \in \mathcal{G}^*}
    (-1)^{w_2(\epsilon)} x_\epsilon
   \; , \quad
  \xi_g = \frac{1}{8}  \sum_{e \in \mathcal{G}}
    (-1)^{w_2(e)} \tilde{x}_e \,  .
\end{align}
By Lemma~\ref{Lemma:scalprod:gamma} and (\ref{lemma:w_2}),
 $w_2$ is a non-singular quadratic
form on $\mathcal{G}$ and also on  $\mathcal{G}^*$. Thus
by Lemma~\ref{Lemma:commutator:holomorph},
$\xi_\gamma$ and  $\xi_g$ are involutions in the holomorph
$\mathfrak{H}(Q_x)$, so they both normalize $Q_x$.

$\xi_\gamma$ and  $\xi_g$ are in $\mathfrak{H}(Q_x) $ but
not in $G(4096_x)$.  In section~\ref{subsection_op_xi_24}
we will show that the  product $\xi = \xi_\gamma \xi_g$ 
is in $G(4096_x)$ as required.
The relevant relations for  $\xi_\gamma$ and
 $\xi_g$ are given by:

\begin{Lemma}
\label{Lemma:conjugation:xi}    
\begin{align*}
   \xi_\gamma^2 = \xi_g^2 = (\xi_\gamma \xi_g)^3
   = [\tilde{x}_e, \xi_g] = [x_\epsilon, \xi_\gamma] = 1 \, ,
    \quad &
      e \in \mathcal{C}, \; \epsilon \in \mathcal{C}^* \, ; 
       \\
      [\tilde{x}_h, \xi_\gamma] = [x_\eta, \xi_g] = 1 \, ,
     \quad &
     h \in \mathcal{H}, \; \eta \in \mathcal{H}^* \, ;    
       \\
     \tilde{x}_d^{\xi_\gamma}
       =  x_\delta^{\xi_g}
       = (-1)^{w_2(d)} \tilde{x}_d x_\delta, 
       \quad  [\tilde{x}_d, x_\delta] = 1 \, ,
     \quad &
       d= \mathcal{G} ,  \;
     \delta= \gamma(d)  \in \mathcal{G}^* \,.
\end{align*}
\end{Lemma}   
 
\fatline{Proof}

We have already shown $\xi_\gamma^2 = \xi_g^2 = 1$.
Put 
$\tilde{x}_{r} = \tilde{x}_{e \cdot \epsilon}$ for 
$r = (e,\epsilon) \in \mathcal{C} \times \mathcal{C}^*$.
We define a scalar product $\left<.,.\right>$ on
$(\mathcal{C} \times  \mathcal{C}^*)
\times (\mathcal{C} \times  \mathcal{C}^*)$
by decreeing  
$\left<(e,\epsilon),(f,\varphi)\right>$ =
$ \left<e,\varphi\right>+\left<f,\epsilon\right>$.
Then  (\ref{relations:Qx:extraspecial}) implies
\[
 [\tilde{x}_{r}, \tilde{x}_s] = 
(-1)^{\left<r, s \right>}  \, , 
      \qquad   \mbox{for} \quad
       r,s \in  \mathcal{C} \times \mathcal{C}^* \; .
\]
Thus $\tilde{x}_r$, $r \in   \mathcal{H} \times \mathcal{C}^*$,
commutes with every term $\tilde{x}_s$, 
$s \in  \mathcal{G}^* \;$ in the sum that defines
$\xi_\gamma$. Hence $\tilde{x}_r$ commutes with $\xi_\gamma$.
A similar argument shows that
$\tilde{x}_r$, $r \in   \mathcal{C} \times \mathcal{H}^*$, 
commutes with $\xi_g$.

Next we show  the
formula for $\tilde{x}_d^{\xi_\gamma}$.
Let $\phi_\gamma: Q_x/\{\pm 1\} \rightarrow \mathcal{G}^*$
be  the linear mapping given by
 $\phi_\gamma(\tilde{x}_{g_n}) = \gamma_n$, 
$\phi_\gamma(\tilde{x}_r) = 0 $, 
$r \in  \mathcal{H} \times   \mathcal{C}^*$.
Then 
\[
 [\tilde{x}_{r}, x_\delta] = 
 (-1)^{\left<r, \delta\right>}  =
  (-1)^{\biscalar{\phi_\gamma(\tilde{x}_r),\delta)}} 
 \qquad  \, \mbox{for} \quad 
 r \in \mathcal{C} \times \mathcal{C}^* \; .
\]
Since $\biscalar{.,.}$ is the bilinear form associated
with $w_2$, the last equation and
Lemma~\ref{Lemma:commutator:holomorph} imply:
\[
   [\tilde{x}_d, \xi_\gamma] = (-1)^{w_2(\phi_\gamma(\tilde{x}_d))}
       x_{\phi_\gamma(\tilde{x}_d)} 
       = (-1)^{w_2(d)} x_\delta    \; .      
\]

The proof of the formula for ${x}_\delta^{\xi_g}$ is similar.
Let $\phi_g: Q_x/\{\pm 1\} \rightarrow \mathcal{G}$
be  the linear mapping given by
$\phi_g(\tilde{x}_{\gamma_n}) = g_n$, 
$\phi_g(\tilde{x}_r) = 0 $, 
$r \in  \mathcal{C} \times   \mathcal{H}^*$. Then
\[
[\tilde{x}_{r}, x_d] = 
(-1)^{\left<r, d\right>}  
   = (-1)^{\biscalar{\phi_g(\tilde{x}_r),d)}} 
\qquad  \, \mbox{for} \quad 
r \in \mathcal{C} \times \mathcal{C}^* \; .
\]
This equation together with
Lemma~\ref{Lemma:commutator:holomorph} implies:
\[
[x_\delta, \xi_\gamma] = (-1)^{w_2(\phi_g(x_\delta))}
x_{\phi_g(\tilde{x}_r)} 
= (-1)^{w_2(\delta)} \tilde{x}_d  
= (-1)^{w_2(d)} \tilde{x}_d    \; .      
\]

Since $\biscalar{.,.}$ is associated with the quadratic form
$w_2$ and hence alternating we have:
\[
  [\tilde{x}_d, x_\delta]= (-1)^{\left<d,\delta\right>} 
     =  (-1)^{\biscalar{d,d}} = 1   \; .    
\]

It remains to show $ (\xi_\gamma \xi_g)^3 = 1$. This follows
from $\xi_\gamma^2 = \xi_g^2 = 1$ and
\[
 8 \xi_g^{\xi_\gamma}
    =  \sum_{d \in \mathcal{G}}
        (-1)^{w_2(d)} \tilde{x}_d^{\xi_\gamma}
    = \sum_{d \in \mathcal{G}}  \tilde{x}_d x_\delta
    =  \sum_{\delta \in \mathcal{G}^*}  x_\delta \tilde{x}_d
     =  \sum_{\delta \in \mathcal{G}^*} 
              (-1)^{w_2(\delta)} x_\delta^{\xi_g}
    =   8 \xi_\gamma^{\xi_g} \; .       
\]

\proofend

By Lemma \ref{Lemma:conjugation:xi}, 
$\xi = \xi_\gamma \xi_g$ operates on
 $\tilde{x}_d$, $\in \mathcal{G}$ by conjugation as follows:  
\begin{align}
\label{transf:xi:gray}   
\tilde{x}_d   \stackrel{\xi}{\longmapsto}
{x}_{\gamma(d)}  \stackrel{\xi}{\longmapsto}
(-1)^{w_2(d)}   \tilde{x}_d  {x}_{\gamma(d)} 
\stackrel{\xi}{\longmapsto}  \tilde{x}_d 
\; .
\end{align}

\begin{figure}[!h]
    \centering
    \begin{tikzpicture}
    \matrix (m) [matrix of math nodes,row sep=1em,
    column sep=4em,minimum width=1.5em]
    {
        |[name=k1]|    \tilde{x}_d \vphantom{1^d}
        & |[name=k2]| 
        (-1)^{w_2(d)}  \,  \tilde{x}_d \, x_{\gamma(d)}  
        %  $ (-1)^{\binom{w(d)}{2}}  \;  (-1)^{w(d)(w(d)-1)/2} $     
        & |[name=k3]| x_{\gamma(d)}  \vphantom{1^d}  \\
    };
    \draw[<->] (k1) edge node[auto] {$\xi_\gamma$}  (k2) ;
    \draw[<->] (k2) edge node[auto] {$\xi_g$}  (k3);
    \draw[<->,rounded corners] (k3.north) |-
    node[right,pos=(2.0)] {$\xi_\gamma$}
    +(0.3,0.4) -| +(0.9,0.1) |- (k3.east);
    \draw[<->,rounded corners] (k1.north) |-
    node[left,pos=(2.0)] {$\xi_g$}
    +(-0.3,0.4) -| +(-0.9,0.1) |- (k1.west);
    \end{tikzpicture} 
    \caption{Action of  $\xi_\gamma$ and  $\xi_g$ on
        $\tilde{x}_d \in \mathcal{G}$ and on
        $x_\delta = \gamma(x_d)  \in \mathcal{G}^*$. }  
    \label{figure:xi}    
\end{figure}

\begin{Lemma}
\label{Lemma:x_dh:xi}        
Let $d,e \in \mathcal{G}$, $h \in \mathcal{P}_\mathcal{H}$,
$\eta \in  \mathcal{H}^*$. Put $\delta = \gamma(d)$,
 $\epsilon = \gamma(e)$. Then
$\xi$ operates by conjugation on $Q_x$ as follows:
\begin{align}
   \label{Lemma:x_dh:xi_1} 
 \tilde{x}_d \tilde{x}_h x_\epsilon x_\eta
    &  \stackrel{\xi}{\longmapsto}
 (-1)^{w_2(e) + \biscalar{d,e}}  
     \tilde{x}_e \tilde{x}_h x_{\delta\epsilon} x_\eta \\ 
   \label{Lemma:x_dh:xi_2} 
 \tilde{x}_d \tilde{x}_h x_\epsilon x_\eta
    &  \stackrel{\xi^2}{\longmapsto}
  (-1)^{w_2(d) + \biscalar{d,e}}  
     \tilde{x}_{de} \tilde{x}_h x_{\delta} x_\eta  \\
   \label{Lemma:x_dh:xi_3} 
   x_h  x_{\gamma(h)}  
       &  \stackrel{\xi}{\longmapsto}
   x_h  x_{\gamma(h)}       
\end{align}
\end{Lemma}

\fatline{Proof}

By  Lemma~\ref{Lemma:conjugation:xi}, conjugation with
$\xi$ fixes $\tilde{x}_h$ and $x_{\eta}$. So     
(\ref{Lemma:x_dh:xi_1}) and (\ref{Lemma:x_dh:xi_2}) follow
from (\ref{transf:xi:gray}) and from the commutator rules in 
(\ref{relations:Qx:extraspecial}).

$\xi$ does not depend on the cocycle $\theta$, provided that
$\theta$ satisfies Lemma~\ref{Lemma:cocycle:property},
So by Corollary~\ref{Corol:cocycle:colored} we may assume
$\theta(h) = \gamma(h)$ and hence
 $\tilde{x}_h =  x_h x_{\gamma(h)} $. Then
(\ref{Lemma:x_dh:xi_3}) follows from (\ref{Lemma:x_dh:xi_1}).
 
\proofend

Every element $x_r$ of $Q_x$ is of shape
$\pm \tilde{x}_d \tilde{x}_h x_\epsilon x_\eta$ as given
by Lemma~\ref{Lemma:x_dh:xi}. So that Lemma    
allows us to conjugate any such  $x_r$ with $\xi$
or $\xi^2$. $\xi$ operates on the basis vector $X_r$
of the representation $98280_x$ in the same way as it
operates on  $x_r$  by conjugation. So 
Lemma~\ref{Lemma:x_dh:xi} also gives us the operation
of $\xi$ and $\xi^2$ on  $98280_x$.

%%%%%%%%%%%%%%%%%%%%%%%%%%%%%%%%%%%%%%%%%%%%%%%%%%%%%%
\subsection{The operation of $\xi$ on $24_x$}
\label{subsection_op_xi_24}
%%%%%%%%%%%%%%%%%%%%%%%%%%%%%%%%%%%%%%%%%%%%%%%%%%%%%%%%

We also define a linear transformation 
$\xi_{24}=\xi_{24a} \xi_{24b}$ on $24_x$, where
$\xi_{24a}$ and $\xi_{24b}$ are $24 \times 24$ 
matrices operating on $24_x$ by right multiplication. 
Matrices  $\xi_{24a}$ and
$\xi_{24b}$ consist of six identical $4 \times 4$--blocks  
$\xi_{4a}$ and $\xi_{4b}$, respectively, where each of the
$4 \times 4$--blocks transforms the four basis vectors of
$24_x$  labelled by a column of the MOG, as given
in (\ref{elements:MOG}). We put:
\begin{align}
\label{def:xi4ab}
\xi_{4a} = \frac{1}{2} 
\left( \begin{array}{cccc}
\scriptstyle 1 & \scriptstyle -1 & \scriptstyle -1 & \scriptstyle -1 \\
\scriptstyle -1 &\scriptstyle  1 &\scriptstyle  -1 &\scriptstyle  -1 \\
\scriptstyle -1 &\scriptstyle  -1 &\scriptstyle  1 &\scriptstyle  -1 \\
\scriptstyle -1 &\scriptstyle  -1 &\scriptstyle  -1  &\scriptstyle  1 \\
\end{array} \right) 
\; , \quad
\xi_{4b} = \left( \begin{array}{cccc}
\scriptstyle -1 & \scriptstyle 0 & \scriptstyle 0 & \scriptstyle 0 \\
\scriptstyle 0 & \scriptstyle 1& \scriptstyle 0 & \scriptstyle 0  \\
\scriptstyle 0 &\scriptstyle  0 & \scriptstyle 1 & \scriptstyle 0 \\
\scriptstyle 0 & \scriptstyle 0 &\scriptstyle  0 & \scriptstyle 1 \\
\end{array} \right) 
\quad .
\end{align}
We define $\Lambda^E =
 \{\lambda \in \Lambda \mid 
  \left<\lambda_{\tilde{\Omega}},\lambda\right> = 
\left<\lambda_\omega,\lambda\right> = 0 \pmod{2} \} $, with
$\tilde{\Omega} \in \mathcal{G}, \omega \in \mathcal{G}^*$ 
as in section~\ref{section:Golay}. Then  $\Lambda^E$ is a
sublattice of $\Lambda$ of index~4. 

$\xi_{4a}$ and $\xi_{4b}$ are orthogonal by construction.
Thus $\xi_{24a}, \xi_{24b}$ and $\xi_{24}$ are also
orthogonal. Direct calculation shows that 
$\xi_{4a}  \xi_{4b}$ and hence also 
$\xi_{24}$ has order~3.
Alternatively, we may check that  matrix
$\xi_{4a}  \xi_{4b}$  has trace~1 and eigenvectors
$(0,1,-1,0)$ and  $(0,1,0,-1)$ with eigenvalue~1. So
orthogonality forces the other eigenvalues to 
$\frac{-1 \pm \sqrt{-3}}{2}$ and hence order~3.

\begin{Lemma}
\label{Lemma:invar:LambdaE}    
$\Lambda^E$ is invariant under $\xi_{24a}$ and 
 $\xi_{24b}$.
\end{Lemma}

\fatline{Proof}

Using the isomorphism in Theorem~\ref{Theorem:Q:Leech}	
it is easy to see that $\Lambda^E$ is
generated by the vectors
\begin{align*}
\lambda^{\pm}_{ij} &:
(4_{\SmallTextMbox{on}\,i}, \pm 4_{\SmallTextMbox{on}\,j}, 
0_{\SmallTextMbox{else}}), \qquad i, j \in \tilde{\Omega}
\, , \\
\lambda^E_{d}  &:
(2_{\SmallTextMbox{on}\,d}, \; 0_{\SmallTextMbox{else}} ), 
\qquad d \in \mathcal{C},  \quad \left<d,\omega\right> \; 
\mbox{even}  \, .  
\end{align*}
Here the condition $\left<d,\omega\right>=0 \pmod{2}$ implies
that $\lambda^E_{d}$ has an even number of entries $2$ in
each column and also in row~0 of the MOG.

The operation of $\xi_{24b}$ is negation of row~0 in the MOG.
Thus $\lambda^{\pm}_{ij} - \xi_{24b}(\lambda^{\pm}_{ij})$
has entries divisible by $8$ and 
$\lambda^{E}_{d} - \xi_{24b}(\lambda^{E}_{d})$ has
an even number of entries divisible by 4 in row 0 of
the MOG and zeros elsewhere. 
All these differences are in $\Lambda^E$,
so  $\Lambda^E$ is invariant under $\xi_{24b}$.

Matrix $1-\xi_{4a}$ contains entry $\frac{1}{2}$
everywhere, so the operation of $1-\xi_{4a}$ may be
described as follows:
For each column of the MOG calculate half the sum of its
entries and write the result into each entry of that column.
Performing this operation on  $\lambda^{\pm}_{ij}$ we obtain
either one column with equal entries divisible by 4 or
two columns with equal entries $\pm2$. These results are
all in $\Lambda^E$. 

Any $\lambda^{E}_{d}$ has an even number of entries $2$
in each MOG column and the total number of entries $2$ is
divisible by $4$. Thus there is a vector 
$e \in \Lambda^E$ with an 
even number of nonzero entries, which are all equal to 4,
such that $\lambda^{E}_{d} - e$ has the same number of
entries $2$ and $-2$ in each column of the MOG.
Hence  $\lambda^{E}_{d} - e$ is invariant under
$\xi_{4a}$. We have already shown 
$\xi_{4a}(e) \in \Lambda^E$, thus
$\xi_{24a}(\lambda^{E}_{d} )\in \Lambda^E$.
Hence $\Lambda^E$ is also invariant under  $\xi_{24a}$.

\proofend

%\vspace{1em}

By (\ref{odd:gray:Golay:elem}) and 
Theorem~\ref{Theorem:Q:Leech} we have the following Leech
lattice vectors in MOG coordinates:
\stepcounter{Theorem}
\begin{eqnarray}
\lambda_{g_{0}} =  \small
\begin{array}{|r|r|r|r|r|r|}
\hline   0 & 2 & 2 & 2 & 2 & 2  \\
\hline   2 & 0 & 0 & 0 & 0 & 0  \\
\hline   2 & 0 & 0 & 0 & 0 & 0  \\
\hline   2 & 0 & 0 & 0 & 0 & 0  \\
\hline
\end{array}     \quad , \qquad
\lambda_{\gamma_{0}} =  \small
\begin{array}{|r|r|r|r|r|r|}
\hline  -3 & 1 & 1 & 1 & 1 & 1  \\
\hline   1 & 1 & 1 & 1 & 1 & 1  \\
\hline   1 & 1 & 1 & 1 & 1 & 1  \\
\hline   1 & 1 & 1 & 1 & 1 & 1  \\
\hline
\end{array}     \quad ,
\end{eqnarray}
where $g_0,\ldots,g_5$ is the natural basis of
$\mathcal{G}$  and $\gamma_0,\ldots,\gamma_5$ is the natural basis 
of $\mathcal{G}$. $g_n$ and $\gamma_n$ are obtained from 
$g_0$ and $\gamma_0$ by exchanging MOG column~$0$ with 
column~$n$. We have:
\[
{ \small \left(
    \begin{array}{rrrrrr}
    -3 &  3 &  0 &  1 &  1 & -2 \\
     1 &  1 & -2 &  1 & -1 &  0 \\
     1 &  1 & -2 &  1 & -1 &  0 \\
     1 &  1 & -2 &  1 & -1 &  0 \\    \end{array}
    \right) }^\top \cdot \xi_{4a}  \cdot \xi_{4b}  = 
{ \small \left(
    \begin{array}{rrrrrr}
     3 &  0 & -3 &  1 & -2 &  1 \\
     1 & -2 &  1 & -1 &  0 &  1 \\
     1 & -2 &  1 & -1 &  0 &  1 \\
     1 & -2 &  1 & -1 &  0 &  1 \\    \end{array}
    \right) }^{\bf{\top}} \; ,
\]
and hence
\begin{equation}
\label{eqn:transform_xi24}
  \lambda_{\gamma_i} \stackrel{\xi_{24}}{\longmapsto}
   \lambda_{g_i} - \lambda_{\gamma_i}
  \stackrel{\xi_{24}}{\longmapsto}  
 - \lambda_{g_i} \stackrel{\xi_{24}}{\longmapsto}
     \lambda_{\gamma_i} \; .        
\end{equation}
Since $\Lambda$ is generated by $\Lambda^E$, $g_0$ and
$\gamma_0$, Lemma~\ref{Lemma:invar:LambdaE} and
(\ref{eqn:transform_xi24}) imply:  
\begin{Corollary}
    \label{Corol:invar:Lambda}    
The Leech lattice  $\Lambda$ is invariant under $\xi_{24}$.
\end{Corollary}

\begin{Lemma}
\label{Lemma:xi:xi24}    
The isomorphism 
$Q_x / Z(Q_x) \rightarrow \Lambda / 2\Lambda $ given
by Theorem~~\ref*{Theorem:Q:Leech} maps conjugation with
$\xi$ on $Q_x$ to the automorphism
$\xi_{24}$ of the Leech lattice $\Lambda$ 
(modulo $2 \Lambda$).   
\end{Lemma}    

\fatline{Proof}

We also write $g_i$ for the preimage of  $g_i$ in  $\mathcal{P}$
with positive sign.   
We have $\theta(g_i) = 0$ by Lemma~\ref{Lemma:cocycle:property} and
hence  $\tilde{x}_{g_i} =  x_{g_i}$ by
(\ref{eqn:def:tilde:xd}).
So by (\ref{transf:xi:gray}), $\xi$ operates on
$x_{g_i}$ and $x_{\gamma_i}$  by conjugation 
as follows:
\begin{equation*}
x_{\gamma_i} \stackrel{\xi}{\longmapsto}
\   x_{\pm g_i}  x_{\gamma_i}
\stackrel{\xi}{\longmapsto}  
 x_{g_i} \stackrel{\xi}{\longmapsto}
x_{\gamma_i} \; .        
\end{equation*}
Comparing this operation of $\xi$ with the operation of
of $\xi_{24}$ on $\lambda_{g_i}$ and  $\lambda_{\gamma_i}$
given by  (\ref{eqn:transform_xi24})
we see that these operations are compatible. So by linearity
the operation of $\xi_{24}$ on $\lambda_r$ (modulo $2 \Lambda$)
is the compatible with the operation of  $\xi$ of $x_r$ on
$Q_x$ (modulo the center of $Q_x$) for all $r$ in 
$\mathcal{G} \oplus \mathcal{G}^*$.   

By Lemma~\ref{lemma:gen:coloured} the space $\mathcal{H}^*$
is generated by vectors $ij \in \setF_2^{24}$, where
$i$ and $j$ are in the same column of the MOG, and not in
row~0. By Theorem~\ref*{Theorem:Q:Leech}, the vector
$\lambda_{ij}$ has entries $4$ and $-4$ in the corresponding
positions. Thus $\lambda_{ij}$ is not changed by
$\xi_{24a}$ or by $\xi_{24b}$. By
Lemma~\ref{Lemma:conjugation:xi}  conjugation with 
$\xi_{k}$ or $\xi_{g}$ does not change $x_{ij}$.
So by linearity the operations of  $\xi_{24}$ on 
$\lambda_\delta$ and of $\xi$ on  $x_\delta$ are trivial
for $\delta \in \mathcal{H}^*$ and and hence compatible. 
Hence the operations of $\xi$ on $x_r$ and of $\xi_{24}$ on
$\lambda_r$ are compatible
for all $r \in \mathcal{G} \oplus \mathcal{C}^*$.

For any $h \in \mathcal{P}_\mathcal{H}$ of weight 8 let 
$\hat{x}_h = x_h x_{\gamma(h)}$. Then  $\hat{x}_h$ is invariant 
under $\xi$ by Lemma~\ref{Lemma:x_dh:xi}.
Let  $\hat{\lambda}_h \in \Lambda$ be any representative
of the image of $\hat{x}_h$ under the mapping
$Q_x \rightarrow \Lambda/2\Lambda$ given by
Theorem~\ref*{Theorem:Q:Leech}. Since the elements 
$\hat{x}_h, h \in \mathcal{P}_\mathcal{H}$, $|h|=8$ and 
$ x_{r}, r \in \mathcal{P}_\mathcal{G} \oplus  \mathcal{C}^*$
generate $Q_x$, it suffices to show  
$\xi_{24}(\hat{\lambda}_h) \in \hat{\lambda}_h + 2 \Lambda$.
(Note that 45 of the 64 elements of $\mathcal{H}$ have weight 8; 
so these elements generate  $\mathcal{H}$.)

By  Theorem~\ref*{Theorem:Q:Leech} for every
$h \in \mathcal{P}_\mathcal{H}$ there is a representative
$\lambda_h$ of the image of  $x_h$ with precisely no or 
two nonzero entries of value $2$ in each column of the MOG 
and zeros in row~0. There is also a representative
$\lambda_{\gamma(h)}$ of the the image of $x_{\gamma(h)}$
with entries~4 in row~0 in the columns where   $\lambda_h$ 
does no vanish, and zero entries elsewhere. Thus 
$\hat{\lambda}_h = \lambda_h - \lambda_{\gamma(h)}$ is a
representative of the image of $\hat{x}_h$.

For $\hat{\lambda}_h$ in each column of the MOG  the sum of 
the entries is zero. Thus $\hat{\lambda}_h$ is invariant 
under $\xi_{24a}$. Since $\lambda_h$ has zeros on row~0 of
the MOG and $\lambda_{\gamma(h)}$ has nonzero entries in 
row~0 of the MOG only, $\xi_{24b}$ maps 
$\hat{\lambda}_h$ = $\lambda_h - \lambda_{\gamma(h)}$ to
$\lambda_h + \lambda_{\gamma(h)}$
= $\hat{\lambda}_h  + 2 \lambda_{\gamma(h)}$.

\proofend

An immediate consequence of 
 Corollary~\ref{Corol:invar:Lambda}    
and  Lemma~ \ref{Lemma:xi:xi24}  is:

\begin{Theorem}
	\label{Theorem:xi}
$\xi  \in G(4096_x), \; \xi_{24} \in G(24_x), \;
   \xi \otimes \xi_{24} \in G_{x0}$.
\end{Theorem}

Using \ref{def:xi4ab} and Theorem~\ref{Theorem:xi}, and
identifying the basis vector $i_x$ of $24_x$ with the position of 
entry $i$ in the  MOG, we  obtain the following action
of $\xi$ and $\xi^2$  on the space $24_x$: 

\begin{Corollary}
\label{Corr:op:xi:24x}	
For the $n$-th	column vector $c_n$ in the MOG we have:
\begin{equation*}
c_n^\top \stackrel{\xi}{\longrightarrow} \frac{1}{2} c_n^\top \cdot
\left( \begin{array}{cccc}
\scriptstyle -1 & \scriptstyle -1 & \scriptstyle -1 & \scriptstyle -1 \\
\scriptstyle 1 &\scriptstyle  1 &\scriptstyle  -1 &\scriptstyle  -1 \\
\scriptstyle 1 &\scriptstyle  -1 &\scriptstyle  1 &\scriptstyle  -1 \\
\scriptstyle 1 &\scriptstyle  -1 &\scriptstyle  -1  &\scriptstyle  1 \\
\end{array} \right) 
\; , \quad
c_n^\top \stackrel{\xi^2}{\longrightarrow} \frac{1}{2} c_n^\top \cdot
\left( \begin{array}{cccc}
\scriptstyle -1 & \scriptstyle 1 & \scriptstyle 1 & \scriptstyle 1 \\
\scriptstyle -1 &\scriptstyle  1 &\scriptstyle  -1 &\scriptstyle  -1 \\
\scriptstyle -1 &\scriptstyle  -1 &\scriptstyle  1 &\scriptstyle  -1 \\
\scriptstyle- 1 &\scriptstyle  -1 &\scriptstyle  -1  &\scriptstyle  1 \\
\end{array} \right) 
\; .
\end{equation*}
\end{Corollary}

%%%
%%%%%%%%%%%%%%%%%%%%%%%%%%%%%%%%%%%%%%%%%%%%%%%%%%%%%%%%%%
\subsection{The operation of $\xi$ on $4096_x$} 
%%%%%%%%%%%%%%%%%%%%%%%%%%%%%%%%%%%%%%%%%%%%%%%%%%%%%%%%%%%%%

Let $\xi_\gamma$ and $\xi_g$ be as in (\ref{eqn:def:xi}). We will now 
derive the operation of $\xi = \xi_\gamma \xi_g$ on $4096_x$. Therefore
we first specify a suitable basis of  $4096_x$.

\begin{Definition}
	\label{Def:P_G0}
\hphantom{a}
\begin{StatementList}{\ref{Def:P_G0}}
	\item 
	   Let $g_0, \ldots, g_5$ be the standard basis of the grey 
	   subspace $\mathcal{G}$ of the  Golay code as in 
	   section \ref{subsection:gray:col}.
	   We also write $g_i$ for the preimage of  $g_i$ in  $\mathcal{P}$
	   with positive sign.   
   \item   
	   Let
	   $
	   \mathcal{P}_{\mathcal{G}}^{0} = \left\lbrace 
	   \big( \,  {\textstyle   \prod_{1 \leq i \leq 5} \,  
	   	\tilde{g}_i^{\alpha_i}  } , 0   \,  \big)  
	   \, \middle| \; 
	   {\textstyle    \sum_{1 \leq i \leq 5} \,  {\alpha_i}} 
	   = 0  \pmod{2}
	   \right\rbrace  \; .
	   $ 
	\item 
	    For $d \in \mathcal{G}$ and $\sigma \in \setZ$  let
	    $d_1^{[\sigma]} = d_1^+$ if $\sigma$ is even and 
	    $d_1^{[\sigma]} = d_1^-$ if $\sigma$ is odd.
\end{StatementList}	
\end{Definition}

Then  $\mathcal{P}_{\mathcal{G}}^{0}$ is a subset of 
$ \mathcal{P}_{\mathcal{G}}$ with 16 even elements. 
Any element $e$ of $\mathcal{P}_{\mathcal{G}}$ has a unique 
decomposition $e = \pm \Omega^{\sigma} g_0^{\kappa} d$,
$d \in  \mathcal{P}_{\mathcal{G}}^{0}$, 
$\sigma, \kappa \in \{0,1\}$.
For the unit vectors $f_1^{[\sigma]}$ in $4096_x$, $f \in \mathcal{P}$ 
we have $(\Omega f)_1^{[\sigma]} = (-1)^\sigma f_1^{[\sigma]}$. 
It is easy to see that we have:

\begin{Lemma}
	\label{basis:vectors:4096x}
	The vectors $(g_0^\kappa d h)_1^{[\sigma]}$, 
	 $\sigma, \kappa = 0,1$, 
	$d \in \mathcal{P}_{\mathcal{G}}^{0}$, 
	$h \in \mathcal{P}_{\mathcal{H}}$ form a basis of $4096_x$, 
	with  $w(g_0^\kappa d) = \kappa \pmod{2}$,
	$w_2(g_0^\kappa d) = w_2(d)$, and $A(g_0, d, h) = \theta(g_0,d) =
	\theta(g_0,h) = \theta(d,h) = 0$.
\end{Lemma}

So the the decomposition of a basis vector of $4096_x$  in 
Lemma~\ref{basis:vectors:4096x} is easy to compute and independent 
of the association of the factors given in the Lemma.

\begin{Lemma}
    \label{eqn:transfrom_dpm_xigamma}
Let  $d \in  \mathcal{P}_{\mathcal{G}}^0$  and
    $h \in  \mathcal{P}_{\mathcal{H}}$. Then $\xi_\gamma$ maps
$(g_0^\kappa d h)_1^{[\sigma]}$ to 
$(-1)^{w_2(d) + 1} (g_0^\kappa d h)_1^{[\sigma + \kappa + 1]}$.    
\end{Lemma}

\fatline{Proof}

For any $f \in \mathcal{P}$ put 
 $f' = \frac{1}{2}f_1^+  +  \frac{1}{2}f_1^-$.   From
(\ref{eqn:d:Omega}) we obtain:
\begin{equation}
\label{eqn:transfrom:d_1:d_prime}
    f_1^+ = f' + (\Omega f)' \; , \quad
    f_1^- = f' - (\Omega f) \; , \quad \mbox{and hence} \quad
     f_1^{[\sigma]} =  f' + (-1)^\sigma (\Omega f)'
\end{equation}
Let $e \in \mathcal{G}, \varphi \in \mathcal{G}^* $. 
From  Table~\ref{table:action:196884x} we see that $x_{\varphi}$ maps 
$(eh)'$ to $ (-1)^{\left<eh, \varphi \right>}  (eh)' $.
We have $\left<h, \varphi\right> = 0$.
So by (\ref{eqn:def:xi}) and Lemma~\ref{Lemma:scalprod:gamma},
we obtain:
\[
(eh)' \stackrel{\xi_\gamma}{\longrightarrow}
{\textstyle\frac{1}{8}} (eh)' S(e) \; , \quad \mbox{with} \quad
 S(e) =  \sum_{\varphi \in \mathcal{G}^*}
(-1)^{w_2(\varphi) + \left<e,\varphi\right>} 
=
 \sum_{f \in \mathcal{G}}
(-1)^{w_2(f) + \smallbiscalar{e,f}} \;.
\] 
By Lemma~\ref{Lemma:scalprod:gamma} we have
$
   w_2(f) + \smallbiscalar{e,f} =    w_2(e f)  + w_2(e) 
$
and hence:
\[
  (-1)^{w_2(e)} S(e) = \sum_{f \in \mathcal{G}} (-1)^{w_2(e f)}
     = \sum_{f \in \mathcal{G}} (-1)^{w_2(f)}
     = \sum_{m=0}^6 {\displaystyle \binom{6}{m}}  (-1)^{m(m-1)/2} 
     = -8 \; .
\]
This proves:
\begin{align}
\label{eqn:act:xik:dprime}
(eh)'   \stackrel{\xi_\gamma}{\longmapsto} 
(-1)^{w_2(e) + 1}     (eh)' \; .
\end{align}
We have  $w_2(\Omega e)  =  w_2(e) + w(e) + 1$ 
 by~(\ref{w_2_plus_Omega}) and hence:  
\begin{align}
\label{eqn:w2:sign:Omega}
(\Omega e h)'  \stackrel{\xi_\gamma}{\longmapsto} 
   (-1)^{w_2(e) + w(e)}  (\Omega e h)'   \; .
\end{align}  
From \ref{eqn:transfrom:d_1:d_prime},
(\ref{eqn:act:xik:dprime}) and   (\ref{eqn:w2:sign:Omega}) we obtain:
\[
    (eh)_1^{[\sigma]}  \stackrel{\xi_\gamma}{\longmapsto} 
      (-1)^{w_2(e) + 1} (eh)_1^{[\sigma + \kappa + 1]}  \; , \quad
  \mbox{where $\kappa$ is the parity of $e$}.    
\]
Putting $e = g_0^\kappa d$, the lemma now follows from 
Lemma~\ref{basis:vectors:4096x}.

\proofend

 The following
Lemma states the operation of $\xi_g$ on these basis
vectors:

\begin{Lemma}
	\label{eqn:transfrom_dpm_xig}
Let  $d \in  \mathcal{P}_{\mathcal{G}}^0$  and
$h \in  \mathcal{P}_{\mathcal{H}}$.	Then $\xi_g$ maps 
$(g_0^\kappa dh)_1^{[\sigma]}$ to
\begin{align*}
{\textstyle\frac{1}{4}}   \,
\sum_{e \in \mathcal{P}_{\mathcal{G}}^{0}}  
(-1)^{w_2(de)} (g_0^{\kappa + \sigma + 1} eh)_1^{[\sigma]} \;  .
\end{align*}
\end{Lemma}

\fatline{Proof}

Let $\mathcal{P}_{\mathcal{G}}^+$ =
$\{(\tilde{f},0) \mid f \in {\mathcal{G}} \}$. Then 
$\mathcal{P}_{\mathcal{G}}^0 \subset \mathcal{P}_{\mathcal{G}}^+
\subset \mathcal{P}_{\mathcal{G}} $.
%Let $e, f \in \mathcal{P}_{\mathcal{G}}^+$. 
By (\ref{eqn:def:xi}) and (\ref{relations:Qx:extraspecial})	we have:
\[
\xi_g = 
{\textstyle\frac{1}{8}} \sum_{ e \in \mathcal{P}_{\mathcal{G}}^+}
(-1)^{w_2(fe)} \tilde{x}_{f}  \tilde{x}_{e}  \; , \quad
   \mbox{for any} \quad f \in \mathcal{P}_\mathcal{G}^+ \; .
\]
We have $\left<h, \theta(e)\right> = \left<h, \theta(f)\right>  = 0 $ 
by Lemma~\ref{Lemma:cocycle:property}. Using	
Table~\ref{table:action:196884x} and 
$e^2 = (-1)^{\left<e,\theta(e)\right>}$ we obtain:
\[
    (fh)_1^{[\sigma]}  \stackrel{x_f}{\longrightarrow}
    h_1^{[\sigma]}  \stackrel{x_{\theta(f)}}{\longrightarrow}
      h_1^{[\sigma]} \stackrel{x_e}{\longrightarrow}
    (\bar{e}h)_1^{[\sigma]}  \stackrel{x_{\theta(e)}}{\longrightarrow}
    (-1)^{\left<e, \theta(e)\right>} (\bar{e}h)_1^{[\sigma]}  
      =  (eh)_1^{[\sigma]} \; .
\]
Thus $\tilde{x}_{f}  \tilde{x}_{e}$ maps $ (fh)_1^{[\sigma]}$ to
$ (eh)_1^{[\sigma]}$ and we obtain:
\begin{align}
\label{proof:transfrom_dpm_xig}
(fh)_1^{[\sigma]}
\; & \stackrel{ \xi_g }{\longmapsto}  \; 
{\textstyle\frac{1}{8}}   \,
\sum_{e \in \mathcal{P}_{\mathcal{G}}^+ }  
(-1)^{w_2(f e) } (eh)_1^{[\sigma]}   \; .
\end{align}
To finish the proof, we put $f = g_0^\kappa d$. 

\noindent {Case $\sigma = 1$}

If $fe$ is
odd then $w_2(f e  \Omega) = w_2(fe)$ by~(\ref{w_2_plus_Omega}), 
and we have  $(eh)_1^- = -(\Omega eh)_1^-$ by
(\ref{eqn:d:Omega}), so that the terms for 
$e$ and $\Omega e$ in the sum (\ref{proof:transfrom_dpm_xig}) 
cancel.  If $fe$ is even, these two terms are equal. 
The set $\mathcal{P}_{\mathcal{G}}^{0}$ contains
exactly one of the two elements $d, \Omega d$ for each even
$d \in \mathcal{P}_{\mathcal{G}}^+$, so that summing 
 (\ref{proof:transfrom_dpm_xig})  over the even elements 
of  $\mathcal{P}_{\mathcal{G}}^{+}$ proves the Lemma.
 We remark that $w_2(g_0 d e) = w_2(de)$
by Lemma~\ref{basis:vectors:4096x}, for 
$d, e \in  \mathcal{P}_{\mathcal{G}}^{0}$.

\noindent{Case $\sigma = 0$}

If $fe$ is 
even then $w_2(f e  \Omega) = -w_2(fe)$ by~(\ref{w_2_plus_Omega}), 
and we have  $(eh)_1^+ = (\Omega eh)_1^+$ by
(\ref{eqn:d:Omega}), so that the terms for 
$e$ and $\Omega e$ in the sum (\ref{proof:transfrom_dpm_xig}) 
cancel.  If $fe$ is odd, these two terms are equal. Now the
same argument as in case $\sigma = 0$ shows that
summing over the odd elements of  $\mathcal{P}_{\mathcal{G}}^{+}$
in (\ref{proof:transfrom_dpm_xig}) proves the Lemma.

\proofend

We have $\xi = \xi_\gamma\xi_g $.
Combining Lemmas~\ref{eqn:transfrom_dpm_xigamma} 
and \ref{eqn:transfrom_dpm_xig}, and using
Lemma~\ref{Lemma:scalprod:gamma}  we obtain:

\begin{Corollary}
	\label{eqn:transfrom_dpm_xi}
	Let $d \in \mathcal{P}_{\mathcal{G}}^{0}, 
	h \in \mathcal{P}_{\mathcal{H}}$. Then:	
	\begin{align*}
	(g_0^\kappa dh)_1^{[\sigma]}&
	\; \stackrel{ \xi }{\longmapsto}  \; 
	{\textstyle\frac{1}{4}}   \,
	\sum_{e \in \mathcal{P}_{\mathcal{G}}^{0}}  
    (-1)^{\biscalar{d,e}+w_2(e)+1} 
  	   (g_0^{\sigma} eh)_1^{[\kappa+\sigma+1]} \; ,
  	\\
	(g_0^\kappa dh)_1^{[\sigma]}&
    \; \stackrel{ \xi^2 }{\longmapsto}  \; 
    {\textstyle\frac{1}{4}}   \,
    \sum_{e \in \mathcal{P}_{\mathcal{G}}^{0}}  
    (-1)^{\biscalar{d,e}+w_2(d)+1} 
        (g_0^{\kappa+\sigma+1} eh)_1^{[\kappa]} \; .
	\end{align*}
\end{Corollary}

The following facts are easy to show and of some practical 
value for the implementation.
 
The set ${{\mathcal{G}}^0}$ =
$\{\tilde{e} \mid e \in {\mathcal{P}_{\mathcal{G}}^0}\}$ is a
4-dimensional subspace of $\mathcal{G}$. Any four different elements 
$b_i, \, i = 1, \ldots, 4$, of   ${\mathcal{G}}^0$ with
$w(b_i) = 4$ form a basis of   ${\mathcal{G}}^0$ with 
$\biscalar{b_i,b_j} = \delta_{i,j}$. If $d$ is a sum of $k$
different basis vectors $b_i$ then $w_2(d) =\binom{k}{2} \pmod{2}$.
Thus a suitable basis of the grey part $\mathcal{G}$ of the Golay
code is $(g_0, b_1, b_2, b_3, b_4, \Omega)$.

% !TeX spellcheck = en_GB

\pagebreak

\section{The Griess algebra}
\label{section:Griess:algebra}

This paper is essentially a rewrite of Conway's construction 
\cite{Conway:Construct:Monster} of the Monster $\setM$, 
with an emphasis on providing a computer construction of $\setM$.
In \cite{Conway:Construct:Monster}, the {\em Griess algebra} is defined
as an algebra on the representation $\rho$ of $\setM$, and it is shown
to be invariant under the action of $\setM$.
The prove of the invariance is an elementary, but lengthy calculation.
There are rather subtle differences between our construction
of  $\rho$ and $\setM$ and the corresponding construction in
\cite{Conway:Construct:Monster}, particularly concerning the signs.
In this section we define the Griess algebra in our notation,
and we prove its invariance under the action of  $\setM$

\subsection{Definition of the Griess algebra}
\label{subsection:Griess:algebra}
\stepcounter{Theorem}

Recall from Section~\ref{subsection:196884} that the representation
$196884_x$ decomposes into subspaces as follows:
\begin{equation}
196884_x = 300_x \oplus 98280_x \oplus 98304_x \, ,
\label{eq:Griess:subspaces}
\end{equation}
where $ \, 300_x =  24_x \otimes_{\mbox{\tiny sym}} 24_x  , \;
98304_x = 4096_x \otimes 24_x  ,$
with basis vectors
\begin{equation}
\begin{array}{lrllll}
\mbox{for} & 300_x: & 
(ii)_1 = i_1 \otimes i_1  \, , & i \in \tilde{\Omega}, &  
\mbox{of norm 1}, \\    
& &
(ij)_1 = i_1 \otimes j_1 + j_1 \otimes i_1, 
& i,j \in \tilde{\Omega}, i \neq j, &
  \mbox{of norm } \sqrt{2}, \\    
\mbox{for} & 98280_x: &  X_r &
r \in Q_{x0}, r \;  \mbox{short},
&  \mbox{of norm } 1, \\ 
\mbox{for} & 98304_x: & 
d^\pm \otimes_1 i, &
d \in \mathcal{P}, i \in \tilde{\Omega}, 
&  \mbox{of norm } 1. 
\end{array}
\label{eqn:basis_Griess}
\end{equation}
All these basis vectors are orthogonal, except when equal or
opposite. For $u, v \in 196884_x$ we write $(u, v)$ for the scalar
product of $u$ and $v$ with respect to this norm.
Thus $\left( (ij)_1, (ij)_1 \right)$ is equal to $1$ if $i=j$  and
to $2$ otherwise. The scalar product is invariant under $G_{x0}$
and $\tau$, hence also under $\mathbb{M}$. 

In the sequel we write $A, A', A''$ for
typical elements of $300_x$, which have a natural interpretation as
real $24 \times 24$ matrices acting on the space
$24_x = \Lambda \otimes_{\setZ} \setR$ spanned by the Leech lattice
$\Lambda$. We write $\lambda, \lambda'$ for typical vectors in $24_x$,
and $i_1, j_1, \ldots; i, j, \ldots \in \tilde{\Omega}$ for the 
standard basis vectors of  $24_x$.
\Skip{
 Note that a basis vector $i_1$
of $24_x$  has norm $1/\sqrt{8}$; thus $(i_1, j_1)$ is equal
to $1/8$ if $i=j$ and to $0$ otherwise. Then short vectors in
$\Lambda$ have squared norm 4 in $24_x$ as usual. Symmetric tensor
products of vectors in $24_x$ have a natural interpretation
as vectors in $300_x$. As a consequence, we have e.g.
$( i_1 \otimes i_1 ,  j_1 \otimes j_1)$ =
$64 (i_1, j_1)^2$. 
}
We write  $(\lambda, \lambda')$ for the scalar product
of   $\lambda$ and $\lambda'$;
so $(i_1, j_1)$ is equal to 1 if $i_1 = j_1$ and to 0 otherwise.
Symmetric tensor products of vectors in $24_x$ have a natural
interpretation as vectors in $300_x$; so we have e.g.
$( i_1 \otimes i_1 ,  j_1 \otimes j_1)$ = $(i_1, j_1)^2$. 

We write $q, q'$ for typical elements of $4096_x$. For a short
element $r$ of $Q_{x0}$ we write $q^r$ for the image of $q$
under the action of $r$; and we write $\lambda_r$ for the short
vector in the Leech lattice $\Lambda$ corresponding to $r$.
$\lambda_r$ is defined up to sign only. Short vectors  in
$\Lambda$ have norm 4.
E.g. in case $X_r = X_{ij}$ the vector $\sqrt{8}\lambda_r$ has
co-ordinates 4 and -4 at positions $i$ and $j$, or vice versa;
and zero elsewhere. 

We define a symmetric trilinear $(.,.,.)$ form on $196884_x$,
referred to as the {\em Griess algebra form}, by specifying its
restrictions to the triples of the subspaces of $196884_x$ 
stated in (\ref{eq:Griess:subspaces}). The restrictions of the
form to such triples of subspaces that are nonzero are:  
\begin{align} 
(A, A', A'')
&=
4 \tr(A A' A'') \, 
\label{eq:GriessAlg:AAA} \\
(X_r, X_s, A) &= \left\{\begin{array}{cll}
\pm  \lambda_r A \lambda_r^\top & &
 (\mbox{if} \; \; rs = \pm1) \\
0 & & (\mbox{otherwise})
\end{array} \right. \, ,
\label{eq:GriessAlg:rsA} \\
(X_r, X_s, X_t) &=
\left\{ \begin{array}{cll}
\pm1 & & (\mbox{if} \; \; rst = \pm1)  \\
0  & \phantom{m} & (\mbox{otherwise})
\end{array}
\right. \, ,
\label{eq:GriessAlg:rst} \\
(q \otimes_1 \lambda, q' \otimes_1 \lambda', A)
&= (q, q') 
[\textstyle{\frac{1}{8}}(\lambda, \lambda') \tr A
     + \lambda A \lambda'^\top] \, ,
\label{eq:GriessAlg:qqA} 
\\
(q \otimes_1 \lambda, q' \otimes_1 \lambda', X_r)
&=
 \textstyle{\frac{1}{8}} (q^r, q') [(\lambda, \lambda') -
  2 (\lambda, \lambda_r) (\lambda', \lambda_r) ] \, .
\label{eq:GriessAlg:XYZ} 
\end{align}
Here the form is the same on any permutation of the three arguments.
The restriction of that form to any other triples of such subspaces
is zero. Note that (\ref{eq:GriessAlg:rsA}) and 
(\ref{eq:GriessAlg:XYZ}) are well defined, since changing the sign
of $\lambda_r$ does not change the right-hand side of the equations.
The left-hand side of (\ref{eq:GriessAlg:XYZ}) is symmetric in the
first two arguments, as it should be, since the element
$r$ of $Q_{x0}$ has order two.
\Skip{
The nonzero parts of the Griess algebra form are also shown in the
following figure. Here a  self-loop at a node denoting a subspace $U$
means that the restriction of the form to $U \times U \times U$ is
nonzero. An arrow from subspace $U$ to subspace $V$ means that the
restriction of the form to (any permutation of) 
$U \times U \times V$ is nonzero.

\usetikzlibrary {shapes.misc} 

\begin{figure}[H]
	\centering

\begin{tikzpicture}
% Draw vertices
\node[draw, rounded rectangle] (Z) at (0,0) {$ 98304_x$};
\node[draw, rounded rectangle] (T) at (2,0) {$ 98280_x$};
\node[draw, rounded rectangle] (A) at (4,0) 
  {$ \phantom{2} 300_x  \phantom{.} $ };
\draw[->, >=latex] (Z) -- (T);
\draw[->, >=latex] (T) -- (A);
\draw[->, >=latex] (Z) to[out=290, in=250, looseness=0.5] (A);
\draw (A) to[out=110, in=70, looseness=5] (A);
\draw (T) to[out=110, in=70, looseness=5] (T);

\end{tikzpicture}
    \caption{Griess algebra form at the orbits on
    	$196884_x$ under the action	of $G_{x0}.$ }  
\label{figure:orbits:Gx0}    
\end{figure}
} % end Skip

We define the Griess algebra '$*$' as an algebra on $196884_x$ by
$u * v = w$ for $u, v \in 196884_x $, with $w$ the unique element
of $196884_x$ satisfying $(w, x) = (u, v, x)$ for all
$x \in  196884_x $.

The vectors  $\lambda_r$ (containing a scalar factor $1/\sqrt{8}$)
occur in pairs in (\ref{eq:GriessAlg:rsA}) and \ref{eq:GriessAlg:XYZ};
so the Griess algebra is defined as an algebra on the
$\setZ[\frac{1}{2}]$ lattice spanned by the basis vectors
of $198884_x$.

By definition the Griess algebra form (and hence also the Griess
algebra) is visibly invariant under $G_{x0}$. In the remainder of
this section we will show:

\addtocounter{Theorem}{-1}

\begin{Theorem}
The Griess algebra is invariant under the action of $\mathbb{M}$.
\end{Theorem}

For the proof of the theorem is suffices to show that the
Griess algebra form is invariant
under the action of the triality element $\tau$. In principle, we can
repeat the corresponding argument in \cite{Conway:Construct:Monster} 
almost verbatim, but with (rather subtle) modifications due to
different sign conventions. We will give a self-contained proof here.

%\subsection{Terms of the Griess algebra form on orbits under $N_{xyz0}$}
\subsection{Action of the triality element on the Griess algebra}
\label{subsection:act:Griess}

Let $\mathfrak{V} = \{V_A, V_B, V_C, V_D, V_T, V_X, V_Y, V_Z\}$, where
the elements of  $\mathfrak{V}$ are the subspaces of $196884_x$ in 
Table~\ref{table:dictionary} in Section~\ref{subsection:action:tau}.
Each of these subspaces in invariant under the action of the subgroup
$N_{xyz0}$ of $G_{x0}$ and $N_{x0}$. The group $N_0$ contains
the triality element $\tau$ and permutes the
subspaces $V_A, V_B, V_C$, and also the subspaces $V_X, V_Y, V_Z$;
and it fixes  $V_D$ and $V_T$.
\Skip{
 More specifically,
$\tau$ cyclically exchanges the ordered pairs
$(V_A, V_X)$, $(V_B, V_Y)$, and $(V_C, V_Z)$ of subspaces;
and an element $x_\delta, \delta \, \mbox{odd}$, exchanges
$(V_B, V_Y)$ with $(V_C, V_Z)$, and fixes $V_A$ and $V_X$.
}

It suffices to show the
invariance of the Griess algebra form under the action of $\tau$
for the basis vectors of these subspaces. Here we take
$X_{ij} - X^+_{ij}$ and $X_{ij} + X^+_{ij}$ as basis vectors of the
subspaces $V_B$ and $V_C$, respectively. Thus all basis vectors
of the spaces $V_A, V_B, V_C$ have norm $\sqrt{2}$.

The subgroup $R$ of $N_{xyz0}$ generated by the involutions
$x_{-1}, y_{-1}, z_{-1}, x_\delta, \delta \in \mathcal{C^*}$,
$\delta \, \mbox{even},$ is an elementary Abelian 2 group of order
$2^{13}$. Hence it defines a grading of the Griess algebra by
decomposing the space $196884_x$ into the $2^{13} $ common
eigenspaces under the actions of the elements of $R$. 
The correspondence between an invariant elementary Abelian 2
group and such a grading on the Griess algebra is just a
multidimensional version of statement (4) after the first
paragraph in \cite{Conway:Construct:Monster}, Appendix~5.   
The group $R$ is isomorphic to the additive group of the vector space 
$\setF_2^2 \oplus \{\delta \in\mathcal{C^*} | \delta \, \mbox{even}\}$
over $\setF_2$, with the elements of the first factor denoted by
$\{0,x,y,z\}$. The grade of such a subspace can be labelled with an
element of the Pontryagin dual $\hat{R}$ of $R$.
So  $\hat{R}$  is the group of the one-dimensional characters of $R$,
and isomorphic to the additive group of the dual space
$\setF_2 \oplus (\mathcal{C} / \{0, \tilde{\Omega}\})$ of the space
mentioned above. Here the scalar product on
$\setF_2^2 \times \setF_2^2$ given by the determinant,
which identifies $\setF_2^2$ with its dual space. Since
$\{0, \tilde{\Omega}\} \cong \mathcal{P} / Z(\mathcal{P})$, we may
denote a grade by $(\alpha, d + \mathcal{P} / Z(\mathcal{P}))$
for $\alpha \in\{0,x,y,z\},  d \in \mathcal{P}$, which we will
abbreviate to $(\alpha, d)$. So e.g. the subspace of grade $(y, d)$  
is the subspace of $196884_x$  fixed by $y_{-1}$, and negated by
$x_{-1}, z_{-1}$, on which every even $\delta \in \mathcal{C^*}$
operates as multiplication by
$(-1)^{\langle d, \delta \rangle}$. 
It turns out that all basis vectors and also their images under
$\tau^{\pm1}$ have a defined grade. More specifically, we
obtain from Table~\ref{table:action:196884x}:
\[
\begin{array}{c|c|c|c|c|c}
\mbox{vector or subspace} & V_A, V_B, V_C, V_D
  & X^+_{d \cdot\delta}, \delta \, \mbox{even} &
   X^+_{d \cdot i} & d ^- \otimes_1 i & d ^+ \otimes_1 i \\
\hline
  \mbox{grade} &  (0,0) & (0,d)  & (x,d) & (y,d) & (z,d)
\end{array} \; .
\]
%\stepcounter{Theorem}
\begin{Lemma}
\label{Griess:grading:zero}	
	\hphantom{a}
	\label{Lemma:Griess:grade}
	\begin{StatementList}{\ref{Lemma:Griess:grade}}
		\item
		The triality element $\tau$ cyclically exchanges the 
		subspaces with grades $(x, d)$, $(y, d)$, and $(z, d)$; 
		and it fixes the subspaces with grades $(0,d)$.
		\item
        If $u, v, w \in 196884_x$ all have a defined grade and
        $(u,v,w) \neq 0$ then the sum of the grades of
        $u$, $v$, and $w$ is $(0, 0)$.		
		\item
        If $u, v, w \in 196884_x$ all have a defined grade then
        the sum of the grades of $u^\tau$, $v^\tau$ and
        $w^\tau$ is $(0, 0)$ if and only if the sum of the
        grades of $u$, $v$ and $w$ is $(0, 0)$.   		
	\end{StatementList}
\end{Lemma}

\fatline{Proof}

(\ref{Lemma:Griess:grade}.1)  follows from (\ref{triality:ii}),
(\ref{triality:ij}), (\ref{triality:Xdi}), and
(\ref{eqn:operations:V_T}). If linear transformation $g$ on
$196884_x$ fixes or negates the vectors $u, v, w \in 196884_x$,
and $(u,v,w) \neq 0$,  then $g$ negates an even number of the
vectors $u,v,w$. Applying this to all $g \in R$ yields
(\ref{Lemma:Griess:grade}.2).
(\ref{Lemma:Griess:grade}.3) is a consequence of
(\ref{Lemma:Griess:grade}.1) and (\ref{Lemma:Griess:grade}.2) . 

\proofend

For $u, v, w \in 196884_x$ and $g \in \mathbb{M}$ we abbreviate 
$(u^g, v^g, w^g)$ to $(u, v, w)^g$.
For a property $Q$ we put $\Delta_Q = 1$ if $Q$ is true, and
$\Delta_Q = 0$ if $Q$ is false. So for the Kronecker delta we
write  $\Delta_{i=j}$ instead of the more common notation 
$\delta_{i,j}$. 

In the next two sections we will compute the Griess algebra form on
triples of the basis vectors of the spaces in $\mathfrak{V}$ and
on their images under $\tau^{\pm 1}$. This will show the invariance
of that form under $\tau$. As in \cite{Conway:Construct:Monster}
we may reduce the number of cases by observing: 
%\stepcounter{Theorem}
%\renewcommand{\labelenumi}{(\arabic{section}.\arabic{Theorem}.\arabic{enumi})}
\begin{enumerate}
	\item[-]
	We have $(u,v,w) = (u',v', w')$ for any permutation
	$(u',v', w')$ of $(u,v,w)$.   
	\item[-]
    It suffices to consider the action of $\tau^{\pm 1}$ on one
    of the triples $(u,v,w)$, $(u,v,w)^{\tau^{\pm 1}}$. 
	\item[-]
	Since $\tau$ normalizes the subgroup $N_{xyz0}$ of $G_{x0}$,
	we may replace the case $(u,v,w)$ by the case  $(u^g,v^g,w^g)$
	for any $g \in N_{xyz0}$.
\end{enumerate}
We wil use the phrase {\em by triality} if we reduce the number of
cases by the first or second of these methods, and the phrase 
{\em by symmetry} if we replace $(u,v,w)$ by $(u^g,v^g,w^g)$,
$g \in N_{xyz0}$.

For $r \in Q_{x0}$ we define
$x_r^+$ to be $x_{\Omega \cdot r}$ if this is short and
$x_r$ otherwise, as in Section~\ref{subcection:short:vectors}.
If $X_r^+$ is a basis vector of $98280_x$ then
$x_r^+$ is the short element of $Q_{x0}$ corresponding to $X_r^+$.

\subsection{Cases involving any of the subspaces $V_X, V_Y, V_Z$}

\newcommand{\dstackrel}[3]{\stackrel{\stackrel{\scriptstyle #1}{#2}}{#3}}
\stepcounter{Theorem}

By triality we may assume that element $u$ in the product $(u, v, w)$
is in $V_X$. Then by Lemma~\ref{Griess:grading:zero} the nonzero cases
are	 $XYZ$ and $XXU$, $U \in \{A,B,C, D, T\}$.

\subsubsection{Case $XYZ$}
We show the invariance of $(u,v,w)$ =
$(X_{d \cdot i}, e^- \otimes_1 j, f^+ \otimes_1 k)$ under the action
of $\tau^{\pm 1}$ for $d, e, f \in \mathcal{P}$ and basis vectors
$i, j, k \in \tilde{\Omega}$. By Lemma~\ref{Griess:grading:zero} 
this expression and its images under $\tau^{\pm 1}$ are nonzero
in case $def \in  \{\pm1, \pm \Omega\}$ only.
% We assume that this is the case.
In this case  element  $y_{\bar{f}} z_{\bar{e}}$ of $N_{xyz0}$
maps $e^-$ to $1^-$, $f^+$ to $1^+$, and $d$ to an element of
$\{\pm1, \pm \Omega\}$. Since $\lambda_{d,i}$ is short, we must
have $d = \pm1$. So by symmetry it suffices to consider the case
$(X_{i}, 1^- \otimes_1 j, 1^+ \otimes_1 k)$. We have
\begin{align*}
8(X_{i}, 1^- \otimes_1 j, 1^+ \otimes_1 k)^\tau
& \stackrel{(\ref{triality:Xdi})}{=} 
8(1^- \otimes_1 i, 1^+ \otimes_1 j, X_k)
\stackrel{(\ref*{eq:GriessAlg:XYZ})}{=}
(i,j) - 2(i, \lambda_k)(j, \lambda_k) \, ,
\\
8(X_{i}, 1^- \otimes_1 j, 1^+ \otimes_1 k)^{\tau^{-1}}
& \stackrel{(\ref{triality:Xdi})}{=} 
8(1^+ \otimes_1  i, X_j, 1^- \otimes_1 k)
\stackrel{(\ref*{eq:GriessAlg:XYZ})}{=}
(i,k) - 2(i, \lambda_j)(k, \lambda_j) \, ,
\\
8(X_{i}, 1^- \otimes_1 j, 1^+ \otimes_1 k)
& \stackrel{(\ref*{eq:GriessAlg:XYZ})}{=}
(j,k) - 2(j, \lambda_i)(k, \lambda_i) \, .
\end{align*}
We have $(i,j) = \Delta_{i=j}$. One of the vectors
$\pm\sqrt{8}\lambda_j$ has entry 3 at the co-ordinate labelled
by $i$ and entry $-1$ elsewhere. Hence
$(i,\lambda_k)(j,\lambda_k)$ = 
$\frac{1}{8} (4\Delta_{i=k} - 1) (4\Delta_{j=k} - 1)$.
So we obtain:
\[
(i,j) - 2(i, \lambda_k)(j, \lambda_k)  =
  -{\textstyle\frac{1}{4}} 
  +\Delta_{i=j} +\Delta_{i=k} + \Delta_{j=k}
  - 4\Delta_{i=j=k} \; ,
\]
 which is symmetric under all permutations of $i$, $j$, and $k$.

\proofend

\subsubsection{Case $XXT$}

We show the invariance of $(u,v,w)$ =
$(X_{e\cdot i}, X_{f\cdot j}, X^+_{d\cdot \delta})$ under the action of
$\tau$ for $d, e, f \in \mathcal{P}$, $|d|=8$, $\delta \subset d$ even,
and $i, j \in \tilde{\Omega}$. The element  $y_{\bar{f}}$ of $N_{xyz0}$
maps $X_{f\cdot j}$ to $X_j$. By Lemma~\ref{Griess:grading:zero}
the product $(X_{e\cdot i}, X_{j}, X_{d\cdot \delta})$ or its images under
$\tau^{\pm1}$ are nonzero in case $de \in \{\pm1, \pm\Omega\}$
only. So by symmetry it suffices to deal with the case
$(u,v,w)$ = $(X^+_{d\cdot i}, X_{j}, X^+_{d\cdot \delta})$, $\delta$ even.
Evaluating the '$+$' exponents we obtain:
\Skip{
$x_{d\cdot i}$ is short; hence $\langle d, i\rangle = 0$, i.e. 
$i \notin d$.
$x_{d\cdot \delta}$ is also short,  so by (\ref{eqn:short:vectors})
we have  $\delta \in  A(d,\mathcal{C})$; and $|d|= 8$ if
$|\delta/2|$ is even,  $|d|= 16$ if
$|\delta/2|$ is odd.  
Thus the relevant cases for $(u,v,w)$ are:
}
\begin{align}
(X_{d\Omega^{\langle d,  i \rangle} \cdot i}, X_{j},
 X_{d\Omega^{|\delta|/2} \cdot \delta})
& \stackrel {(\ref{eq:GriessAlg:rst})}{=} \Delta_{ij = \delta}
  \Delta_{\langle d, i \rangle = |ij|/2} \, ; \;
\mbox{for} \; |d| = 8, \; \delta \in  A(d,\mathcal{C})\, .
 \label{Griess:case:rd:octad}
\end{align}
We have
\begin{align}
\noalign{
$(X_{d\Omega^{\langle d, i \rangle}\cdot i}, X_{j},
X_{d\Omega^{|\delta|/2}\cdot \delta})^\tau$  
} \nonumber \\
&\dstackrel{(\ref{triality:Xdi}),}{(\ref*{eqn:operations:V_T})\phantom{,}}{=}
{\frac{1}{8}} { \sum_{\epsilon \in A(d,\mathcal{C})}}
 (-1)^{|\delta \cap \epsilon| + |\epsilon|/2 } 
\left(
(d \Omega^{\langle d,i \rangle })^- \otimes_1 i, 1^-\otimes_1 j,
 X_{d \Omega^{|\epsilon|/2} \cdot \epsilon} \right) \nonumber\\
& \stackrel{(\ref*{eq:GriessAlg:XYZ})}{=}
\frac{1}{64} 
{ \sum_{\epsilon \in A(d,\mathcal{C})}}
 (-1)^{|\delta \cap \epsilon| + |\epsilon|/2 } \, 
(-1)^{{\langle d, i \rangle} + |\epsilon|/2 }\,
\left((i,j) -2 (i, \lambda_{d,\epsilon})  (j, \lambda_{d,\epsilon})
 \right) \nonumber \\
&\phantom{m} =  \phantom{m} \frac{(-1)^{\langle d, i \rangle}}{64} {
\sum_{\epsilon \in A(d,\mathcal{C})}}
\, (-1)^{|\delta \cap \epsilon|} \,  
\Delta_{i=j} - (-1)^{|ij\delta \cap \epsilon|}
 \Delta_{\{i,j\} \subset d} \nonumber  \\
& \qquad \qquad \qquad
( {\scriptstyle \mbox{\scriptsize since} \, (i,j) = \Delta_{i=j} , \;
 (i, \lambda_{d,\epsilon}) 
 = \sqrt{1/2} \cdot (-1)^{i \cap \epsilon}  \nonumber	
}) \\
&\phantom{m} =  \phantom{m} (-1)^{\langle d, i \rangle} \left(
 \Delta_{\delta=0} \Delta_{i=j} 
 - \Delta_{\delta=ij} \Delta_{\{i,j\} \subset d}
\right)	\; \nonumber \\
& \phantom{m} =  \phantom{m} \left\{
\begin{array}{ccl}
1 & & \mbox{if } i = j; \, \delta = 0; \; \{ i, j \} \not \subset d  \\
1 & & \mbox{if } i \neq j; \, \delta = ij; \{ i, j \} \subset d  \\
0 & & \mbox{otherwise}
\end{array}	
\right.   \nonumber \\
& \phantom{m} =  \phantom{m} \Delta_{\delta=ij} 
  \Delta_{\langle d, i \rangle = |ij|/2} \quad
({\scriptstyle \mbox{\scriptsize since }
	\langle d, ij \rangle =0  \mbox{\scriptsize \ in case }
\delta = ij}) \, . 
\label{Griess:case:rd:octad_sum1} 
\end{align}
For the computation of $\lambda_{d, \epsilon}$ (up to sign) we
remark that the co-ordinate  of $\sqrt{8}\lambda_{d, \epsilon}$
at position $i$ is equal to -2 if $i \in \epsilon$, to 2 if
$i  \in d \setminus \epsilon$, and to 0 otherwise.
\Skip{
As in \cite{Conway:Construct:Monster},
it turns out that the nonzero cases are:
\[
\begin{array}{ccccc}
\delta  = ij & \land & i = j & \land & i \notin d \, , \\ 
\delta  = ij & \land & i \neq j & \land  &i,j \in d   \, .
\end{array}
\]
}
To finish this case, we also have deal with the operation of
$\tau^{-1}$. A computation similar to 
\ref{Griess:case:rd:octad_sum1} yields:
\begin{align}
\noalign{
	$(X_{d\Omega^{\langle d, i \rangle}\cdot i}, X_{j},
	X_{d\Omega^{|\delta|/2}\cdot \delta})^{\tau^{-1}}$  
} \nonumber \\
&\dstackrel{(\ref{triality:Xdi},(\ref{eqn:operations:V_T}),}
{(\ref*{eq:GriessAlg:XYZ})}{=}
\frac{1}{64} 
{ \sum_{\epsilon \in A(d,\mathcal{C})}}
 (-1)^{ |\delta|/2 + |\delta \cap \epsilon| } \, 
\left((i,j) -2 (i, \lambda_{d,\epsilon})  (j, \lambda_{d,\epsilon})
\right) \nonumber \\
&\stackrel{(\ref{Griess:case:rd:octad_sum1})}
{\phantom{mm} = \phantom{mm}}
(-1)^{ |\delta|/2 + {\langle d, i \rangle}}
 \Delta_{\delta=ij} 
\Delta_{\langle d, i \rangle = |ij|/2} \; \nonumber \\
&\phantom{mmi} = \phantom{mmi} \Delta_{\delta=ij} 
  \Delta_{\langle d, i \rangle = |ij|/2} \, , \quad
( {\scriptstyle \mbox{\scriptsize since} \, ij = \delta \subset d
	\mbox{\ \scriptsize implies}
	(-1)^{ |\delta|/2 + \langle d, i \rangle} =1} ) \;.
\label{Griess:case:rd:octad_sum2}
\end{align}
\proofend 

\subsubsection{Cases $XXA, XXB, XXC$}

We show the invariance of $(u,v,w)$ with 
$u = X_{d\cdot i}$, $v = X_{e\cdot j}$,
$w \in \{(k,l)_1, X_{kl} \pm X_{\Omega\cdot kl}\}$, $k \neq l$ 
under the action of $\tau^{\pm1}$.  The element
$y_{\bar{d}}$ of $N_{xyz0}$
maps $X_{d\cdot i}$ to $X_i$. By Lemma~\ref{Griess:grading:zero}
the product $(X_{i}, X_{e\cdot j}, X_{kl})$ or its images under
$\tau^{\pm1}$ are nonzero in case $e \in \{\pm1\}$ only.
(Note that $X_{\pm \Omega\cdot  j}$ is not short.)
 So by symmetry it suffices to deal with the cases
$u = X_{i}$, $v = X_{j}$, and $w$ as above. 
By simplifying (\ref{triality:Xdi}) we obtain:
\begin{align}
X_{i} 
\, \stackrel{\tau}{\longmapsto} \,
 1^- \otimes_1 i
\, \stackrel{\tau}{\longmapsto} \,
1^+ \otimes_1 i
\, \stackrel{\tau}{\longmapsto} \,
X_{i} \, .
\label{triality:1di} 
\end{align}
The following statement is just (\ref{triality:ij}):
\begin{align}
\label{triality:ij:copy}
(ij)_1  \, & \stackrel{\tau}{\longmapsto} \,
X_{ij} - X_{ij}^+   \, \stackrel{\tau}{\longmapsto} \,
X_{ij} + X_{ij}^+   \, \stackrel{\tau}{\longmapsto} \, (ij)_1  \, ,
\quad i \neq j  \; .
\end{align}
To finish this case, we just have to compute the operation of
$\tau^{\pm 1}$ on $(X_i, X_j, w)$, where $w$ is any of the
terms occurring in (\ref{triality:ij:copy}).
Put $d^{[0]} = d^+, d^{[1]} = d^-$. 

\vspace{1ex}
\fatline{Subcase $i \neq j, \; \{i, j\} \neq \{k, l\}$ }

From (\ref{eq:GriessAlg:rsA})--(\ref{eq:GriessAlg:XYZ})
we see that $(X_i,X_j,w)^g = 0$ holds for $g = 1, \tau^{\pm 1}$  
in this case. 

\vspace{1ex}
\fatline{Subcase $i  \neq j, \; \{i, j\} = \{k, l\}$ }

The relevant
calculations can be done by using (\ref{triality:1di}),
(\ref{triality:ij:copy}) and the following formulas:
\begin{align}
\left(X_i, X_j, (ij)_1 \right) 
&\stackrel{(\ref{eq:GriessAlg:rsA})}{=} 0 ,
\nonumber \\
\left(X_i, X_j, X_{ij} \pm X_{\Omega\cdot ij}  
\right) & 
\stackrel{(\ref{eq:GriessAlg:rst} )}{=}
1,
\nonumber \\
\left(
1^{[\sigma]} \otimes_1 i,
1 ^{[\sigma]} \otimes_1 j,
(ij)_1 
\right) & 
\stackrel{(\ref{eq:GriessAlg:qqA} )}{=}
1 , 
\nonumber \\
\left(
1^{[\sigma]} \otimes_1 i,
1 ^{[\sigma]} \otimes_1 j,
X_{ij} + (-1)^\nu X_{\Omega\cdot  ij} 
\right) & 
\stackrel{(\ref{eq:GriessAlg:XYZ} )}{=}
 \Delta_{\sigma \neq \nu},
 \qquad\mbox{for} \; i \neq j \; .
\nonumber  
%\label{eq:Griess:case:delta:2}
\end{align}
Put $X_r = X_{\Omega^\mu\cdot ij}$.
For the last of these formulas we remark that one of the
vectors $\pm\sqrt{8}\lambda_{\Omega^\nu,ij}$ has co-ordinate $-4$
at position $i$, co-ordinate $4 \cdot (-1)^\nu$ at position $j$,
and zero elsewhere; hence
$(i, \lambda_r )(j, \lambda_r) = 4 \cdot (-1)^{\mu}$.
We also have $(1^{[\sigma]})^r = (-1)^{\sigma \mu} \cdot 1^{[\sigma]}$,
i.e $(1^{[\sigma]},(1^{[\sigma]})^r) =  (-1)^{\sigma \mu}$.

\vspace{1ex}
\fatline{Subcase $i =j$, $k \neq l$ }
 
The relevant
calculations can be done by using (\ref{triality:1di}),
(\ref{triality:ij:copy}) and the following formulas:
\begin{align}
\left(X_i, X_i, (kl)_1 \right) 
&\stackrel{(\ref{eq:GriessAlg:rsA})}{=}
{\textstyle \frac{1}{4}} - \Delta_{i \in \{k,l\} } ,
\nonumber \\
\left(X_i, X_i, X_{kl} \pm X_{\Omega\cdot  kl}  
\right) & 
\stackrel{(\ref{eq:GriessAlg:rst} )}{=}
0,
\nonumber \\
\left(
1^{[\sigma]} \otimes_1 i,
1 ^{[\sigma]} \otimes_1 i,
(kl)_1 
\right) & 
\stackrel{(\ref{eq:GriessAlg:qqA})}{=}
0, 
\nonumber \\
\left(
1^{[\sigma]} \otimes_1 i,
1 ^{[\sigma]} \otimes_1 i,
X_{kl} + (-1)^\nu X_{\Omega\cdot  kl} 
\right) & 
\stackrel{(\ref{eq:GriessAlg:XYZ} )}{=}
  \Delta_{\sigma = \nu}
   \left({\textstyle \frac{1}{4}} - \Delta_{i \in \{k,l\} }\right)
\, , \qquad\mbox{for} \; k \neq l \; .
\nonumber  
\end{align}
For the first of these formulas we remark that one of the
vectors $\pm \sqrt{8}\lambda_i$ has co-ordinate $-3$ at
position $i$ and $1$ elsewhere. For the last of these formulas a similar
calculation as in the previous subcase yields
$(i,\lambda_r)^2 = 2 \Delta_{i\in \{k,l\}}$,
$(1^{[\sigma]},(1^{[\sigma]})^r) =  (-1)^{\sigma \mu}$, 
with $X_r = X_{\Omega^\mu\cdot kl}$. 
 
\proofend

\subsubsection{Case $XXD$}

We show the invariance of $(u,v,w)$ with 
$u = X_{d\cdot i}$, $v = X_{e\cdot j}$, $w = (k,k)_1$ under the action
of $\tau^{\pm1}$. As in the previous cases it suffices to deal with
$u = X_{i}$, $v = X_{j}$, and $w$ as above. In case $i \neq j$
we obtain  $(u,v,w)=  (u,v,w)^{\tau^{\pm1}} = 0$ from
(\ref{eq:GriessAlg:rsA}) and (\ref{eq:GriessAlg:qqA}).
A similar computation as in the previous case yields:
\begin{align}
\left(X_i, X_i, (jj)_1 \right) 
&\stackrel{(\ref{eq:GriessAlg:rsA})}{=}
{\textstyle \frac{1}{8}} + \Delta_{i= j } \, ,
\nonumber \\
\left(
1^{[\sigma]} \otimes_1 i,
1 ^{[\sigma]} \otimes_1 i,
(jj)_1 
\right) & 
\stackrel{(\ref{eq:GriessAlg:qqA})}{=}
{\textstyle \frac{1}{8}} + \Delta_{i= j } \, .
\nonumber
\end{align}
 
\proofend

\subsection{Cases involving subspace $V_T$, but not subspaces $V_X, V_Y, V_Z$}

\stepcounter{Theorem}

Put $V_{ABCD} = V_A\oplus V_B\oplus V_C\oplus V_D$, and
$V_+ = V_{ABCD} \oplus V_T$.
We decompose the triality element $\tau$ into a product
$\tau = y_\tau x_\tau$ of involutions $x_\tau$ and $y_\tau$
acting on  $V^+$.
Here the action of $x_\tau$ and $y_\tau$ on $V_T$ is given
by (\ref{eqn:operations:V_T}). $x_\tau$ and $y_\tau$ act on
$V_{ABCD}$ in the same way as 
$x_i$ and $y_i$ for any $i \in \tilde{\Omega}$, respectively.
From Table~\ref{table:action:196884x} we see that the action	   
of $x_i$ and $y_i$ on $V_{ABCD}$ is the same for all
$i \in \tilde{\Omega}$, and we obtain:
\begin{align}
x_\tau \, \mbox{exchanges} \; X_{ij} +  X_{\Omega \cdot ij}  \;
\mbox{with} \; X_{ij} -  X_{\Omega \cdot ij},  \; \mbox{and fixes} \;
(ij)_1 \; \mbox{and} \; (ii)_1 \; , 
\label{eqn:Griess:TTA:xtau} \\
y_\tau \, \mbox{exchanges} \; (ij)_1  \;
\mbox{with} \; X_{ij} +  X_{\Omega \cdot ij},  \; \mbox{and fixes} \;
X_{ij} -  X_{\Omega \cdot ij} \; \mbox{and} \; (ii)_1 \; . 
\label{eqn:Griess:TTA:ytau}
\end{align} 

\Skip{
\medskip
For these cases we use a grading of $V_{+}$ onto $\mathcal{C}$
different from that in Lemma~\ref{Griess:grading:zero}:	
\begin{equation} 
\label{eqn:grading:vplus}       
\begin{array}{c|c|c|c|c}
\mbox{vector or subspace} & V_A,  V_D & X_{ij} &  X^+_{ij}
& X^+_{d \cdot\delta}, |d| = 8, \delta \subset d \, \mbox{ even}  \\
\hline
\mbox{grade} &  0 & 0 & \Omega & d \Omega^{|\delta|/2} 
\end{array} \; .
\end{equation}        
\begin{Lemma}
	%\label{Vplus:grading:zero}	
	%\hphantom{a}
	\label{Lemma:Vplus:grade}
	    Let the grading on $V_+$ with values in $\mathcal{C}$  be
	    given by (\ref{eqn:grading:vplus}).  
		If $u, v, w \in 196884_x$ all have a defined grade and
		$(u,v,w) \neq 0$ then the sum of the grades of
		$u$, $v$, and $w$ is $0$.		
\end{Lemma}

\fatline{Proof}

The elementary Abelian 2  group $\{x_\delta, \mid \delta \in  C^*\}$
acts on $V^+$ and lets the restriction of the Griess algebra
product to $V^+$ invariant. So the same argument as in the proof
of Lemma~\ref{Griess:grading:zero} shows that $V_+$ has a
grading with values in the Pontryagin dual of  $\mathcal{C}^*$,
which is $\mathcal{C}$. The grades in (\ref{eqn:grading:vplus})
have been  computed from Table~\ref{table:action:196884x}.

\proofend
}

\subsubsection{Case $TTT$}
 
We show the invariance of $(u,v,w)$ =
$(X^+_{d\cdot \delta}, X^+_{e\cdot \epsilon}, X^+_{f\cdot \varphi})$
for octads $d, e, f \in \mathcal{P}$ and
$\delta \in A(d,\mathcal{C})$, $\epsilon \in A(e,\mathcal{C})$, 
$\varphi \in A(f,\mathcal{C})$,   under the action of
$x_\tau$ and $y_\tau$. By Lemma~\ref{Lemma:Griess:grade} we have
for any $n \in \setZ$:
\begin{equation}
\label{eqn:grade:TTT}
 (X^+_{d\cdot \delta}, X^+_{e\cdot \epsilon}, X^+_{f\cdot \varphi})^{\tau^n} 
 \neq 0 \; \implies  def \in \{\pm 1,  \pm \Omega\}   \, .
\end{equation}
With $d, e, f, \delta, \epsilon ,\varphi$ as above we obtain:
\begin{Lemma}
\label{lemma:Griess:TTT}	
\[	
 (X^+_{d\cdot \delta}, X^+_{e\cdot \epsilon}, X^+_{f\cdot \varphi}) =
 \left\{
 \begin{array}{cl}
    \pm(-1)^{ \left<d, \epsilon \right>} & \mbox{ if }
        def = \pm1, \; \mbox{ and }  \; \delta \epsilon \varphi = \omega 
           \pmod{\mathcal{C}}\\
     0 & \mbox{ otherwise }     
 \end{array}
 \right. \; ;
\] 	
where $\omega = d \cap e \in \mathcal{C}^*$. If that product is nonzero,
we also have  $\omega = d\cap f = e \cap f \pmod{\mathcal{C}} $, 
$ |\delta|/ 2 + |\epsilon|/ 2 + |\varphi|/ 2 = 0 \pmod{2}$, and
$\left<d, \epsilon\right> = \left<f, \epsilon\right>$ =
$\left<d, \phi\right> = \left<e, \phi\right>$ =
$\left<e, \delta\right> = \left<f, \delta\right>$.
\end{Lemma}
\fatline{Proof}

By (\ref{eqn:grade:TTT}) it suffices to consider the case
$def = \Omega^m$. A calculation shows:
\begin{equation}
\label{eqn:grade:TTT01}
x^+_{d\cdot \delta} x^+_{e\cdot \epsilon} x^+_{f\cdot \varphi} =
x_{\omega \delta \epsilon \varphi}  x_{-1}^\sigma x_\Omega^\nu\, , 
\quad \mbox{ if } \; def =  \Omega^m \, , \; \mbox{ where } 
\end{equation}
\[
\omega = d \cap e = d\cap f = e \cap f \in \mathcal{C}^*, \;
\sigma  = \langle d, \epsilon \rangle , \;
\nu = m + |\delta|/ 2 + |\epsilon|/ 2 + |\varphi|/ 2 \, .
\]

In case $def = \Omega$ we have $\omega = 0$ and the vectors
$\delta, \epsilon, \varphi$ have
disjoint support; hence  $|\omega\delta \epsilon \varphi| = 2 \pmod{4}$
by (\ref{eqn:grade:TTT01}), excluding the possibility  
$\omega\delta \epsilon \varphi = 0 \pmod{\mathcal{C}}$. Thus
(\ref{eq:GriessAlg:rst}) implies
$ (X^+_{d\cdot \delta}, X^+_{e\cdot \epsilon}, X^+_{f\cdot \varphi}) = 0$
in case $def = \pm \Omega$.

So we may assume $def=1$. We also assume 
$\delta \epsilon \varphi = \omega$, since otherwise the product in the
lemma is zero by (\ref{eq:GriessAlg:rst}  and  (\ref{eqn:grade:TTT01}).
Put $\omega_d = e \cap f, \omega_e = d\cap f, \omega_f = d\cap f$,
and $\hat{\omega} = d \cup e \cup f$. Since $d, e, f$ are octads,
we have $|\omega_d| = |\omega_e| = |\omega_f| = 4$; and
$\hat{\omega}$ is the disjoint union of  $\omega_d$, $\omega_e$,
and  $\omega_f$. The representatives of the Golay cocode element
$\omega$ that are subsets of  $\hat{\omega}$ are
 $\omega_d$, $\omega_e$, $\omega_f$, and $\hat{\omega}$.

Considering $\delta, \epsilon$ and $\varphi$ as subsets of
$d, e$, and $f$, respectively, we have 
$\delta+ \epsilon +\varphi \in
\{\omega_d, \omega_e, \omega_f, \hat{\omega} \}$.
By changing at most one of the
sets $\delta, \epsilon$, or $\varphi$ to its complement in
$d, e$, or $f$, we may assume
$\delta + \epsilon +\varphi = \hat{\omega}$. Since $\delta \subset d$,
$\epsilon \subset e$, $\varphi \subset f$, and every element
of $\hat{\omega}$ is contained in exactly two of the sets 
$d, e$, and $f$, the sets $\delta, \epsilon$, and $\varphi$
form a partition of $\hat{\omega}$. Since $|\hat{\omega}| = 12$,
we obtain $ |\delta|/ 2 + |\epsilon|/ 2 + |\varphi|/ 2 = 0 \pmod{2}$.

We have $\left<d, \epsilon\right>$ =  $\left<de, \epsilon\right>$ 
= $\left<f, \epsilon\right>$ and 
$\left<d, \epsilon\right>$ = $\left<d, \delta\phi\omega\right>$
= $\left<d, \phi\right>$. The proof of the equations
$\left<d, \epsilon\right>$ = $\left<e, \phi\right>$ =
$\left<e, \delta\right>$ = $\left<f, \delta\right>$ is similar.

\proofend

Now the invariance of $(u, v, w)$ under $x_\tau$  follows
from (\ref{eqn:operations:V_T}) and Lemma~\ref{lemma:Griess:TTT}.	

\medskip

To show the invariance of $(u, v, w)$ under $y_\tau$, we will make one
more simplification.  We can find an $h \in \mathcal{P}$
with $h \cap f = \varphi$. Transformation with $y_h$ changes
$X^+_{f\cdot \varphi}$ to $X^+_{f}$, up to sign; and another
transformation with $y_d$ changes  $X^+_{f}$ to
$X^+_{f \cdot \omega}$, with $\omega = d \cap f$.
Since  $y_h y_d \in N_{xyz0}$, we may assume by symmetry that we have
$(u,v,w)$ =
$(X^+_{d\cdot \delta}, X^+_{e\cdot \epsilon}, X^+_{f\cdot \omega})$.

By (\ref{eqn:operations:V_T}), $y_\tau$ maps $X_{d,\delta}$
to a linear combination of vectors  $X_{d,\delta'}$, 
$\delta' \in A(d,\mathcal{C})$; and a similar statement holds for 
$X_{e,\epsilon}$ and $X_{f,\varphi}$. Thus Lemma~\ref{lemma:Griess:TTT}	
implies:
\begin{equation}
\label{eqn:ttt:map:Omega1}
(X^+_{d\cdot \delta}, X^+_{e\cdot \epsilon}, X^+_{f\cdot \omega})^{\tau^n} = 0
\, , \quad \mbox{ for }  \;  def \neq \pm 1 \, , \; n \in \setZ .
\end{equation} 

So we may also assume $def = 1$.
We have
\begin{align}
\noalign{$
	(X^+_{d\cdot \delta}, X^+_{e\cdot \epsilon}, X^+_{f\cdot \omega}
	)^{y_\tau}	
	$}
& \stackrel{\phantom{mm}(\ref{eqn:operations:V_T})\phantom{mm}}{=}
{\textstyle\frac{1}{512}}
\sum_{\delta' \in A(d, \mathcal{C}), \epsilon' \in A(e, \mathcal{C}),
	\varphi' \in A(f, \mathcal{C})}   \!\!\!
(-1)^{|\delta \cap \delta'| + |\epsilon \cap \epsilon'|
	+ |d \cap f \cap   \varphi'|}
\left(X^+_{d \cdot \delta'}  X^+_{e \cdot \epsilon'}
X^+_{f \cdot \varphi'} \right) \nonumber
 \\  
& \stackrel{(\rm Lemma~\ref{lemma:Griess:TTT})}{=}
{\textstyle\frac{1}{512}}
\sum_{\stackrel{\scriptstyle \delta' \in A(d, \mathcal{C}), 
 \epsilon' \in A(e, \mathcal{C}), \varphi' \in A(f, \mathcal{C}),}
		{\delta' \epsilon'  \varphi' = \omega } \!\!\!
}
(-1)^{|\delta \cap \delta'| + |\epsilon \cap \epsilon'|
  + |d \cap \varphi'|+ \langle d, \varphi' \rangle} \nonumber
\\
& \stackrel{\phantom{mmmmmmm}}{=}
{\textstyle\frac{1}{512}}
\sum_{
	\delta' \in A(d, \mathcal{C}), \epsilon' \in A(e, \mathcal{C}),
	\delta'  \epsilon' \in A(f, \mathcal{C})} \!\!\!
(-1)^{|\delta \cap \delta'| + |\epsilon \cap \epsilon'|} \;\; .
\label{eqn:Griess:TTT20} 
\end{align}
Here the condition  $\varphi' \in A(f, \mathcal{C})$,
$\delta' \epsilon' \varphi' = \omega $ in the second sum 
may be simplified to $\delta'  \epsilon' \in A(f, \mathcal{C})$,
since $\omega \in  A(f, \mathcal{C})$.

\fatline{Subcase 1:} $\delta = \epsilon \pmod{\mathcal{C}}$

Then $(X^+_{d\cdot \delta}, X^+_{e\cdot \epsilon}, X^+_{f\cdot \omega}) = 1$
by Lemma~\ref{lemma:Griess:TTT}. The Golay cocode word
$\delta = \epsilon \in A(d,\mathcal{C})  \cap  A(e,\mathcal{C})$
has a representative $\delta_0$  that is a subset of  $d \cap e$. 
Each exponent in the sum (\ref{eqn:Griess:TTT20}) is equal
to $\delta_0 \cap \delta'\epsilon'$ for some 
$\delta'\epsilon' \in  A(f, \mathcal{C})$, i.e. $\delta'\epsilon'\subset f$.
Since $d \cap e \cap f = \emptyset$, all exponents in that sum are zero.
Hence
\[
(X^+_{d\cdot \delta}, X^+_{e\cdot \epsilon}, X^+_{f\cdot \omega})
^{y_\tau} =
c \cdot 
(X^+_{d\cdot \delta}, X^+_{e\cdot \epsilon}, X^+_{f\cdot \omega}) \;
\mbox{ for a constant }  c > 0 \, .
\]
We could compute $c$ by counting the terms in the sum
(\ref{eqn:Griess:TTT20}); but we prefer a simpler argument at the
end of this subsection.

\fatline{Subcase 2:} $\delta \neq \epsilon \pmod{\mathcal{C}}$

Then $(X^+_{d\cdot \delta}, X^+_{e\cdot \epsilon}, X^+_{f\cdot \omega}) = 0$
by Lemma~~\ref{lemma:Griess:TTT}.
Let $W$ be the vector space $A(d, \mathcal{C}) \oplus A(e, \mathcal{C})$
over $\setF_2$. The mapping from $W$ to  $\setF_2$ that maps
$(\delta', \epsilon') \in W$
to the exponent in the sum  in (\ref{eqn:Griess:TTT20}) is linear.
Since that sum runs over a linear subspace of $W$ and contains
negative terms, it is zero. Hence
$(X^+_{d\cdot \delta}, X^+_{e\cdot \epsilon}, X^+_{f\cdot \omega})^{y_\tau} 
= 0$.

\proofend

Obviously, the constant $c$ in subcase 1 is the same for all
triples  $X^+_{d\cdot \delta}, X^+_{e\cdot \epsilon}, X^+_{f\cdot \varphi}$
with
$(X^+_{d\cdot \delta}, X^+_{e\cdot \epsilon}, X^+_{f\cdot \varphi}) \neq 0$.
Hence $(v, v, w)^{y_\tau} = c (u,v,w)$ holds for all $u, v, w \in V_T$.
Since $c > 0$ and $y_\tau$ is an involution acting on $V_T$, we conclude $c=1$.

\subsubsection{Cases $TTA$, $TTB$, $TTC$}

We show the invariance of $(u,v,w)$ with 
$u = X^+_{d\cdot \delta}$, $v = X^+_{e\cdot \epsilon}$,
$\delta \subset d, \epsilon \subset e, \delta, \epsilon$ even,
$w \in \{(ij)_1, X_{ij} \pm X_{\Omega\cdot ij}\}$, $i \neq j$, 
under the action of the involutions $x_\tau$ and $y_\tau$. 
By Lemma~\ref{Griess:grading:zero}
the product  $(u,v,w)$ or its images under
$\tau^{\pm1}$ are nonzero in case $d = e$ only. From
(\ref{eq:GriessAlg:rsA}) and (\ref{eq:GriessAlg:rst}) we see
that $(u,v,w) \neq 0$ holds in case $i, j \in d$ only. Then we have:
\begin{align}
(X^+_{d,\delta}, X^+_{d,\epsilon}, (ij)_1) &
\stackrel{(\ref{eq:GriessAlg:rsA})}{=}
\Delta_{\delta=\epsilon} \cdot
  (-1)^{|\delta \cap ij|} \; , \quad 
  i, j \in d,   \; \, i \neq j \; ,
 \label{eqn:Griess:TTA:ij} 
 \\
(X^+_{d,\delta}, X^+_{d,\epsilon}, 
 X_{ij} + (-1)^\nu X_{\Omega\cdot  ij}  ) & 
\stackrel{(\ref{eq:GriessAlg:rst})}{=}
 \Delta_{\delta\epsilon = ij}
  (-1)^{\nu(|\delta|/2 + |\epsilon|/2 ) } , \quad 
  i, j \in d,   \; \, i \neq j \; .
 \label{eqn:Griess:TTA:X} 
\end{align}
For equation (\ref{eqn:Griess:TTA:ij}) we remark that one of the
vectors $\pm\sqrt{8}\lambda^+_{d, \delta}$ has co-ordinate $2$
at positions $i \in \delta$, co-ordinate $-2$ at positions 
$i \in d \setminus \delta$, and zero elsewhere. 
	   
For $w = (ij)_1$ the invariance of $(u,v,w)$ under $x_\tau$ follows
from (\ref{eqn:operations:V_T}), (\ref{eqn:Griess:TTA:xtau}),
and (\ref{eqn:Griess:TTA:ij}). For $w = X_{ij} \pm X_{\Omega\cdot ij}$
the invariance of $(u,v,w)$ under $x_\tau$ follows
from (\ref{eqn:operations:V_T}), (\ref{eqn:Griess:TTA:xtau}), and
(\ref{eqn:Griess:TTA:X}). It remains to show the invariance of
$(u,v,w)$ under $y_\tau$. We have
\begin{align}
(X^+_{d,\delta}, X^+_{d,\epsilon}, (ij)_1)^{y_\tau} 
& \dstackrel{(\ref{eqn:operations:V_T}),}{(\ref{eqn:Griess:TTA:ytau})}{=}
{\textstyle\frac{1}{64}}
\sum_{\delta' \in A(d, \mathcal{C}), \epsilon' \in A(d, \mathcal{C})}
\! \!
(-1)^{|\delta \cap \delta'| + |\epsilon \cap \epsilon'|}
(X^+_{d,\delta'}, X^+_{d,\epsilon'}, X_{ij}+X_{\Omega \cdot ij})
\nonumber  \\
&  \stackrel{(\ref{eq:GriessAlg:rst})}{=} 
{\textstyle\frac{1}{64}}
\sum_{\delta' \in A(d, \mathcal{C}), \epsilon' \in A(d, \mathcal{C}),
\delta'\epsilon' = ij}
 (-1)^{|\delta \cap \delta'| + |\epsilon \cap \epsilon'|} 
\nonumber  \\
 & \phantom{m} = \phantom{m}
 {\textstyle\frac{1}{64}}
 \sum_{\delta' \in A(d, \mathcal{C})}
  (-1)^{|\delta \epsilon \cap \delta'| + |\epsilon \cap ij|} 
\nonumber  \\
 & \phantom{m} = \phantom{m}
\Delta_{\delta=\epsilon} \cdot  (-1)^{|\delta \cap ij|}  \; .
\label{eqn:Griess:TTA:calc}
\end{align}
This proves the invariance of $(u,v,w)$ under $y_\tau$ for all basis
vectors $w$ of $V_A$. Since images of these basis vectors under
the involution  $y_\tau$ form a basis of  $V_B$, the invariance of
of $(u,v,w)$ under $y_\tau$ has also been shown for the basis vectors
$w =  X_{ij}+X_{\Omega \cdot ij}$ of $V_B$.
\Skip{                                           
Todo: Show invariance under  $y_\tau$ for
 $w =   X^+_{d,\epsilon'}, X_{ij}-X_{\Omega \cdot ij}$ !!!
}

Using $|\alpha|/2 + |\beta|/2 + |\alpha\beta|/2 = |\alpha \cap \beta|
\pmod{2}$, $\alpha, \beta \in A(d, \mathcal{C})$, several times,
 we obtain:
\begin{align}
(X^+_{d,\delta}, X^+_{d,\epsilon},  X_{ij}\!-\!X_{\Omega \cdot ij})
^{y_\tau}
& \dstackrel{(\ref{eqn:operations:V_T}),}{(\ref{eqn:Griess:TTA:ytau})}{=}
\! {\textstyle\frac{1}{64}} \! \! 
\sum_{\delta' \in A(d, \mathcal{C}), \epsilon' \in A(d, \mathcal{C})}
\! \! \!
(-1)^{|\delta \cap \delta'| + |\epsilon \cap \epsilon'|}
(X^+_{d,\delta'}, X^+_{d,\epsilon'}, X_{ij}\!-\!X_{\Omega \cdot ij})
\nonumber  \\
&  \stackrel{(\ref{eq:GriessAlg:rst})}{=} 
 {\textstyle\frac{1}{64}}
\sum_{\delta' \in A(d, \mathcal{C}), \epsilon' \in A(d, \mathcal{C}),
	\delta'\epsilon' = ij}
(-1)^{|\delta \cap \delta'| + |\epsilon \cap \epsilon'|
  + |\delta'|/2 + |\epsilon'|/2 }	
\nonumber \\
 & \phantom{m} = \phantom{m}
 {\textstyle\frac{1}{64}}
\sum_{\delta' \in A(d, \mathcal{C})}
(-1)^{|\delta \cap \delta'| + |\epsilon \cap \delta' ij|
	+ |\delta'|/2 + |\delta' ij|/2 }	
\nonumber \\
 & \phantom{m} = \phantom{m}
 {\textstyle\frac{1}{64}}  \cdot
(-1)^{ |\epsilon ij|/2 + |\epsilon|/2} \cdot
  \sum_{\delta' \in A(d, \mathcal{C})}
(-1)^{|\delta \epsilon ij \cap \delta'|} 	 
 \label{eqn:Griess:TTA:calc2}
\nonumber  \\
 & \phantom{m} = \phantom{m}
\Delta_{\delta\epsilon = ij}  (-1)^{|\delta|/2 + |\epsilon|/2}  \; .
\end{align}
Here the last sum in (\ref{eqn:Griess:TTA:calc2}) is 64 if
$\delta\epsilon = ij$ and 0 otherwise. So 
$(X^+_{d,\delta}, X^+_{d,\epsilon},  X_{ij}\!-\!X_{\Omega \cdot ij})$
is invariant under the operation of  $y_\tau$ by
 (\ref{eqn:Griess:TTA:X}) and (\ref{eqn:Griess:TTA:calc2}).

\subsubsection{Case $TTD$}

As in the previous cases we see that 
$(X^+_{d,\delta}, X^+_{e,\epsilon}, (ii)_1)^g$, $g = 1, \tau^{\pm 1}$
is nonzero in case $d = e$ only. By definition both $x_\tau$ and
$y_\tau$ fix $(ii_1)$.  We have
\begin{align*}
(X^+_{d,\delta}, X^+_{d,\epsilon}, (ii_1)) &
\stackrel{(\ref{eq:GriessAlg:rsA})}{=}
{\textstyle \frac{1}{2}} \Delta_{\delta = \epsilon} \; .
\end{align*}
This is invariant under $x_\tau$ by  (\ref{eqn:operations:V_T}). 
A calculation similar to that in (\ref{eqn:Griess:TTA:calc2}) yields:
\begin{align*}
(X^+_{d,\delta}, X^+_{d,\epsilon}, (ii_1))^{y_\tau} & =
{\textstyle\frac{1}{64}} \cdot {\textstyle\frac{1}{2}} \cdot
\sum_{\delta' \in A(d, \mathcal{C}), \epsilon' \in A(e, \mathcal{C}),
	\delta' = \epsilon'}
(-1)^{|\delta \cap  \delta'| + |\epsilon \cap \epsilon'|} \\
&= {\textstyle \frac{1}{2}} \Delta_{\delta = \epsilon} \; .
\end{align*}

\subsection{Cases involving  subspaces $V_A, V_B, V_C, V_D$ only }

\stepcounter{Theorem}

Here it
suffices to consider the restriction of the Griess algebra form
to $V_{ABCD}$, referred to as the {\em form} in the remainder of 
this subsection. We will show the invariance the form under
$x_\tau$ and $y_\tau$, with $x_\tau, y_\tau$ as in the previous
subsection. $x_\tau$ acts on  $V_{ABCD}$ is the same way as the
element  $x_i$, $i \in \tilde{\Omega}$ of $G_{x0}$, proving the
invariance of the form under $x_\tau$.

To reduce the number of cases to be considered for the action
of $y_\tau$, we will define a finer grading on $V_{ABCD}$ than
the grading given in Section~\ref{subsection:act:Griess}. The
group $R$ defined in that section acts trivially on  $V_{ABCD}$.
Let $R'$ be the group generated by
$x_d, y_d, z_d, d \in \mathcal{P}$, and $R$. By
Theorem~\ref{Thm:relations:N} the commutator group of $R'$ 
is a subgroup of $R$; so the Abelian group $R_1 = R' / R$ has
a well-defined action on  $V_{ABCD}$.  $R_1$ is isomorphic
to the additive group of the vector space 
$\setF_2^2 \otimes (\mathcal{C} / \{1, \tilde{\Omega}\})$ over
$\setF_2$, with the elements of the first factor denoted by
$\{0,x,y,z\}$. Similar to the construction in
Section~\ref{subsection:act:Griess}, the algebra on
$V_{ABCD}$ has a grading with grades labelled by the 
Pontryagin dual $\hat{R_1}$ of $R_1$. The  group $\hat{R_1}$ is
isomorphic to the additive group of the dual space
$\setF_2^2 \otimes \{\delta \in \mathcal{C^*} \mid \delta \,
\mbox{even} \}$ of the space mentioned above. Here we identify
$\setF_2^2$ with its dual space as in 
Section~\ref{subsection:act:Griess}.
E.g. a vector $v \in V_{ABCD}$ has grade $x \otimes \delta$ if
$x_d$ fixes $v$ and $y_e, z_e$ multiply $v$ with
$(-1)^{\langle e, \delta \rangle}$. From 
Table~\ref{table:action:196884x} we see that our choosen basis
vectors of $V_{ABCD}$ have the following grades: 
\[
\begin{array}{c|c|c|c|c}
\mbox{basis vector} & (ii)_1 & (ij)_1  & X_{ij} - X_{\Omega \cdot ij}
&  X_{ij} + X_{\Omega \cdot ij} \\
\hline
\mbox{in space} & V_D & V_A & V_B & V_C \\
\hline
\mbox{grade} &  0 \otimes 0 = 0 &  x \otimes ij & y \otimes ij
  & z \otimes ij
\end{array} \; \;  .
\]
The same argument as in the proof of Lemma~\ref{Griess:grading:zero}	
shows the the grades of three basis vectors $u,v,w$ of  $V_{ABCD}$
must sum up to 0 in case $(u,v,w) \neq 0$.
By (\ref{eqn:Griess:TTA:ytau}) the involution $y_\tau$ exchanges
$V_A$ with $V_C$, and fixes $V_B$ and $V_D$ pointwise. So we are left
with the cases $(u,v,w) \neq 0$, where $V_A$ or $V_C$ is involved,
and the grades of $u, v$, and $w$ sum up to 0:
\begin{align}
((ij)_1, X_{ij} + X_{\Omega \cdot ij}, X_{ij} - X_{\Omega \cdot ij}) 
& \stackrel{(\ref{eq:GriessAlg:rsA})}{=} -8 , \quad i \neq j \;, 
\nonumber \\
((ii)_1, (jk)_1, (jk)_1) 
& \stackrel{(\ref{eq:GriessAlg:AAA})}{=} 4(\Delta_{i=k} + \Delta_{j=k}) 
, \quad j \neq k \; , \label{Griess_ABCD:1} \\
((ii)_1, X_{jk} + X_{\Omega \cdot jk}, X_{jk} + X_{\Omega \cdot jk}) 
& \stackrel{(\ref{eq:GriessAlg:rsA})}{=} 4(\Delta_{i=k} + \Delta_{j=k})
, \quad j \neq k \; , \label{Griess_ABCD:2} \\
((ij)_1, (ik)_1, (jk)_1) 
& \stackrel{(\ref{eq:GriessAlg:AAA})}{=}
4,  \quad |\{i,j,k\}| = 3 \; , \label{Griess_ABCD:3}  \\
( X_{ij} + X_{\Omega \cdot ij}, X_{ik} + X_{\Omega \cdot ik},
 X_{jk} + X_{\Omega \cdot jk}) 
& \stackrel{(\ref{eq:GriessAlg:rsA})}{=}
4,
\quad  |\{i,j,k\}| = 3 \; . \label{Griess_ABCD:4} 
\end{align}
Note that one of the
vectors $\pm\sqrt{8}\lambda_{\Omega^\nu,ij}$ has co-ordinate $-4$
at position $i$, co-ordinate $4 \cdot (-1)^\nu$ at position $j$,
and zero elsewhere. The invariance of the form under the action of
$y_\tau$ follows from the equality of (\ref{Griess_ABCD:1}) and 
(\ref{Griess_ABCD:2}), and the equality of (\ref{Griess_ABCD:3})
and  (\ref{Griess_ABCD:4}).

\proofend

%\subfile{Formulas}

% !TeX spellcheck = en_GB

\clearpage

%%%%%%%%%%%%%%%%%%%%%%%%%%%%%%%%%%%%%%%%%%%%%%%%%%%%%%%%%%%%%%%%%%%%%%%
\section*{Notation}
\label{setion:Notation}
%%%%%%%%%%%%%%%%%%%%%%%%%%%%%%%%%%%%%%%%%%%%%%%%%%%%%%%%%%%%%%%%%%%%%%%

\begin{tabular}{cll}
Symbol & Description 
\hphantom{012345678901234567890123456789012345678901234567}
       & \hspace{-2.25em} Section \vspace{1ex} \\
\multicolumn{2}{l} {$a,b,c,d,e,f,h \quad$
 Elements of the Parker loop $\mathcal{P}$ or of
 the Golay code $\mathcal{C}$}
       &  \ref{subsection:Parker:loop}  \\
$A(d,e,f)$ &
Associator map of the elements $d,e,f$ of
the Parker loop $\mathcal{P}$
       & \ref{subsection:Parker:loop}  \\
$\AutStP$ &
    The group of standard automorphisms of the Parker
    loop $\mathcal{P}$    
       & \ref{section:auto:Parker}  \\ 
$C(d,e)$ &
Commutator map of the elements $d$ and $e$ of 
the Parker loop $\mathcal{P}$
       & \ref{subsection:Parker:loop}  \\
$\mathcal{C}$, $\mathcal{C}^*$ &
$\mathcal{C}$ is the 12-dimensional Golay code in $\Field_2^{24}$,
 $\mathcal{C}^*$  its cocode  $\Field_2^{24}/\mathcal{C}$
       & \ref{subsection:Golay}  \\
$\delta, \epsilon, \varphi, \eta$ &
    Elements of the Golay cocode  $\mathcal{C}^*$
       & \ref{subsection:Parker:loop}  \\    
\Skip{       
$g$  &
   A function  $\mathcal{G^*} \cup \mathcal{C} \rightarrow
   \mathcal{P}_\mathcal{G}^+$  with $g(\gamma_i) = g_i$,
   $g(d) = g(\gamma(d))$,  $d \in    \mathcal{C}$.
      &   \ref{subsection:cocycle}  \\   
} % \Skip       
$g_0,\ldots,g_5$ &
Standard  basis of the ''grey'' subspace $\mathcal{G}$
of the Golay code $\mathcal{C}$ ; \\
      &
     $g_i, i=0,\ldots,5,$ is also considered as the
     element $(g_i,0)$ of  $\mathcal{P}_\mathcal{G}$. 
       & \ref{subsection:gray:col}  \\
$\mathcal{G}$, $\mathcal{G}^*$  &
The  subspaces of the ''grey'' elements of 
 $\mathcal{C}$ and $\mathcal{C}^*$, respectively 
       & \ref{subsection:gray:col}  \\
$G_{x0} $ &
  A maximal subgroup of $\setM$ of structure $2^{1+24}_+.\mbox{Co}_1$ 
       & \ref{section:N0}  \\      
$\gamma$ &
A specific function $\mathcal{C} \rightarrow \mathcal{G}^*$,  
satisfying $\gamma(g_i) = \gamma_i$, $i=0,\ldots,5$
       &  \ref{subsection:cocycle}  \\
$\gamma_0,\ldots,\gamma_5$ &      
Standard  basis of the ''grey'' subspace $\mathcal{G}^*$
of the Golay cocode $\mathcal{C}^*$             
       & \ref{subsection:gray:col}  \\              
$\mathcal{H}$, $\mathcal{H}^*$   &
The  subspace of the ''coloured'' elements of 
$\mathcal{C}$ and $\mathcal{C}^*$, respectively 
       & \ref{subsection:gray:col}  \\
$i,j,k$ &
Elements of $\tilde{\Omega}$,
also considered  as elements of  $\mathcal{C}^*$  of weight 1
       &  \ref{subsection:Parker:loop}  \\       
$ij$ &
Shorthand for $i \cup j$,  considered  as an element of 
$\mathcal{C}^*$ of weight 2
       &  \ref{subsection:Parker:loop}  \\  
$\Lambda$ &
The 24-dimensional Leech lattice. 
       & \ref{subsection:Qx:Leech} \\
$\setM$ &
The monster group, i.e. the largest sporadic simple group   
       &  \ref{section:Introduction}  \\ 
$M_{24}$ &
  Mathieu group,  acts on $\tilde{\Omega}$ as
  the automorphism group of $\mathcal{C}$.
   & \ref{subsection:Golay} \\      
MOG &     
Miracle Octad Generator; a tool for calculations in 
 $\mathcal{C}  \subset \setF_2^{24}$.   
    & \ref{subsection:Golay}   \\                 
$N$ &
A fourfold cover of the maximal subgroup $N_0$ of  $\setM$ 
       &  \ref{section:N0}  \\  
$N_0$ &
A maximal subgroup of $\setM$ of structure
  $2^{2+11+22}.(M_{24} \times S_3)$
       &  \ref{section:N0}  \\     
$N_{x0}$ &       
  Subgroup of structure of $\setM$ 
    with $G_{x0} \cap N_0 = N_{x0}$
  & \ref{section:N0}  \\                
${\Omega} $   &
The  element $(\tilde{\Omega},0)$ of the Parker loop $\mathcal{P}$
       & \ref{subsection:Golay}    \\    
$\tilde{\Omega} $
 &  A set of size 24 used for labeling the basis vectors of
   $\Field_2^{24}$,  &  \\
 & its power set $2^{\tilde{\Omega}} $ is identified
   with $\Field_2^{24}$, 
     and we have  $\mathcal{C} \subset \Field_2^{24}$
       & \ref{subsection:Golay}    \\    
$\omega$ &
A specific ''grey'' element in the subset  $\mathcal{G}^*$ of
the  cocode  $\mathcal{C}^*$  
& \ref{subsection:gray:col}  \\
$P(d)$ &
The  squaring map 
in $\mathcal{P}$, with $d^2 = (0,P(d))$ 
for $d \in \mathcal{P}$
       & \ref{subsection:Parker:loop}  \\
$\mathcal{P}$ &
The Parker loop, any $d \in \mathcal{P}$ has the form
$(\tilde{d}, \mu)$, $\tilde{d} \in \mathcal{C}, \mu \in \Field_2$
       & \ref{subsection:Parker:loop}  \\
$\mathcal{P}_\mathcal{G}, \mathcal{P}_\mathcal{H}$ &
Subsets of $\mathcal{P}$:
  origins of $\mathcal{G}$ and $\mathcal{H}$ of the mapping
$\tilde{\hphantom{.}}$ : 
 $\mathcal{P} \rightarrow \mathcal{C}$.
    &  \ref{subsection:cocycle}  \\
\Skip{    
$\mathcal{P}_\mathcal{G}^+$ &
The subset $\{ (d,0) | d \in \mathcal{G}\}$ 
 of $\mathcal{P}_\mathcal{G}$.
     &  \ref{subsection:cocycle}  \\
}     
$\pi, \pi', \pi''$ &
Standard automorphisms of the Parker Loop  $\mathcal{P}$ 
    in $\AutStP$ 
       & \ref{section:auto:Parker}  \\   
$Q_{x}$ &
Subgroup  of structure   $2_+^{1+24}.\mbox{Co}_1$ of the group $G_{x0}$ 
      &   \ref{section:N0} \\  
$\mbox{sign}(d)$ &      
       ''Sign'' of an element $d$ of the Parker loop  $\mathcal{P}$.
     & \ref{subsection:Parker:loop} \\
$\theta$ &
Cocycle  of Parker loop $\mathcal{P}$, with
  $(\tilde{d},0) \cdot (\tilde{e},0)
     = (\tilde{d}+\tilde{e},\theta(\tilde{d},\tilde{e})) $
       &  \ref{subsection:any:cocycle}   \\ 
$w(d)$  &
Weight of vector $d \in \mathcal{G}$ with respect to the
basis $g_0,\ldots,g_5$               
       & \ref{subsection:gray:col}    \\    
$w(\delta)$  &
Weight of vector $\delta \in \mathcal{G}^*$ with respect to the
basis $\gamma_0,\ldots,\gamma_5$               
      & \ref{subsection:gray:col}    \\
$w_2(d), w_2(\delta)$ &
   Equal to $\binom{w(d)}{2}$,  $\binom{w(\delta)}{2}$ 
     modulo 2, for $d \in \mathcal{C}, \delta d \in \mathcal{C}^*$
      & \ref{subsection:cocycle}       \\  
$Z(G)$ &  
   The center of a group or a loop $G$ 
    &     \ref{section:auto:Parker}  \\
\Skip{
\multicolumn{2}{l} {      
$\omega_0,\ldots,\omega_5, \omega_\infty \quad$ 
Specific representatives of the Golay cocode word  $\omega$
in $\Field_2^{24}$ }
& \ref{subsection:gray:col}  \\  
}  %% \Skip     
$\tilde{d}$ &
Image of $d \in \mathcal{P}$ in $\mathcal{C}$ under
the natural homomorphism $\mathcal{P} \rightarrow \mathcal{C}$
       & \ref{subsection:Parker:loop}  \\  
$\tilde{x}_d$ &
    The element $x_{-1}^{\mbox{\scriptsize sign}(d)} x_{d} x_{\theta(d)} $
    of  the group $Q_{x0}$, for
     $d \in \mathcal{P}$ or  $d \in \mathcal{C}$ 
    & \ref{section:N0} \\  
$\bar{d}$ &
       Inverse of $d$ in the Parker loop  $\mathcal{P}$:
          $\bar{d} = d^{-1}$  
       & \ref{subsection:Parker:loop}  \\
$|d|$, $|\delta|$ &
Weight of code word $d \in \mathcal{C}$, min. weight
of cocode word  $\delta\in \mathcal{C}^*$
      & \ref{subsection:Golay}    \\
$\left<.,.\right>$ & 
   The scalar product, e.g. on $\mathcal{C} \times \mathcal{C^*}$
     &  \ref{subsection:Golay}     \\
$\smallbiscalar{d,e}$ &    
   equal to $\left<d, \gamma(e)\right>$ and to
   $w_2(de) + w_2(d) + w_2(e)$ for  $d,e \in \mathcal{G}$
      & \ref{subsection:cocycle}       \\ 
'$ * $' & The product in the Griess algebra
   & \ref{subsection:Griess:algebra}
                
\end{tabular}

\pagebreak

\end{document}